\documentclass[ctagsplt]{article}
\usepackage{amssymb}

\usepackage{graphicx}
\usepackage{amsmath}

\newtheorem{theorem}{Theorem}

\newtheorem{corollary}[theorem]{Corollary}
\newtheorem{definition}[theorem]{Definition}
\newtheorem{example}[theorem]{Example}

\newtheorem{proposition}[theorem]{Proposition}
\newtheorem{remark}[theorem]{Remark}

\input{tcilatex}

\begin{document}

\title{$ADE$-bundle over rational surfaces, \\
configuration of lines and rulings}
\author{Naichung Conan Leung\thanks{%
This paper is partially by a NSF grant, DMS-9803616.} \\
School of Mathematics, \\
University of Minnesota, \\
Minneapolis, MN 55455, U.S.A..\\
and\\
IHES, \\
91440 Bures-Sur-Yvette, France.}
\maketitle

\begin{abstract}
To each del Pezzo surface $X_{n}$ (resp. ruled surface, ruled surface with a 
section), we describe a natural Lie algebra bundle of type $E_{n}$ (resp. $%
D_{n-1}$, $A_{n-2}$) over it . 

Using lines and rulings on any such surface, we describe various representation
bundles corresponding to fundamental representations of the corresponding
Lie algebra.

When we specify a geometric structure on the surface to reduce the Lie
algebra to a smaller one, then the classical geometry of the configuration
of lines and rulings is encoded beautifully by the branching rules in Lie
theory. We discuss this relationship in details.

When we degenerate the surface to a non-normal del Pezzo surface, we
discover that the configurations of lines and rulings are also governed by
certain branching rules. However the degeneration theory of the bundle is
not fully understood yet. 

\newpage 
\end{abstract}

\tableofcontents

\pagebreak

\section{Introduction and results}

\bigskip

\footnote{%
The content of this paper overlaps partially with a recent preprint of
Friedman and Morgan 'Exceptional groups and del Pezzo surfaces'
(math.AG/0009155), they construct $E_{n}$ bundles over del Pezzo surfaces,
possibly with rational double point singularities. They also discuss
reducing the structure groups. We pay more emphasis on the interplays
between representation theory of exceptional Lie algebra, reduction of these
bundles and classical geometry of configurations of lines and rulings on
rational surfaces. We also learn from their preprint that our bundles are
really conformal bundles as in their paper.}\textbf{Lines on del Pezzo
surfaces}

Studying of lines on a cubic surface is a fascinating subject and it has a
very long history \cite{Segre}. It is a classical fact that there are $27$
lines on a cubic surface. Moreover the configuration of these lines
possesses a high degree of symmetry which is closely related to the Weyl
group of $E_{6}$. Similar results hold true for other del Pezzo surfaces
(see e.g. \cite{Manin}). Recall that a del Pezzo surface is a smooth surface
with an ample anti-canonical divisor. Such a surface is either $\mathbb{P}%
^{1}\times \mathbb{P}^{1}$ or a blowup of $\mathbb{P}^{2}$ at $n$ generic
points with $n\leq 8$.We call the latter surface $X_{n}$. When $n\leq 6$ the
anti-canonical linear system embeds $X_{n}$ in $\mathbb{P}^{d}$ as a degree $%
d$ surface and $d=9-n$. A line on $X_{n}$ is then equivalent to a divisor $l$
with $l^{2}=l\cdot K=-1$. We will continue to call such divisor a line even $%
n=7$ or $8$. Analogous to the cubic surface case, lines on $X_{n}$ possess a
high degree of symmetry which is closely related to the Weyl group of $E_{n}$%
. 
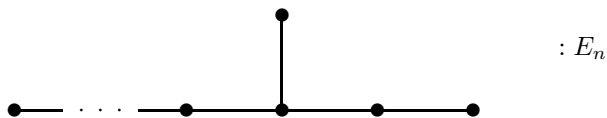
\begin{equation*}
\end{equation*}
\begin{picture}(3.4,1)(1.4,.4)
\thicklines
\put(4.3,1){\circle*{.075}}
\put(3.8,1){\circle*{.075}}
\put(3.8,1){\line(1,0){.5}}
\put(3.3,1){\circle*{.075}}
\put(3.3,1){\line(0,1){.5}}
\put(3.3,1.5){\circle*{.075}}
\put(3.3,1){\line(1,0){.5}}
\put(2.8,1){\circle*{.075}}
\put(2.8,1){\line(1,0){.5}}
\put(2.55,1){\line(1,0){.25}}
\put(2.45,1){\circle*{.02}}
\put(2.35,1){\circle*{.02}}
\put(2.25,1){\circle*{.02}}
\put(1.9,1){\circle*{.075}}
\put(1.9,1){\line(1,0){.25}}
\put(4.75,1.2){\makebox(.25,.25){$: E _ n$}}
\end{picture}

When $n\leq 5$, $E_{n}$ coincides with classical Lie algebras as follows, $%
E_{5}=\mathbf{so}\left( 10\right) $, $E_{4}=\mathbf{sl}\left( 5\right) $, $%
E_{3}=\mathbf{sl}\left( 3\right) \times \mathbf{sl}\left( 2\right) $, $E_{2}=%
\mathbf{sl}\left( 2\right) \times \mathbf{u}\left( 1\right) $ and $E_{1}=%
\mathbf{u}\left( 1\right) $ (or $\mathbf{sl}\left( 2\right) $ which we
denote $\overline{E}_{1}$).

If we look at the fundamental representation of the Lie algebra $E_{n}$
corresponding the left end node of its Dynkin diagram and call it $\mathbf{L}%
_{n}$. Then we see an interesting relationship between $\dim \mathbf{L}_{n}$
and the number of lines on $X_{n}:$%
\begin{equation*}
\begin{tabular}{|c|l|l|l|l|l|l|l|l|}
\hline
$n$ & 1 & 2 & 3 & 4 & 5 & 6 & 7 & 8 \\ \hline
$\dim \mathbf{L}_{n}$ & 1 & 3 & 6 & 10 & 16 & 27 & 56 & 248 \\ \hline
$\text{\# lines on }X_{n}$ & 1 & 3 & 6 & 10 & 16 & 27 & 56 & 240 \\ \hline
\end{tabular}
\end{equation*}
Namely $\dim \mathbf{L}_{n}\,$equals the number of lines on $X_{n}$ except
when $n$ $=8$. In that exceptional case $\mathbf{L}_{8}$ coincides with the
adjoint representation of $E_{8}$ and $\dim \mathbf{L}_{8}$ equals the
number of lines on $X_{8}$ plus the rank of $E_{8}=8$. Motivated from this
we define a holomorphic bundle $\mathcal{L}_{n}$ on $X_{n}$ using its lines
as follows: 
\begin{eqnarray*}
\mathcal{L}_{n} &=&\bigoplus_{l\text{:line on }X_{n}}O\left( l\right)
\,\,\,\,\,\,\,\,\,\text{when }n\leq 7, \\
\mathcal{L}_{8} &=&\bigoplus_{l\text{:line on }X_{8}}O\left( l\right)
+O\left( -K\right) ^{\oplus 8}.
\end{eqnarray*}
We will show that $E_{n}$ is the structure group for $\mathcal{L}_{n}$.
Moreover different specializations of $X_{n}$ would reduce its structure
group to various subgroups of $E_{n}$. On these specializations of $X_{n}$,
special configurations of lines would translate into branching rules for the
representation $\mathbf{L}_{n}$ under various subgroups of $E_{n}$.

\bigskip

\bigskip

\textbf{Rulings on del Pezzo surfaces}

Besides lines, another geometric structure on $X_{n}$ is a ruling on $X_{n}$%
, which is a fibration of $X_{n}$ over $\mathbb{P}^{1}$ whose generic fiber
is a smooth rational curve. Any fiber $R$ of a ruling satisfies $R^{2}=0$
and $R\cdot K=-2$. Conversely given any such divisor, it is the fiber class
of a unique ruling on $X_{n}$. For example there are two rulings on $\mathbb{%
P}^{1}\times \mathbb{P}^{1}$ and one ruling on $X_{1}$, the blowup of $%
\mathbb{P}^{2}$ at a point. As we will see, these two cases correspond to
standard representations of $\overline{E}_{1}=\mathbf{sl}\left( 2\right) $
and $E_{1}=\mathbf{u}\left( 1\right) $ respectively.

If we look at the fundamental representation of the Lie algebra $E_{n}$
corresponding to the right end node of its Dynkin diagram and call it $%
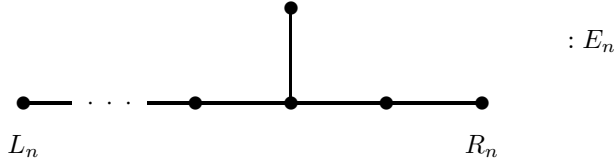
\mathbf{R}_{n}$. 
\begin{equation*}
\end{equation*}
\begin{picture}(3.4,1)(1.4,.4)
\thicklines
\put(4.3,1){\circle*{.075}}
\put(3.8,1){\circle*{.075}}
\put(3.8,1){\line(1,0){.5}}
\put(3.3,1){\circle*{.075}}
\put(3.3,1){\line(0,1){.5}}
\put(3.3,1.5){\circle*{.075}}
\put(3.3,1){\line(1,0){.5}}
\put(2.8,1){\circle*{.075}}
\put(2.8,1){\line(1,0){.5}}
\put(2.55,1){\line(1,0){.25}}
\put(2.45,1){\circle*{.02}}
\put(2.35,1){\circle*{.02}}
\put(2.25,1){\circle*{.02}}
\put(1.9,1){\circle*{.075}}
\put(1.9,1){\line(1,0){.25}}
\put(1.775,.65){\makebox(.25,.25){$L _ n$}}
\put(4.175,.65){\makebox(.25,.25){$R _ n$}}
\put(4.75,1.2){\makebox(.25,.25){$: E _ n$}}
\end{picture}

Then we see an interesting relationship between $\dim \mathbf{R}_{n}$ and
the number of rulings on $X_{n}:$%
\begin{equation*}
\begin{tabular}{|c|c|c|c|c|c|c|c|}
\hline
$n$ & $1$ & $2$ & $3$ & $4$ & $5$ & $6$ & $7$ \\ \hline
$\dim \mathbf{R}_{n}$ & $1$ & $2$ & $3$ & $5$ & $10$ & $27$ & $133$ \\ \hline
$\text{\# rulings on }X_{n}$ & $1$ & $2$ & $3$ & $5$ & $10$ & $27$ & $126$
\\ \hline
\end{tabular}
\end{equation*}
Namely $\dim \mathbf{R}_{n}\,$equals the number of rulings on $X_{n}$ for $%
n\leq 6$. When $n=7,$ $\mathbf{R}_{7}$ coincides with the adjoint
representation of $E_{7}$ and $\dim \mathbf{R}_{7}$ equals the number of
rulings on $X_{7}$ plus the rank of $E_{7}$, which is 7. Motivated from this
we define a holomorphic bundle $\mathcal{R}_{n}$ on $X_{n}$ using its
rulings as follows: 
\begin{eqnarray*}
\mathcal{R}_{n} &=&\bigoplus_{R\text{:ruling on }X_{n}}O\left( R\right)
\,\,\,\,\,\,\,\,\,\text{when }n\leq 6, \\
\mathcal{R}_{7} &=&\bigoplus_{R\text{:ruling on }X_{7}}O\left( R\right)
+O\left( -K\right) ^{\oplus 7}.
\end{eqnarray*}
The construction of $\mathcal{R}_{8}$ includes more than rulings on $X_{8}$
and we shall not discuss it here. Again we will see that the structure group
for $\mathcal{R}_{n}$ is given by $E_{n}$. Different specializations of $%
X_{n}$ would reduce its structure group to various subgroups of $E_{n}$. On
these specializations of $X_{n},$ special configurations of rulings would
translate into branching rules for the representation $\mathbf{R}_{n}$ under
various subgroups of $E_{n}$.

\bigskip

\textbf{Cubic surfaces}

For example in the cubic surface $X_{6}$ case, if we fix any line $L$ on it,
then there are 10 other lines intersecting $L$ and 16 other lines disjoint
from $L$. In terms of $\mathcal{L}_{6}$, the choice of a line reduces the
structure group from $E_{6}$ to $E_{5}$ and we have a decomposition of
representation bundles. 
\begin{equation*}
\mathcal{L}_{6}=\pi ^{\ast }\mathcal{L}_{5}+\pi ^{\ast }\mathcal{R}%
_{5}\otimes O\left( -L\right) +O\left( L\right) ,
\end{equation*}
where $\pi :X_{6}\rightarrow X_{5}$ is the blowdown of $L$. This
decomposition corresponds to the branching rule $\mathbf{L}_{6}|_{E_{5}}=%
\mathbf{L}_{5}+\mathbf{R}_{5}+\mathbf{1}$, where $\mathbf{1}$ denotes the
trivial representation. The component $\pi ^{\ast }\mathcal{L}_{5}$
corresponds to pullback of lines from $X_{5}$ and hence disjoint from $L.$
Each direct summand of the component of $\pi ^{\ast }\mathcal{R}_{5}\otimes
O\left( -L\right) $ corresponds to the pullback of a particular member of a
ruling which passes through the blowup point. Therefore such lines on $X_{6}$
would intersect $L$ at one point. Moreover these 10 lines intersecting $L$
divides into 5 pairs, each pair forms a triangle with $L$. This corresponds
to the fact the $E_{5}$ equals $\mathbf{so}\left( 10\right) $ and $\mathbf{R}%
_{5}$ is the standard representation of $\mathbf{so}\left( 10\right) $.

In general, if we fix a ruling on $X_{n}$ then the structure group of $%
\mathcal{L}_{n}$ or $\mathcal{R}_{n}$ would reduce to $D_{n-1}=\mathbf{so}%
\left( 2n-2\right) $. If we further choose a line section of this ruling
then the structure group would further reduce to $A_{n-2}=\mathbf{sl}\left(
n-1\right) $. various geometric properties of lines and rulings on $X_{n}$
related to corresponding branching rules will be discussed in section three
and four. One should notice that we can discuss $D_{n-1}$-bundle (resp. $%
A_{n-2}$-bundle) on blowup of $\mathbb{P}^{2}$ at $n$ generic points
together with a ruling (resp. a ruling with a section) without any
restriction on $n$.

For rulings on a cubic surface, they are in one-to-one correspondence with
lines. Namely, given any ruling on $X_{6}$, there is a unique line which is
a 2-section of the ruling. This is because $\mathcal{R}_{6}=\mathcal{L}%
_{6}^{\ast }\otimes O\left( -K\right) $ which reflects the isomorphism $%
\mathbf{R}_{6}=\mathbf{L}_{6}^{\ast }$ between $E_{6}$-representations.

When we degenerate the cubic surface to a union of three planes, then the
structure Lie algebra would reduce from $E_{6}$ to $\mathbf{sl}\left(
3\right) \times \mathbf{sl}\left( 3\right) \times \mathbf{sl}\left( 3\right) 
$. When we degenerate it to a union of a hyperplane and a smooth quadratic
surface, then the structure Lie algebra would reduce from $E_{6}$ to $%
\mathbf{sl}\left( 2\right) \times \mathbf{sl}\left( 6\right) $. Details of
these cases will be given in section six. These two are examples of
degenerations into nonnormal del Pezzo surfaces \cite{Reid1}, the
discussions on exactly how $\mathcal{L}_{n}$ and $\mathcal{R}_{n}$
degenerate in such situations is not completed here. We hope to come back to
this problem in the future.

In general $\mathcal{L}_{n}$ and $\mathcal{R}_{n}$ are representation
bundles of a $E_{n}$-bundle $\mathcal{E}_{n}$ over $X_{n}$ which we now
describe. 
\begin{equation*}
\begin{tabular}{|ll|}
\hline
(i) $\mathcal{E}_{4}$ & $\,$is the automorphism bundle of $\mathcal{R}_{4}$
preserving \\ 
&  \\ 
& $\Lambda ^{5}\mathcal{R}_{4}\cong O\left( -2K\right) .$ \\ \hline
(ii) $\mathcal{E}_{5}$ & $\,$is the automorphism bundle of $\mathcal{R}_{5}$
preserving \\ 
&  \\ 
& $\,q_{5}:\mathcal{R}_{5}\otimes \mathcal{R}_{5}\rightarrow O\left(
-K\right) .$ \\ \hline
(iii) $\mathcal{E}_{6}$ & $\,$is the automorphism bundle of $\mathcal{R}_{6}$
and $\mathcal{L}_{6}$ preserving \\ 
&  \\ 
& $\,c_{6}:\mathcal{L}_{6}\otimes \mathcal{L}_{6}\rightarrow \mathcal{R}_{6},
$ and \\ 
& $c_{6}^{\ast }:\mathcal{R}_{6}\otimes \mathcal{R}_{6}\rightarrow \mathcal{L%
}_{6}\otimes O\left( -K\right) .$ \\ \hline
(iv) $\mathcal{E}_{7}$ & is the automorphism bundle of $\mathcal{L}_{7}$
preserving \\ 
&  \\ 
& $\,f_{7}:\mathcal{L}_{7}\otimes \mathcal{L}_{7}\otimes \mathcal{L}%
_{7}\otimes \mathcal{L}_{7}\rightarrow O\left( -2K\right) .$ \\ \hline
(v) $\mathcal{E}_{8}$ & $\,$is the automorphism bundle of $\mathcal{L}_{8}$
preserving \\ 
&  \\ 
& $\mathcal{L}_{8}\wedge \mathcal{L}_{8}\rightarrow \mathcal{L}_{8}\otimes
O\left( -K\right) .$ \\ \hline
\end{tabular}
\end{equation*}

\textbf{Physics motivations}

These $E_{n}$-bundles are related to F-theory in physics. They are also
defined and studied in a very recent paper by Friedman and Morgan \cite{FM}.
If $\Sigma $ is an anti-canonical curve in $X_{n}$, then $\Sigma $ is of
genus one. Restricting $\mathcal{E}_{n}$ to $\Sigma $ for various complex
structures on $X_{n}$ and various embeddings of $\Sigma $ into $X_{n}$ gives
the moduli space of flat $E_{n}$-bundles over the elliptic curve $\Sigma $
(see \cite{Donagi} and \cite{FMW}). Reversing this process one would
construct a degenerated K3 surface from each $E_{8}\times E_{8}$ flat bundle
over $\Sigma $. To globalize this construction, given a stable $E_{8}\times
E_{8}$-bundle over an elliptically fibered Calabi-Yau threefold, we expect
to obtain a degenerated Calabi-Yau fourfold with a K3 fibration. Duality
from physics predicts that certain string theory or F-theory on these spaces
are equivalent to one another.

\bigskip

On the other hand, if $X_{n}$ is embedded in a Calabi-Yau threefold $M$ and
we choose a family of Ricci flat metrics on $M$ so that the size of $X_{n}$
goes to zero (i.e. blowing down $X_{n}$ inside $M$). Then M-theory duality
predicts that enhanced gauge symmetry exists in the limit. This phenomenon
might be closely related to our bundle $\mathcal{E}_{n}$, $\mathcal{L}_{n}$
and $\mathcal{R}_{n}$ and their various reductions.

\bigskip

Now we come back to the focus of this paper. To discuss the reduction of $%
\mathcal{E}_{n}$ into various smaller subgroups, it is easier to deal with
the associated Lie algebra bundle $L\mathcal{E}_{n}$ which equals $L\mathcal{%
E}_{n}=O_{X}^{\oplus n}\bigoplus_{D}O_{X}\left( D\right) $ where the
summation is over those divisors $D$ satisfying $D^{2}=-2$ and $D\cdot K=0$.
Throughout the paper, we discuss several ways to construct the fiberwise Lie
algebra structure on $O_{X}^{\oplus n}\bigoplus_{D}O_{X}\left( D\right) $ .

\pagebreak

For the rest of this introduction, we will describe decompositions of $L%
\mathcal{E}_{n}$ and its representation bundles associated to various
specializations or degenerations of $X_{n}$. Readers should consult
individual sections for discussions of their geometric meanings. First we
consider the situation where we blow down a line $L$ on $X_{n}$. We write
the blowdown morphism as $\pi :X_{n}\rightarrow X_{n-1}$. Then we have the
following results.

\begin{theorem}
If $L$ is a line on $X_{n}$, then we have the following decompositions as
representation bundles of $\pi ^{\ast }L\mathcal{E}_{n-1}$ (or $\pi ^{\ast }L%
\mathcal{E}_{7}+L\mathcal{A}_{1}$ when $n=8$). 
\begin{eqnarray*}
\medskip L\mathcal{E}_{n} &=&\pi ^{\ast }L\mathcal{E}_{n-1}+O_{X_{n}}+\pi
^{\ast }\mathcal{L}_{n-1}\otimes O\left( -L\right) +\pi ^{\ast }\mathcal{L}%
_{n-1}^{\ast }\otimes O\left( L\right) ,\text{ for }n\leq 7, \\
L\mathcal{E}_{8} &=&\pi ^{\ast }L\mathcal{E}_{7}+L\mathcal{A}_{1}+\pi ^{\ast
}\mathcal{L}_{7}\otimes O\left( -L\right) \otimes \Lambda _{1}\text{.}
\end{eqnarray*}
Here $\mathcal{A}_{1}$ is the automorphism bundle of $\Lambda _{1}=O+O\left(
L+K\right) $ preserving its determinant.

For $\mathcal{L}_{n}$ we have 
\begin{equation*}
\begin{tabular}{ll}
$\medskip \mathcal{L}_{n}=\pi ^{\ast }\mathcal{L}_{n-1}$ & $+\pi ^{\ast }%
\mathcal{R}_{n-1}\otimes O\left( -L\right) +O\left( L\right) ,\text{ for }%
n\leq 6,$ \\ 
$\medskip \mathcal{L}_{7}=\pi ^{\ast }\mathcal{L}_{6}$ & $+\pi ^{\ast }%
\mathcal{R}_{6}\otimes O\left( -L\right) +O\left( L\right) +O\left(
-K-L\right) ,$ \\ 
$\mathcal{L}_{8}=\pi ^{\ast }\mathcal{L}_{7}\otimes \Lambda _{1}^{\ast }$ & $%
+\pi ^{\ast }\mathcal{R}_{7}\otimes O\left( -L\right) +\mathcal{A}%
_{1}\otimes O\left( -K\right) \text{.}$%
\end{tabular}
\end{equation*}
These decompositions describe the configuration of lines on $X_{n}$ with
respect to the fixed line $L$.
\end{theorem}

In particular, if we restrict our attention to a fiber over any point on $%
X_{n}$, we recover the following Lie algebra facts:

\begin{remark}
$E_{n-1}$ is a Lie subalgebra of $E_{n}$ for $n\leq 7$ whose Dynkin diagram
is obtained by removing the node in the Dynkin diagram of $E_{n}$ which
corresponds to the fundamental representation $\mathbf{L}_{n}$. When $n=8$, $%
E_{8}$ has a Lie subalgebra $E_{7}+\mathbf{\ }A_{1}$. As a representation of 
$E_{n-1}$ (or $E_{7}+A_{1}$), we have the following decomposition of $E_{n}.$%
\begin{eqnarray*}
E_{n} &=&E_{n-1}+\mathbf{1}+\mathbf{L}_{n-1}+\mathbf{L}_{n-1}^{\ast },\text{
for }n\leq 7, \\
E_{8} &=&E_{7}+A_{1}+\mathbf{L}_{7}\otimes \Lambda _{1}\text{.}
\end{eqnarray*}
For $\mathbf{L}_{n}$, we have 
\begin{eqnarray*}
\mathbf{L}_{n} &=&\mathbf{L}_{n-1}+\mathbf{R}_{n-1}+\mathbf{1,}\text{for }%
n\leq 6, \\
\mathbf{L}_{7} &=&\mathbf{L}_{6}+\mathbf{R}_{6}+\mathbf{1+1,} \\
\mathbf{L}_{8} &=&\mathbf{L}_{7}\otimes \Lambda _{1}+\mathbf{R}_{7}+A_{1}%
\text{.}
\end{eqnarray*}
\end{remark}

\newpage

Second we consider the situation where we specify a ruling $R$ on $X_{n}$.
In this case we can consider blowup of $\mathbb{P}^{2}$ at an arbitrary
number of points. That is there is no restriction on $n$ for our surface $%
X_{n}$.We have the following results.

\begin{theorem}
(1) If $R$ determines a ruling on $X_{n}$, then the rank $2n-2$ vector
bundle 
\begin{equation*}
\mathcal{W}_{n-1}=\bigoplus_{\substack{ C^{2}=-1 \\ CK=-1 \\ CR=0}}%
O_{X}\left( C\right) ,
\end{equation*}
carries a fiberwise non-degenerate quadratic form 
\begin{equation*}
q_{n-1}:\mathcal{W}_{n-1}\otimes \mathcal{W}_{n-1}\rightarrow
O_{X_{n}}\left( R\right) .
\end{equation*}
The automorphism bundle of $\mathcal{W}_{n-1}$ preserving $q_{n-1}$ is a $%
D_{n-1}$-bundle $\mathcal{D}_{n-1}$ over $X_{n}$ whose associated Lie
algebra bundle equals 
\begin{equation*}
L\mathcal{D}_{n-1}=\Lambda ^{2}\mathcal{W}_{n-1}\otimes O_{X_{n}}\left(
-R\right) =O_{X}^{\oplus n-1}\bigoplus_{\substack{ C^{2}=-2 \\ CK=0 \\ CR=0}}%
O_{X}\left( C\right) .
\end{equation*}
The representation bundle of $L\mathcal{D}_{n-1}$ corresponding to spinor
representations are 
\begin{equation*}
\mathcal{S}^{+}=\bigoplus_{\substack{ S^{2}=-1 \\ SK=-1 \\ SR=1}}%
O_{X_{n}}\left( S\right) \,\,\text{ and \thinspace \thinspace }\mathcal{S}%
^{-}=\bigoplus_{\substack{ T^{2}=-2 \\ TK=0 \\ TR=1}}O_{X_{n}}\left(
T\right) .
\end{equation*}
They are related by Clifford multiplication homomorphisms 
\begin{equation*}
\mathcal{S}^{+}\otimes \mathcal{W}_{n-1}^{\ast }\rightarrow \mathcal{S}^{-}%
\text{ and }\mathcal{S}^{-}\otimes \mathcal{W}_{n-1}\rightarrow \mathcal{S}%
^{+}\text{ ,}
\end{equation*}
When $n=2m$ is even, we have isomorphism 
\begin{equation*}
\left( \mathcal{S}^{+}\right) ^{\ast }\otimes O_{X_{2m}}\left( \left(
m-4\right) R-K_{X}\right) \cong \mathcal{S}^{-}.
\end{equation*}
and when $n=2m+1$ is odd, we have isomorphisms 
\begin{eqnarray*}
\medskip \left( \mathcal{S}^{+}\right) ^{\ast }\otimes O_{X_{2m+1}}\left(
\left( m-3\right) R-K_{X}\right)  &\cong &\mathcal{S}^{+}, \\
\left( \mathcal{S}^{-}\right) ^{\ast }\otimes O_{X_{2m+1}}\left( \left(
m-4\right) R-K_{X}\right)  &\cong &\mathcal{S}^{-}.
\end{eqnarray*}

(2) When $n\leq 8$ then $L\mathcal{D}_{n-1}$ is a Lie algebra subbundle of $L%
\mathcal{E}_{n}$. For $n=7$ in fact $L\mathcal{E}_{7}$ has a Lie algebra
subbundle $L\mathcal{D}_{6}+L\mathcal{A}_{1}$ and when $n=8,$ $L\mathcal{E}%
_{8}$ has a Lie algebra subbundle $L\mathcal{D}_{8}$. We can decompose
representation bundles of $L\mathcal{E}_{n}$ under $L\mathcal{D}_{n-1}$ as
follows: For $\mathcal{L}_{n}$ we have 
\begin{eqnarray*}
\medskip \mathcal{L}_{n} &=&\mathcal{W}_{n-1}+\mathcal{S}^{+}\text{ when }%
n\leq 5, \\
\medskip \mathcal{L}_{6} &=&\mathcal{W}_{5}+\mathcal{S}^{+}+O_{X_{6}}\left(
-K_{X}-R\right) , \\
\medskip \mathcal{L}_{7} &=&\mathcal{W}_{6}\otimes \Lambda _{1}+\mathcal{S}%
^{+}, \\
\mathcal{L}_{8} &=&\mathcal{W}_{8}+\mathcal{S}^{+}\text{.}
\end{eqnarray*}
This describes the relationship between the ruling and configuration of
lines on $X_{n}$. For $\mathcal{R}_{n}$ we have 
\begin{eqnarray*}
\medskip \mathcal{R}_{n} &=&O\left( R\right) \left( O+\mathcal{S}^{-}\right)
,\text{ \thinspace \thinspace \thinspace for }n\leq 4, \\
\medskip \mathcal{R}_{5} &=&O\left( R\right) \left( O+\mathcal{S}%
^{-}+O\left( -K-2R\right) \right) , \\
\medskip \mathcal{R}_{6} &=&O\left( R\right) \left( O+\mathcal{S}^{-}+%
\mathcal{W}_{5}\otimes O\left( -K-2R\right) \right) , \\
\mathcal{R}_{7} &=&O\left( R\right) \left( S^{2}\Lambda _{1}+\mathcal{S}%
^{-}\otimes \Lambda _{1}+\Lambda ^{2}\mathcal{W}_{6}\otimes \left(
-K-2R\right) \right) \text{.}
\end{eqnarray*}
This describes the configuration of ruling on $X_{n}$ with respect to a
fixed ruling. For $L\mathcal{E}_{n}$ we have 
\begin{eqnarray*}
\medskip L\mathcal{E}_{n} &=&L\mathcal{D}_{n-1}+O+\mathcal{S}^{-}+\mathcal{S}%
^{-}\otimes O\left( \left( 4-m\right) R+K\right) ,\text{ \thinspace
\thinspace \thinspace for }n=2m+1\leq 5, \\
\medskip L\mathcal{E}_{n} &=&L\mathcal{D}_{n-1}+O+\mathcal{S}^{-}+\mathcal{S}%
^{+}\otimes O\left( \left( 4-m\right) R+K\right) ,\text{ \thinspace
\thinspace \thinspace for }n=2m\leq 6, \\
\medskip L\mathcal{E}_{7} &=&L\mathcal{D}_{7}+L\mathcal{A}_{1}+\mathcal{S}%
^{-}\otimes \Lambda _{1}\otimes O\left( R+K\right) , \\
L\mathcal{E}_{8} &=&L\mathcal{D}_{8}+\mathcal{S}^{+}\text{.}
\end{eqnarray*}
\end{theorem}

In particular, if we restrict our attention to a fiber over any point on $%
X_{n}$, we recover the following Lie algebra facts:

\begin{remark}
(1) The space of infinitesimal automorphisms of a non-degenerate quadratic
form on $W_{n-1}\cong \mathbb{C}^{n-1}$ is a Lie algebra of type $D_{n-1}=%
\mathbf{so}\left( 2n-2\right) $. We can identify this Lie algebra with $%
\Lambda ^{2}W_{n-1}$. There are Clifford multiplication homomorphisms
between the two spinor representation of $D_{n-1}:$%
\begin{equation*}
S^{+}\otimes W_{n-1}^{\ast }\rightarrow S^{-}\text{ and }S^{-}\otimes
W_{n-1}\rightarrow S^{+}\text{ .}
\end{equation*}
When $n=2m$ is even, we have isomorphism 
\begin{equation*}
\left( S^{+}\right) ^{\ast }\cong S^{-}.
\end{equation*}
and when $n=2m+1$ is odd, we have isomorphisms 
\begin{equation*}
\left( S^{+}\right) ^{\ast }\cong S^{+}\text{ and }\left( S^{-}\right)
^{\ast }\cong S^{-}.
\end{equation*}
(2) When $n\leq 8$ then $D_{n-1}$ is a Lie subalgebra of $E_{n}$ whose
Dynkin diagram is obtained by removing the node in the Dynkin diagram of $%
E_{n}$ which corresponds to the fundamental representation $\mathbf{R}_{n}$.
For $n=7$ in fact $E_{7}$ has a Lie subalgebra $D_{6}+A_{1}$ and when $n=8,$ 
$E_{8}$ has a Lie subalgebra $D_{8}$. We can decompose representations of $%
E_{n}$ under $D_{n-1}$ as follows: For $\mathbf{L}_{n}$ we have 
\begin{eqnarray*}
\mathbf{L}_{n} &=&W_{n-1}+S^{+}\text{ when }n\leq 5, \\
\mathbf{L}_{6} &=&W_{5}+S^{+}+\mathbf{1}, \\
\mathbf{L}_{7} &=&W_{6}\otimes \Lambda _{1}+S^{+}, \\
\mathbf{L}_{8} &=&W_{8}+S^{+}\text{.}
\end{eqnarray*}
For $\mathbf{R}_{n}$ we have 
\begin{eqnarray*}
\mathbf{R}_{n} &=&\mathbf{1}+S^{-},\text{ \thinspace \thinspace \thinspace
for }n\leq 4 \\
\mathbf{R}_{5} &=&\mathbf{1}+S^{-}+\mathbf{1}, \\
\mathbf{R}_{6} &=&\mathbf{1}+S^{-}+W_{5}, \\
\mathbf{R}_{7} &=&\mathbf{1\otimes }S^{2}\Lambda _{1}+S^{-}\otimes \Lambda
_{1}+\Lambda ^{2}W_{6}\otimes \mathbf{1}\text{.}
\end{eqnarray*}
For $E_{n}$ we have 
\begin{eqnarray*}
E_{n} &=&D_{n-1}+O+S^{-}+S^{-},\text{ \thinspace \thinspace \thinspace for }%
n=2m+1\leq 5, \\
E_{n} &=&D_{n-1}+O+S^{-}+S^{+},\text{ \thinspace \thinspace \thinspace for }%
n=2m\leq 6, \\
E_{7} &=&D_{7}+A_{1}+S^{-}\otimes \Lambda _{1}, \\
E_{8} &=&D_{8}+S^{+}\text{.}
\end{eqnarray*}
\end{remark}

\newpage

Third we again consider the situation where we specify a ruling $R$ on $%
X_{n} $. However we also fix a line section (resp. ruling section) of this
ruling. This corresponds to choosing a direct summand of $\mathcal{S}^{+}$
(resp. $\mathcal{S}^{-}$). We have the following results.

\begin{theorem}
(1) If $R$ determines a ruling on $X_{n}$ and $O\left( S\right) $ is a
direct summand of $\mathcal{S}^{+}$, then $S$ is a line section of the
ruling. We consider the rank $n-1$ bundle 
\begin{equation*}
\mathcal{\Lambda }_{n-2}=\bigoplus_{\substack{ C^{2}=-1 \\ CK=-1 \\ CR=0\,\,
\\ CS=1\,\,}}O_{X_{n}}\left( C\right) ,
\end{equation*}
which has determinant equals $O_{X_{n}}\left( -K_{X}-2S+\left( n-4\right)
R\right) $. The automorphism bundle of $\mathcal{\Lambda }_{n-2}$ preserving
its determinant is a $A_{n-2}$-bundle $\mathcal{A}_{n-2}$ over $X_{n}$. Then 
$\mathcal{A}_{n-2}$ is a principal subbundle of $\mathcal{D}_{n-1}$. If $%
\mathcal{\Lambda }_{n-2}^{l}$ denotes the $l^{th}$ exterior power of $%
\mathcal{\Lambda }_{n-2}.$ We have the following decompositions of
representation bundles of $L\mathcal{D}_{n-1}$ under $L\mathcal{A}_{n-2}$. 
\begin{eqnarray*}
\medskip L\mathcal{D}_{n-1} &=&L\mathcal{A}_{n-2}+O+\mathcal{\Lambda }%
_{n-2}^{2}\otimes O\left( -R\right) +\left( \mathcal{\Lambda }%
_{n-2}^{2}\right) ^{\ast }\otimes O\left( R\right) , \\
\mathcal{W}_{n-1} &=&\mathcal{\Lambda }_{n-2}+\mathcal{\Lambda }_{n-2}^{\ast
}\otimes O_{X_{n}}\left( R\right)  \\
\medskip  &=&\mathcal{\Lambda }_{n-2}+\mathcal{\Lambda }_{n-2}^{n-2}\otimes
O\left( K+2S+\left( 5-n\right) R\right) , \\
\medskip \mathcal{S}^{+} &=&\sum_{l=0}^{\left[ \frac{n-1}{2}\right] }%
\mathcal{\Lambda }_{n-2}^{2l}\otimes O\left( S-lR\right) , \\
\mathcal{S}^{-} &=&\sum_{l=1}^{\left[ \frac{n}{2}\right] }\mathcal{\Lambda }%
_{n-2}^{2l-1}\otimes O\left( S-lR\right) .
\end{eqnarray*}
When $n\leq 7$, these decompositions describe the configuration of lines and
rulings on $X_{n}$ in relation to the ruling $R$ with a line section $S$.

(2) If $R$ determines a ruling on $X_{n}$ and $O\left( T\right) $ is a
direct summand of $\mathcal{S}^{-}$, then $R+T$ is another ruling on $X_{n}$
and it is also a section of the original ruling $R$. We consider the rank $%
n-1$ bundle 
\begin{equation*}
\mathcal{\bar{\Lambda}}_{n-2}=\bigoplus_{\substack{ C^{2}=-1 \\ CK=-1 \\ %
CR=0\,\, \\ CT=1\,\,}}O_{X_{n}}\left( C\right) 
\end{equation*}
which has determinant equals $O_{X_{n}}\left( -K_{X}-2T+\left( n-5\right)
R\right) $. The automorphism bundle of $\mathcal{\bar{\Lambda}}_{n-2}$
preserving its determinant is a principal $A_{n-2}$-bundle $\mathcal{\bar{A}}%
_{n-2}$ over $X_{n}$. Then $\mathcal{\bar{A}}_{n-2}$ is a principal
subbundle of $\mathcal{D}_{n-1}$. If $\mathcal{\bar{\Lambda}}_{n-2}^{l}$
denotes the $l^{th}$ exterior power of $\mathcal{\bar{\Lambda}}_{n-2}$. We
have the following decompositions of representation bundles of $L\mathcal{D}%
_{n-1}$ under $L\mathcal{\bar{A}}_{n-2}$. 
\begin{eqnarray*}
\medskip L\mathcal{D}_{n-1} &=&L\mathcal{\bar{A}}_{n-2}+O+\mathcal{\bar{%
\Lambda}}_{n-2}^{2}\otimes O\left( -R\right) +\left( \mathcal{\bar{\Lambda}}%
_{n-2}^{2}\right) ^{\ast }\otimes O\left( R\right) , \\
\mathcal{W}_{n-1} &=&\mathcal{\bar{\Lambda}}_{n-2}+\mathcal{\bar{\Lambda}}%
_{n-2}^{\ast }\otimes O_{X_{n}}\left( R\right)  \\
\medskip  &=&\mathcal{\bar{\Lambda}}_{n-2}+\mathcal{\bar{\Lambda}}%
_{n-2}^{n-2}\otimes O\left( K+2T+\left( 6-n\right) R\right)  \\
\medskip \mathcal{S}^{+} &=&\sum_{l=1}^{\left[ \frac{n}{2}\right] }\mathcal{%
\bar{\Lambda}}_{n-2}^{2l-1}\otimes O\left( T-\left( l-1\right) R\right) . \\
\mathcal{S}^{-} &=&\sum_{l=0}^{\left[ \frac{n-1}{2}\right] }\mathcal{\bar{%
\Lambda}}_{n-2}^{2l}\otimes O\left( T-lR\right) .
\end{eqnarray*}
When $n\leq 7$, these decompositions describe the configuration of lines and
rulings on $X_{n}$ in relation to the rulings $R$ and $R+T$ .
\end{theorem}

In particular, if we restrict our attention to a fiber over any point on $%
X_{n}$, we recover the following Lie algebra facts:

\begin{remark}
(1) $A_{n-2}=\mathbf{sl}\left( n-1\right) $ is a Lie subalgebra of $D_{n-1}=%
\mathbf{so}\left( 2n-2\right) $ whose Dynkin diagram is obtained by removing
the node in the Dynkin diagram of $D_{n-1}$ which corresponds to the
fundamental representation $S^{+}$. We have the following decomposition of $%
D_{n-1}$ representations under $A_{n-2}$%
\begin{eqnarray*}
\medskip D_{n-1} &=&A_{n-2}+\mathbf{1}+\Lambda _{n-2}^{2}+\left( \Lambda
_{n-2}^{2}\right) ^{\ast }, \\
\medskip W_{n-1} &=&\Lambda _{n-2}+\Lambda _{n-2}^{\ast }=\Lambda
_{n-2}+\Lambda _{n-2}^{n-2}, \\
S^{+} &=&\sum_{l=0}^{\left[ \frac{n-1}{2}\right] }\Lambda _{n-2}^{2l}\text{
and }S^{-}=\sum_{l=1}^{\left[ \frac{n}{2}\right] }\Lambda _{n-2}^{2l-1}.
\end{eqnarray*}
Here $\Lambda _{n-2}^{l}$ denotes the $l^{th}$ exterior power of the
standard representation of $A_{n-2}$.

(2) $A_{n-2}$ is also a Lie subalgebra of $D_{n-1}$ whose Dynkin diagram is
obtained by removing the node in the Dynkin diagram of $D_{n-1}$ which
corresponds to the fundamental representation $S^{-}$. We have the following
decomposition of $D_{n-1}$ representations under $A_{n-2}$%
\begin{eqnarray*}
\medskip D_{n-1} &=&A_{n-2}+\mathbf{1}+\Lambda _{n-2}^{2}+\left( \Lambda
_{n-2}^{2}\right) ^{\ast }, \\
\medskip W_{n-1} &=&\Lambda _{n-2}+\Lambda _{n-2}^{\ast }=\Lambda
_{n-2}+\Lambda _{n-2}^{n-2}, \\
S^{+} &=&\sum_{l=1}^{\left[ \frac{n}{2}\right] }\Lambda _{n-2}^{2l-1}\text{
and }S^{-}=\sum_{l=0}^{\left[ \frac{n-1}{2}\right] }\Lambda _{n-2}^{2l}.
\end{eqnarray*}
\end{remark}

\newpage

Fourth we consider degenerating $X_n$ to a nonnormal del Pezzo surface. We
have the following result.

\begin{theorem}
We consider a family of smooth del Pezzo surfaces $X_{n}\left( t\right) $
degenerating to a nonnormal surface $X_{n}\left( 0\right) $. Let $Z$ be the
limit of the intersection of $X_{n}\left( t\right) $ and the singular locus
of $X_{n}\left( 0\right) $ as $t$ goes to zero. We write $\mathcal{I}%
=\bigoplus_{q\in Z}\mathcal{I}_{\left\{ q\right\} }$.

(1) When $n=5$ and we degenerate the quartic surface $X_{5}$ inside $\mathbb{%
P}^{4}$ into a union of two quadratic surfaces $X_{5}\left( 0\right)
=Q_{\left( 1\right) }\cup Q_{\left( 2\right) }$. Let $\mathcal{R}_{\left(
i\right) }=\bigoplus_{R}O_{Q_{\left( i\right) }}\left( R\right) $ with $R$
satisfying $R^{2}=0$ and $R\cdot K_{Q_{\left( i\right) }}=-2$, then we write 
\begin{eqnarray*}
\mathcal{L}_{5}\left( 0\right)  &=&\mathcal{I}\otimes \mathcal{R}_{\left(
1\right) }+\mathcal{I}\otimes \mathcal{R}_{\left( 1\right) }, \\
\mathcal{R}_{5}\left( 0\right)  &=&\mathcal{I\wedge I+R}_{\left( 1\right)
}\otimes \mathcal{R}_{\left( 2\right) }.
\end{eqnarray*}
These decomposition describe the limit configurations of lines and rulings
on $X_{5}\left( t\right) $ as $t$ approach zero.

(2) When $n=6$ and we degenerate the cubic surface $X_{6}$ into a union of
three planes $X_{6}\left( 0\right) =H_{1}\cup H_{2}\cup H_{3}$. Then we
write 
\begin{equation*}
\mathcal{L}_{6}\left( 0\right) =\mathcal{I}_{\left( 1\right) }\otimes 
\mathcal{I}_{\left( 2\right) }\otimes O_{H_{3}}\left( 1\right) +\mathcal{I}%
_{\left( 2\right) }\otimes \mathcal{I}_{\left( 3\right) }\otimes
O_{H_{1}}\left( 1\right) +\mathcal{I}_{\left( 3\right) }\otimes \mathcal{I}%
_{\left( 1\right) }\otimes O_{H_{2}}\left( 1\right) .
\end{equation*}

This decomposition describes the limit configuration of lines on $%
X_{6}\left( t\right) $ as $t$ approach zero. Moreover the structure of the
triple product $c_{6}\left( 0\right) :\mathcal{L}_{6}\left( 0\right) \otimes 
\mathcal{L}_{6}\left( 0\right) \otimes \mathcal{L}_{6}\left( 0\right)
\rightarrow O_{X_{6}\left( 0\right) }\left( 1\right) $ can also be described
in very explicit terms. Similar results hold true for $\mathcal{R}_{6}$ via $%
\mathcal{R}_{6}=\mathcal{L}_{6}^{\ast }\otimes O\left( -K\right) $.

(3) When $n=6$ and we degenerate the cubic surface $X_{6}$ into a union of a
plane and a quadratic surface $X_{6}\left( 0\right) =H\cup Q$. Let $\mathcal{%
R}_{Q}=\oplus _{R}O_{Q}\left( R\right) $ for $R$ satisfying $R^{2}=0$ and $%
R\cdot K_{Q}=-2$. Then we write 
\begin{equation*}
\mathcal{L}_{6}\left( 0\right) =\Lambda ^{2}\mathcal{I}\otimes O_{H}\left(
1\right) +\mathcal{I}\otimes \mathcal{R}_{Q}.
\end{equation*}
This decomposition describes the limit configuration of lines on $%
X_{6}\left( t\right) $ as $t$ approach zero. Moreover the structure of the
triple product 
\begin{equation*}
c_{6}\left( 0\right) :\mathcal{L}_{6}\left( 0\right) \otimes \mathcal{L}%
_{6}\left( 0\right) \otimes \mathcal{L}_{6}\left( 0\right) \rightarrow
O_{X_{6}\left( 0\right) }\left( 1\right) ,
\end{equation*}
can also be described in explicit terms. Similar results hold true for $%
\mathcal{R}_{6}$ via $\mathcal{R}_{6}=\mathcal{L}_{6}^{\ast }\otimes O\left(
-K\right) $.

(4) When $n=7$ and we degenerate $X_{7}$ into a union of two copies of $%
\mathbb{P}^{2}$ joining along a conic curve. Then we write 
\begin{equation*}
\mathcal{L}_{7}\left( 0\right) =\Lambda ^{2}\mathcal{I}\otimes O_{\mathbb{P}%
_{\left( 1\right) }^{2}}\left( 1\right) +\Lambda ^{2}\mathcal{I}\otimes O_{%
\mathbb{P}_{\left( 2\right) }^{2}}\left( 1\right) .
\end{equation*}
This decomposition describes the limit configuration of lines on $%
X_{7}\left( t\right) $ as $t$ approach zero. Moreover the structure of the
quartic product 
\begin{equation*}
f_{7}:\mathcal{L}_{7}\left( 0\right) \otimes \mathcal{L}_{7}\left( 0\right)
\otimes \mathcal{L}_{7}\left( 0\right) \otimes \mathcal{L}_{7}\left(
0\right) \rightarrow O\left( -2K\right) ,
\end{equation*}
can also be described in explicit terms.
\end{theorem}

If we restrict our attention to fibers over any smooth point on $X_n\left(
0\right) $, they are related to the following Lie algebra facts:

\begin{remark}
(1) When $n=5$, $E_{5}=D_{5}=\mathbf{so}\left( 10\right) $ has a Lie
subalgebra $A_{3}\times A_{1}\times A_{1}=\mathbf{sl}\left( 4\right) \mathbf{%
+sl}\left( 2\right) \mathbf{+sl}\left( 2\right) $. We have the following
decomposition of $E_{5}$ representations under $A_{3}\times A_{1}\times
A_{1}.$%
\begin{eqnarray*}
\mathbf{L}_{5} &=&\Lambda _{3}\otimes \Lambda _{1}\otimes \mathbf{1}+\Lambda
_{3}^{\ast }\otimes \mathbf{1}\otimes \Lambda _{1}, \\
\mathbf{R}_{5} &=&\Lambda _{3}^{2}\otimes \mathbf{1}\otimes \mathbf{1}+%
\mathbf{1}\otimes \Lambda _{1}\otimes \Lambda _{1}.
\end{eqnarray*}
Here $\Lambda _{n}$ denotes the standard representation of $A_{n}=\mathbf{sl}%
\left( n\right) $.

(2) When $n=6$, $E_{6}$ has a Lie subalgebra $A_{2}\times A_{2}\times A_{2}=%
\mathbf{sl}\left( 3\right) \mathbf{+sl}\left( 3\right) \mathbf{+sl}\left(
3\right) $. We have the following decomposition of $E_{6}$ representation
under $A_{2}\times A_{2}\times A_{2}.$%
\begin{equation*}
\mathbf{L}_{6}=\Lambda _{2}\otimes \Lambda _{2}^{\ast }\otimes \mathbf{1}+%
\mathbf{1}\otimes \Lambda _{2}\otimes \Lambda _{2}^{\ast }+\Lambda
_{2}^{\ast }\otimes \mathbf{1}\otimes \Lambda _{2}.
\end{equation*}
Moreover the triple product $c_{6}:\mathbf{L}_{6}\otimes \mathbf{L}%
_{6}\otimes \mathbf{L}_{6}\rightarrow \mathbb{C}$ can be described
explicitly in terms of this decomposition. Similar structure holds for $%
\mathbf{R}_{6}=\mathbf{L}_{6}^{\ast }$.

(3) When $n=6$, $E_{6}$ has a Lie subalgebra $A_{1}\times A_{5}=\mathbf{sl}%
\left( 2\right) \mathbf{+sl}\left( 6\right) $. We have the following
decomposition of $E_{6}$ representation under $A_{1}\times A_{5}.$%
\begin{equation*}
\mathbf{L}_{6}=\Lambda _{5}^{2}\otimes \mathbf{1}+\Lambda _{5}\otimes
\Lambda _{2}.
\end{equation*}
Moreover the triple product $c_{6}:\mathbf{L}_{6}\otimes \mathbf{L}%
_{6}\otimes \mathbf{L}_{6}\rightarrow \mathbb{C}$ can be described
explicitly in terms of this decomposition. Similar structure holds for $%
\mathbf{R}_{6}=\mathbf{L}_{6}^{\ast }$.

(4) When $n=7$, $E_{7}$ has a Lie subalgebra $A_{7}=\mathbf{sl}\left(
8\right) $. We have the following decomposition of $E_{7}$ representation
under $A_{7}.$%
\begin{equation*}
\mathbf{L}_{7}=\Lambda _{7}^{2}+\Lambda _{2}^{6}.
\end{equation*}
Moreover the triple product $f_{7}:\mathbf{L}_{7}\otimes \mathbf{L}%
_{7}\otimes \mathbf{L}_{7}\otimes \mathbf{L}_{7}\rightarrow \mathbb{C}$ can
be described explicitly in terms of this decomposition.
\end{remark}

\textbf{Notations}: We will simply use $X$ to denote $X_{n}$ when there is
no confusion that might occur. The canonical divisor of $X_{n}$ is called $%
K_{X_{n}}$ or simply $K.$ If $D$ is a divisor on $X$, then there is an
associated line bundle $O\left( D\right) $ on $X$ together with a rational
section of $O\left( D\right) $. Such a section is canonical up to
multiplication by a non-zero scalar. We implicitly fix one such section for
each $D$. $H^{i}\left( O\left( D\right) \right) $ denotes the $i^{th}$
cohomology group of the sheaf of sections of $O\left( D\right) $ whose
dimension is denoted as $h^{i}\left( O\left( D\right) \right) $. Moreover $%
\chi \left( O\left( D\right) \right) $ denotes the Euler characteristic,
namely $\chi \left( O\left( D\right) \right) =\sum_{i=0}^{2}\left( -1\right)
^{i}h^{i}\left( O\left( D\right) \right) $. We often use $+$ to replace $%
\oplus $ to improve the visual effect of our equations.

Regarding notations from Lie theory, a fundamental representation means a
highest weight representation whose highest weight vector is the first
lattice point along an edge of the fundamental Weyl chamber. Such vectors
are in one-to-one correspondence with nodes of the Dynkin diagram. Therefore
we usually speak of a fundamental representation corresponding to some
particular node in the Dynkin diagram. When $n\leq 5$, then $E_n$ is a
classical Lie algebra and the standard representation of it refers to its
defining representation. For example the standard representation of $\mathbf{%
sl}\left( 2\right) $ is of dimension two and the standard representation of $%
\mathbf{so}\left( 2n-2\right) $ is of dimension $2n-2$.

Even though we use the complex number field, most arguments in this paper
work equally well over any algebraically closed field of characteristic zero.%
$\bigskip $\newpage

\section{$E_{n}$-bundles over del Pezzo surfaces}

In this section we study the bundle $\mathcal{L}_{n}$, $\mathcal{R}_{n}$ and 
$L\mathcal{E}_{n}$ over a del Pezzo surface $X_{n}$. Before we do this we
first discuss some general properties of del Pezzo surfaces. Many of these
properties can be found in [Beauville], [Hartshorne] and [Manin].

\subsection{Geometry of del Pezzo surfaces}

As we mentioned in the introduction, a del Pezzo surface has an ample
anti-canonical divisor class. We can represent $X_{n}$ as a blowup of $%
\mathbb{P}^{2}$ at $n$ generic points. When $n\leq 6$ then the complete
linear system of cubics through these $n$ points defines an embedding $%
X_{n}\subset \mathbb{P}^{9-n}.$ It is an anti-canonical embedding of degree $%
d=9-n$. For example $X_{5}$ is a complete intersection of two quadrics in $%
\mathbb{P}^{4}$ and $X_{6}$ is a cubic surface in $\mathbb{P}^{3}$. When $%
n=7 $ this linear system is not very ample on $X_{7}$. Instead it defines a
double cover of $\mathbb{P}^{2}$ branched along a quartic curve. On the
other hand we can embed $X_{7}$ inside the weighted projective space $%
\mathbb{P}\left( 1,1,1,2\right) $ and $X_{8}$ inside $\mathbb{P}\left(
1,1,2,3\right) $.

Suppose we fix a point $p$ on $X_{n}\subset \mathbb{P}^{9-n}$. The
projection from $p$ defines a morphism from the blowup of $X_{n}$ at $p$ to $%
\mathbb{P}^{9-n-1}$. We denote this blowup surface as $X_{n+1}$. Then this
morphism is the anti-canonical morphism for $X_{n+1}$.

Next we want to discuss lines and rulings on $X_{n}$. When $n\leq 6$ we have
the anti-canonical embedding $X_{n}\subset \mathbb{P}^{d}$. A curve in $X_{n}
$ is a line inside $\mathbb{P}^{d}$ if and only if it is an exceptional
curve in $X_{n}$. By abuse of notations we continue to call an exceptional
curve in $X_{n}$ a line even when $n=7$ or $8$. When $X_{n}$ is represented
as a blowup of $\mathbb{P}^{2}$ at generic points $p_{1},...,p_{n}$: 
\begin{equation*}
\pi :X_{n}\rightarrow \mathbb{P}^{2}\text{.}
\end{equation*}
Then we can describe (see \cite{Manin}) all lines in $X_{n}$ as follows: $C$
is a line in $X_{n}$ if and only if (i) $\pi \left( D\right) =p_{i}$; or
(ii) $\pi \left( D\right) $ is a line passes through $p_{i}$ and $p_{j}$; or
(iii) $\pi \left( D\right) $ is a conic passes through five of the $p_{i}$%
's; or (iv) $\pi \left( D\right) $ is a cubic passes through seven of the $%
p_{i}$'s and with one being a double point; or (v) $\pi \left( D\right) $ is
a quartic passes through 8 of the $p_{i}$'s and three being double points;
or (vi) $\pi \left( D\right) $ is a quintic passes through 8 of the $p_{i}$%
's and six being double points; or (vii) $\pi \left( D\right) $ is a sextic
passes through 8 of the $p_{i}$'s and seven being double points and one
triple point. The number of lines on each $X_{n}$ is given in the table in
the introduction. The following proposition (from \cite{Manin}) is a useful
numerical characterization of lines on $X_{n}.$

\begin{proposition}
If $D$ is a divisor on $X_{n}$, then there is a line linearly equivalent to $%
D$ if and only if $D^{2}=-1$ and $D\cdot K=-1$.
\end{proposition}

Proof of proposition: If $D$ is linearly equivalent to a line, then it
implies that $D^2=-1$ and $D\cdot K=-1$ by the adjunction formula.
Conversely if $D$ satisfies these two equalities, then the Riemann-Roch
formula give $\chi \left( O\left( D\right) \right) =1$. On the other hand $%
K-D$ is not an effective divisor because $\left( -K\right) \cdot \left(
K-D\right) <0$ and $-K$ is an ample divisor. By Serre duality, we have $%
h^2\left( O\left( D\right) \right) =0$. So $h^0\left( O\left( D\right)
\right) $ is at least one and we can assume that $D$ is an effective
divisor. In fact $D$ is an irreducible divisor because $D\cdot \left(
-K\right) =1$.

The genus formula showed that $D$ has arithmetic genus equals zero. It
implies that the irreducible divisor $D$ is in fact a smooth rational curve.
Now $D^{2}=-1$ implies that it is a line.

Therefore

\begin{eqnarray*}
\mathcal{L}_{n} &=&\bigoplus_{\substack{ l^{2}=-1  \\ lK=-1}}O\left(
l\right) \,\,\,\,\,\,\,\,\,\text{when }n\leq 7, \\
\mathcal{L}_{8} &=&\bigoplus_{_{\substack{ l^{2}=-1  \\ lK=-1}}}O\left(
l\right) +O\left( -K\right) ^{\oplus 8}.
\end{eqnarray*}

The following notion about certain types of configuration of lines on $X_{n}$
will be used later.

\begin{definition}
If the dual graph of a configuration of $d$ lines on $X_{n}$ is a d-gon,
i.e. a polygon with $d$ edges, then we call the configuration a d-gon.

A d-gon is called a triangle (resp. rectangle, pentagon, hexagon, septagon,
octagon) if $d$ equals three (resp. four, five, six, seven, eight).
\end{definition}

\begin{proposition}
\label{Prop: d-gon}(1) When $n\leq 5$ there is no triangle on $X_{n}$. When $%
n=6$ then every triangle on $X_{6}$ is an anti-canonical divisor.

(2) When $n\leq 4$ there is no rectangle on $X_{n}.$ When $n=5$ then every
rectangle on $X_{5}$ is an anti-canonical divisor.
\end{proposition}

Proof of proposition: The proof of part (1) and (2) are essentially the same
and we will only give one of them. If $C=l_1\cup l_2\cup l_3$ is a triangle
on $X_n$, then $C^2=3$ and $C\cdot K=-3$. Therefore if $n\leq 5$ then $%
K^2=9-n>3$ and hence $\left( C\cdot K\right) ^2=9<\left( C^2\right) \cdot
\left( K^2\right) .$ This violates the Hodge index theorem so no triangle
exists.

When $n=6$, we have $\left( C\cdot K\right) ^2=\left( C^2\right) \cdot
\left( K^2\right) .$ Again by Hodge index theorem $C$ must be linearly
equivalent to a multiple of $K$. Because $C\cdot K=-3$, we have $C\equiv -K$%
. Hence the result.

In fact the same proof shows that when $n<9-k$ there is no k-gon on $X_{n}$.
When $n=9-k$ then every k-gon on $X_{n}$ is an anti-canonical divisor.

\begin{corollary}
If $l_{1}$ and $l_{2}$ are two lines on $X_{n}$. Then $l_{1}\cdot l_{2}\neq
2 $ when $n\leq 6$.

If $l_{1}\cdot l_{2}=2$ and $n=7$, then $l_{2}$ is linearly equivalent to $%
-K-l_{1}$.
\end{corollary}

Next we are going to give a similar numerical characterization of rulings on 
$X_{n}$.

\begin{proposition}
If $D$ is a divisor on $X_{n}$, then $D$ determines a ruling on $X_{n}$ if
and only if $D^{2}=0$ and $D\cdot K=-2$.
\end{proposition}

Proof of proposition: If $D$ determines a ruling on $X$, then it implies
that $D^2=0$ and $D\cdot K=-2$ by the adjunction formula. Conversely if $D$
satisfies these two equalities, then the Riemann-Roch formula give $\chi
\left( O\left( D\right) \right) =2$. By similar arguments as above, we can
assume that $D$ is an effective divisor with at most two irreducible
components and its complete linear system is at least one dimensional. Again
the genus formula showed that $D$ has arithmetic genus equals zero.

We claim that $D$ has no fixed component. Otherwise $D=D_{0}+R$ where $R$ is
the fixed component and $D_{0}$ is the free component. Since $-K$ is ample
and $-K\cdot D=2$, we have $-K\cdot D_{0}=-K\cdot R=1$. Since $p_{a}\left(
D\right) =0$ and $R^{2}<0$, $R$ must be a line on $X$. Using the equality 
\begin{eqnarray*}
0 &=&D^{2} \\
&=&D_{0}^{2}+2D_{0}\cdot R-1,
\end{eqnarray*}
and the fact that $D_{0}^{2}\geq 0$, we must have $D_{0}^{2}=1$ and $%
D_{0}\cdot R=0$. Now the adjunction formula implies that $p_{a}\left(
D_{0}\right) =1.$ This contradicts to $p_{a}\left( D\right) =0$.

This linear system defines a fibration of $X$ whose generic fiber is an
irreducible curve with $p_{a}\left( D\right) =0$ (and hence a smooth
rational curve). That is $D$ determines a ruling on $X$. $\square $

\bigskip

Hence we have 
\begin{equation*}
\mathcal{R}_{n}=\tbigoplus_{\substack{ R^{2}=0  \\ RK=-2}}O\left( R\right)
\,\,\,\,\,\,\,\,\,\text{when }n\leq 6,
\end{equation*}
We also define 
\begin{equation*}
\mathcal{R}_{7}=\tbigoplus_{_{\substack{ R^{2}=0  \\ RK=-2}}}O\left(
R\right) +O\left( -K\right) ^{\oplus 7}.
\end{equation*}

Now we want to describe a vector bundle homomorphism 
\begin{equation*}
c_{n}:\mathcal{L}_{n}\otimes \mathcal{L}_{n}\rightarrow \mathcal{R}_{n}.
\end{equation*}
Suppose that $l_{1}$ and $l_{2}$ are lines on $X_{n}$, then the divisor $%
R=l_{1}+l_{2}$ determines a ruling on $X_{n}$ if and only if $l_{1}$ and $%
l_{2}$ intersect at one point. It is because 
\begin{equation*}
R\cdot K=\left( l_{1}+l_{2}\right) \cdot K=-1-1=-2
\end{equation*}
and 
\begin{equation*}
R^{2}=\left( l_{1}+l_{2}\right) ^{2}=l_{1}^{2}+l_{2}^{2}+2l_{1}\cdot
l_{2}=-2+2l_{1}\cdot l_{2}.
\end{equation*}
Therefore $R^{2}=0$ if and only if $l_{1}\cdot l_{2}=1$. In this case there
is an isomorphism $O_{X}\left( l_{1}\right) \otimes O_{X}\left( l_{2}\right) 
\overset{\simeq }{\rightarrow }O_{X}\left( R\right) $.Combining these
isomorphisms for various intersecting pairs of lines, we obtain the bundles
homomorphism 
\begin{equation*}
c_{n}:\mathcal{L}_{n}\otimes \mathcal{L}_{n}\rightarrow \mathcal{R}_{n},
\end{equation*}
for $n\leq 6$. When $n=7$ two lines $l_{1}$ and $l_{2}$ can intersect at
more than one point. When this happens, $l_{1}+l_{2}$ is an anti-canonical
divisor on $X_{7}$. Combining with this we obtain the homomorphism $\mathcal{%
L}_{7}\otimes \mathcal{L}_{7}\rightarrow \mathcal{R}_{7}$. See later section
for details. One should notice that there could be more than one pair of
lines that determines the same ruling on $X_{n}$. In fact, we will show
later that, each ruling on $X_{n}$ comes from exactly $n-1$ pairs of lines.

\bigskip

\subsection{Construction of $L\mathcal{E}_{n}$ over $X_{n}$}

Next we construct the following holomorphic vector bundle $L\mathcal{E}_{n}$%
: 
\begin{equation*}
L\mathcal{E}_{n}=O^{\oplus n}+\tbigoplus_{\substack{ D^{2}=-2 \\ DK=0}}%
O\left( D\right) 
\end{equation*}
In Chapter four of \cite{Manin}, Manin discussed properties of these
divisors $D$ and their relationship with the root system of $E_{n}$. Here we
simply use them to construct a Lie algebra bundle of type $E_{n}$ over $X_{n}
$.

\begin{remark}
As a holomorphic vector bundle over $X_{n}$, $\mathcal{L}_{n}$ (resp. $%
\mathcal{R}_{n}$ and $L\mathcal{E}_{n}$) is a direct sum of line bundles
which have the same degree with respect to the anti-canonical polarization.
Therefore it is a semi-stable vector bundle over $X_{n}.$
\end{remark}

\qquad In the following proposition we show that the above numerical
criterion on a divisor class $D$ characterizes those divisor classes which
can be written as the difference of two lines on $X_{n}$.

\begin{proposition}
(1) If $D$ is a divisor on $X_{n}$ with $n\leq 7$, then $D^{2}=-2$ and $%
D\cdot K=0$ if and only if $D\equiv l-l^{\prime }$ for some lines $l$ and $%
l^{\prime }$ on $X_{n}$.

(2) If $D$ is a divisor on $X_{8}$, then $D^{2}=-2$ and $D\cdot K=0$ if and
only if $D+K$ is a line.
\end{proposition}

Proof of proposition: Part (2) follows immediately from the previous
characterization of a line on $X_n$ and $K^2=1$ on $X_8$.

Now we assume that $n\leq 7$. Since the Neron-Severi group of $X_n$ is
spanned by the canonical divisor $K$ and the collection of lines on $X_n$, a
divisor with zero intersection with $K$ and every line on $X_n$ would
necessarily be numerically trivial. However $D^2=-2$, there must be at least
one line $l$ on $X_n$ such that $D\cdot l$ is nonzero.

Replacing $D$ by $-D$ if necessary, we can assume that $D\cdot l>0$. If we
write $D^{\prime }=D+l$ then $D^{\prime 2}=-3+2D\cdot l\geq -1$. Also $%
D^{\prime }\cdot K=-1$. Applying the Hodge index theorem, we have $\left(
D^{\prime 2}\right) \left( K^2\right) \leq \left( D^{\prime }\cdot K\right)
^2=1$. Therefore $D^{\prime 2}=-1$, $0$ or $1$. However $D^{\prime 2}$
cannot be zero since $\acute{D}^{^{\prime }2}$ and $D^{\prime }\cdot K$ have
the same parity.

If $D^{\prime 2}=1$, then the above inequality implies that $K^{2}\leq 1$.
This is impossible for $X_{n}$ with $n\leq 7$ because $K_{X_{n}}^{2}=9-n.$
Hence $D^{\prime 2}=-1$ and $D^{\prime }\cdot K=-1$. Namely $D^{\prime }$ is
a line $l^{\prime }$ and we have $D=D^{\prime }-l=l^{\prime }-l$. Hence we
proved the proposition. $\square $

\bigskip

Therefore we have 
\begin{eqnarray*}
L\mathcal{E}_{n} &=&O_{X}^{\oplus n}\bigoplus_{\substack{ D=l-l^{\prime } 
\\ l,l^{\prime }\text{ disjoint lines on }X}}O_{X}\left( D\right) \text{
when }n\leq 7. \\
L\mathcal{E}_{8} &=&\mathcal{L}_{8}\otimes O\left( K\right) \text{.}
\end{eqnarray*}

\qquad Now if we blow down a fixed line $L$ on $X_{n}$ then the above
description permits us to relate $L\mathcal{E}_{n}$ to the bundle $L\mathcal{%
E}_{n-1}$ on the blowdown surface $X_{n-1}$. Let us first denote the blow
down morphism as $\pi :X_{n}\rightarrow X_{n-1}$. Second it is obvious that $%
\pi ^{\ast }L\mathcal{E}_{n-1}+O_{X_{n}}$ is a subbundle of $L\mathcal{E}%
_{n} $. Third $\mathcal{L}_{n-1}$ is constructed in terms of lines on $%
X_{n-1}$ and therefore $\pi ^{\ast }\mathcal{L}_{n-1}\otimes O\left(
-L\right) $ and its dual are also subbundles of $L\mathcal{E}_{n}$. When $%
n\leq 7$ this exhausts the whole bundle $L\mathcal{E}_{n}$. 
\begin{equation*}
L\mathcal{E}_{n}=\pi ^{\ast }L\mathcal{E}_{n-1}+O_{X_{n}}+\pi ^{\ast }%
\mathcal{L}_{n-1}\otimes O\left( -L\right) +\pi ^{\ast }\mathcal{L}%
_{n-1}^{\ast }\otimes O\left( L\right) .
\end{equation*}
When $n=8$ the line bundle $O\left( K+L\right) $ and its dual also lie
inside $L\mathcal{E}_{8}$. We have

\begin{align*}
L\mathcal{E}_{8}& =\pi ^{\ast }L\mathcal{E}_{7}+O_{X_{8}}+\pi ^{\ast }%
\mathcal{L}_{7}\otimes O\left( -L\right) +\pi ^{\ast }\mathcal{L}_{7}^{\ast
}\otimes O\left( L\right) \\
& +O\left( -K-L\right) +O\left( K+L\right) .
\end{align*}

\qquad In general the maximum number of disjoint lines on $X_{n}$ equals $n$%
. If we fix such a maximal collection of lines on $X_{n}$ and label them $%
L_{1},L_{2},...,L_{n}$. Again we denote the blow down morphism of these
lines as $\pi :X_{n}\rightarrow \mathbb{P}^{2}$. Namely the exceptional
locus of $\pi $ is $L_{1}\cup L_{2}\cup ...\cup L_{n}$. We also denote the
pullback of the hyperplane class of $\mathbb{P}^{2}$ by $H$. Then, for $n<8$%
, we have 
\begin{align*}
L\mathcal{E}_{n}& =\bigoplus^{n}O_{X_{n}}\bigoplus_{i\neq j}O\left(
L_{i}-L_{j}\right) \\
& \bigoplus_{i<j<k}O\left( H-L_{i}-L_{j}-L_{k}\right) +O\left(
-H+L_{i}+L_{j}+L_{k}\right) \\
& \bigoplus_{i_{1}<...<i_{6}}O\left( 2H-\sum_{m=1}^{6}L_{i_{m}}\right)
+O\left( -2H+\sum_{m=1}^{6}L_{i_{m}}\right) .
\end{align*}
When $n=8$ we have 
\begin{align*}
L\mathcal{E}_{8}& =\bigoplus^{8}O_{X_{8}}\bigoplus_{i\neq j}O\left(
L_{i}-L_{j}\right) \\
& \bigoplus_{i<j<k}O\left( H-L_{i}-L_{j}-L_{k}\right) +O\left(
-H+L_{i}+L_{j}+L_{k}\right) \\
& \bigoplus_{i_{1}<...<i_{6}}O\left( 2H-\sum_{m=1}^{6}L_{i_{m}}\right)
+O\left( -2H+\sum_{m=1}^{6}L_{i_{m}}\right) \\
& \bigoplus_{i=1}^{8}O\left( 3H-\sum_{j=1}^{8}L_{j}-L_{i}\right) +O\left(
-3H+\sum_{j=1}^{8}L_{j}+L_{i}\right) ,
\end{align*}
Similarly we can decompose $\mathcal{L}_{n}$, when $n<8$, as 
\begin{equation*}
\mathcal{L}_{n}=\bigoplus_{i=1}^{n}O\left( L_{i}\right)
\bigoplus_{i<j}O\left( H-L_{i}-L_{j}\right)
\bigoplus_{i_{1}<...<i_{5}}O\left( 2H-\sum_{m=1}^{5}L_{i_{m}}\right) ,
\end{equation*}
and $\mathcal{L}_{8}=L\mathcal{E}_{8}\otimes O\left( -K\right) $.

\subsection{Fiberwise Lie algebra structure on $L\mathcal{E}_{n}$}

\qquad We remind the reader that the above decomposition of $L\mathcal{E}%
_{n} $ in terms of $L\mathcal{E}_{n-1}$ and $\mathcal{L}_{n-1}$ resembles
the decomposition of the adjoint representation of the Lie algebra $E_{n}$
under its maximal Lie subalgebra $E_{n-1}+\mathbf{u}\left( 1\right) $. When $%
n=8$, $E_{7}+\mathbf{u}\left( 1\right) $ is not a maximal Lie subalgebra of $%
E_{8}$. In fact $E_{7}+\mathbf{u}\left( 1\right) $ is contained in a maximal
Lie subalgebra $E_{7}+\mathbf{sl}\left( 2\right) $. They are 
\begin{equation*}
\begin{array}{ll}
E_{n}=E_{n-1}+\mathbf{u}\left( 1\right) +\mathbf{L}_{n-1}\otimes \eta ^{9-n}+%
\overline{\mathbf{L}_{n-1}\otimes \eta ^{9-n}} & \text{for }n\leq 7, \\ 
E_{8}=E_{7}+A_{1}+\mathbf{L}_{7}\otimes \Lambda _{1}. & 
\end{array}
\end{equation*}
Here $\eta $ is the standard representation of $\mathbf{u}\left( 1\right) $
and $\Lambda _{1}$ is the standard representation of $A_{1}=\mathbf{sl}%
\left( 2\right) $. This highly suggests that the bundle $L\mathcal{E}_{n}$
carries a fiberwise Lie algebra structure of type $E_{n}$ such that $%
\mathcal{L}_{n}$ is one of its representation bundle corresponding to the
representation $\mathbf{L}_{n}$ of $E_{n}$. There are various ways to
describe this fiberwise Lie algebra structure on $L\mathcal{E}_{n}$. A
direct approach would be to first write $L\mathcal{E}_{n}=\Lambda \otimes _{%
\mathbb{Z}}O_{X}\bigoplus_{D}O_{X}\left( D\right) $ where $\Lambda \subset
H^{2}\left( X_{n},\mathbb{Z}\right) $ is the perpendicular complement of $%
K_{X}$ with respect to the intersection pairing and the summation is over
those $D$'s satisfying $D^{2}=-2$ and $D\cdot K=0$. Second first Chern
classes of $D$'s are elements in $\Lambda $ and they form the weight lattice
of $E_{n}$ (see \cite{Manin}). This gives the bracket between $\Lambda
\otimes _{\mathbb{Z}}O_{X}$ and $O_{X}\left( D\right) $. Third the bracket
between elements in $O\left( D_{1}\right) $ and $O\left( D_{2}\right) $ can
be described using homomorphisms $O\left( D_{1}\right) \otimes O\left(
D_{2}\right) \rightarrow O\left( D_{3}\right) $ which equals zero unless $%
D_{1}\cdot D_{2}=1$ and $D_{3}=D_{1}+D_{2},$ in that case, it is an
isomorphism.

\qquad Instead of filling in the details of this algebraic approach, we
discuss several other descriptions of the fiberwise Lie algebra structure on 
$L\mathcal{E}_{n}$ using its decomposition into subbundles. It is because
these latter approaches reflect more about the geometry of lines and rulings
on $X_{n}$. One way to describe the Lie algebra structure on $E_{8}$ is to
break it down using its Lie subalgebra $D_{8}$, see for example \cite{Adams}%
. This approach to describe the fiberwise Lie algebra structure on $L%
\mathcal{E}_{8}$ will be discussed in the last section of this paper. By
restricting this bracket to $\pi ^{\ast }L\mathcal{E}_{7}$ we obtain the
fiberwise Lie algebra structure on $L\mathcal{E}_{7}$. Inductively we obtain
the fiberwise Lie algebra structure on every $L\mathcal{E}_{n}$.

\qquad Two other methods to describe the Lie algebra structure on $E_{n}$
when $n\leq 7$ is to break it down using its maximum Lie subalgebra $D_{n-1}$
and $E_{n-1}$. We discuss using $L\mathcal{E}_{n-1}$ (resp. $L\mathcal{D}%
_{n-1}$) to construct the fiberwise Lie algebra structure on $L\mathcal{E}%
_{n}$ now (resp. in next section).

\qquad First let us recall that $E_{n}=E_{n-1}+\mathbf{u}\left( 1\right) +%
\mathbf{L}_{n-1}\otimes \eta ^{9-n}+\overline{\mathbf{L}_{n-1}\otimes \eta
^{9-n}}$ when $n\leq 7$. Moreover the Lie algebra structure on $E_{n}$ can
be reconstruct from $E_{n-1}$ and $\mathbf{L}_{n-1}$ together with various
pairings between them. Likewise we can construct a Lie algebra bundle of
type $E_{n}$ on $X_{n}$ inductively from the $E_{n-1}$-bundle $L\mathcal{E}%
_{n-1}$ on $X_{n-1}\ $and the vector bundle $\mathcal{L}_{n-1}$. When $n=1$, 
$L\mathcal{E}_{1}$ is just the trivial line bundle on $X_{1}$ with the
Abelian Lie algebra structure on each fiber. $\mathcal{L}_{1}=O\left(
l\right) $ where $l$ is the unique line on $X_{1}$. Moreover the fiberwise
action of $L\mathcal{E}_{1}$ on $\mathbf{L}_{1}$ is given by the isomorphism 
$O\otimes O\left( l\right) \cong O\left( l\right) $.

\qquad By induction on $n,\mathcal{E}_{n-1}$ is an $E_{n-1}$-bundle over $%
X_{n-1}$. This determines a fiberwise Lie bracket homomorphism 
\begin{equation*}
\alpha :L\mathcal{E}_{n-1}\otimes L\mathcal{E}_{n-1}\rightarrow L\mathcal{E}%
_{n-1}.
\end{equation*}
Moreover we have a homomorphism 
\begin{equation*}
\beta :L\mathcal{E}_{n-1}\otimes \mathcal{L}_{n-1}\rightarrow \mathcal{L}%
_{n-1}
\end{equation*}
which is induced from (i) the natural homomorphism $O\left(
l_{1}-l_{2}\right) \otimes O\left( l\right) \rightarrow O\left( l_{1}\right) 
$ which is an isomorphism when $l_{2}=l$ and vanishes otherwise and (ii) the
homomorphism $\Lambda \otimes O\left( l\right) \rightarrow O\left( l\right) $
which send $D\otimes s$ to $\left( D\cdot l\right) s$. Notice that $\beta $
express $\mathcal{L}_{n-1}$ as a representation bundle of $L\mathcal{E}%
_{n-1} $. Similarly we have the dual representation. 
\begin{equation*}
\gamma :L\mathcal{E}_{n-1}\otimes \mathcal{L}_{n-1}^{\ast }\rightarrow 
\mathcal{L}_{n-1}^{\ast }
\end{equation*}

Furthermore we can define a homomorphism 
\begin{equation*}
\delta :\mathcal{L}_{n-1}\otimes \mathcal{L}_{n-1}^{\ast }\rightarrow L%
\mathcal{E}_{n-1}
\end{equation*}
which is induced from the natural homomorphism $O\left( l\right) \otimes
O\left( -l^{\prime }\right) \rightarrow O\left( l-l^{\prime }\right) $ when $%
l$ and $l^{\prime }$ are disjoint. We also have a trace homomorphism 
\begin{equation*}
t:\mathcal{L}_{n-1}\mathcal{\,\otimes L}_{n-1}^{\ast }\rightarrow O_{X}.
\end{equation*}

These homomorphisms defines a fiberwise Lie bracket on $\mathcal{E}_{n}$
using our previous decomposition as follow: 
\begin{equation*}
\begin{array}{cl}
& \left[ \left( a,b,cd,ef\right) ,\left( a^{\prime },b^{\prime },c^{\prime
}d^{\prime },e^{\prime }f^{\prime }\right) \right] \\ 
= & \left( 
\begin{array}{cccc}
\alpha \left( a,a^{\prime }\right) & t\left( c,e^{\prime }\right) \left(
d,f^{\prime }\right) & \beta \left( a,c^{\prime }\right) d^{\prime } & 
\gamma \left( a,e^{\prime }\right) f \\ 
+\delta \left( c,e^{\prime }\right) \left( d,f^{\prime }\right) & -t\left(
c^{\prime },e\right) \left( d^{\prime },f\right) , & -\beta \left( a^{\prime
},c\right) d & -\gamma \left( a^{\prime },e\right) f^{\prime } \\ 
-\delta \left( c^{\prime },e\right) \left( d^{\prime },f\right) &  & 
+c^{\prime }b^{9-n}d^{\prime } & +e^{\prime }b^{9-n}f^{\prime } \\ 
&  & -c^{\prime }b^{9-n}d, & -eb^{\prime 9-n}f
\end{array}
\right) .
\end{array}
\end{equation*}

The fact that this bracket gives $L\mathcal{E}_{n}$ a fiberwise $E_{n}$ Lie
algebra structure over $X_{n}$ follows from a corresponding construction for
the Lie algebra structure on $E_{n}$ \cite{Adams}. If we want to construct
the fiberwise Lie algebra structure on $L\mathcal{E}_{8}$ using this method,
we need to suitable modify it to adopt the $\mathbf{sl}\left( 2\right) $
subalgebra of $E_{8}$.\qquad 

As we have mentioned in this construction that $\mathcal{L}_{n}$ is the
representation bundle of the Lie algebra bundle $L\mathcal{E}_{n}$ over $%
X_{n}$ corresponding to the fundamental representation $\mathbf{L}_{n}$.
Similarly the bundle $\mathcal{R}_{n}$ is a representation bundle of $L%
\mathcal{E}_{n}$ corresponding to the fundamental representation $\mathbf{R}%
_{n}$ of $E_{n}$. To describe its action for $n\leq 6$ we use the
homomorphism $\Lambda \otimes O\left( R\right) \rightarrow O\left( R\right) $
which send $D\otimes s$ to $\left( D\cdot R\right) s$ and the fact that $%
R^{\prime }=D+R$ determines a ruling on $X_{n}$ if and only if $D\cdot R=1.$
In fact this happens precisely when we can write $D=l_{1}-l_{2}$ and $%
R=l_{2}+l_{3}$ such that $l_{1}$ and $l_{3}$ intersects at a point.

When $n=7$, we suppose that $D=l_{1}-l_{2}$ and $R=l_{2}+l_{3}$ then $l_{1}$
and $l_{3}$ can also intersect at two points. In that case they form a 2-gon
and therefore $D+R$ is linearly equivalent to $-K$. If we write 
\begin{equation*}
\mathcal{R}_{7}=\Lambda \otimes _{\mathbb{Z}}O\left( -K\right) +\tbigoplus_{R%
\text{:ruling}}O\left( R\right) .
\end{equation*}
Then in this case the action of $s\in O\left( D\right) \subset L\mathcal{E}%
_{7}$ on $t\in O\left( R\right) \subset \mathcal{R}_{7}$ is given by $\left[
D\right] \otimes \left( st\right) \in \Lambda \otimes _{\mathbb{Z}}O\left(
-K\right) \subset \mathcal{R}_{7}$. Here $\left[ D\right] \in \Lambda $
denotes the divisor class of $D$ and $st$ is an element of $O\left(
-K\right) $ given by the product of $s$ and $t$ under the isomorphism $%
O\left( D\right) \otimes O\left( R\right) \cong O\left( -K\right) $.

\bigskip

Now if $s\in O\left( D\right) \subset L\mathcal{E}_{7}$ and $\left[ D_{0}%
\right] \otimes u\in \Lambda \otimes _{\mathbb{Z}}O\left( -K\right) \subset 
\mathcal{R}_{7}$, then $R=D-K$ satisfies $R^{2}=0$ and $R\cdot K=-2$ and
therefore $R$ determines a ruling on $X_{7}$. The action of $s$ on $\left[
D_{0}\right] \otimes u$ is given by $\left( D\cdot D_{0}\right) \left(
su\right) $ where $D\cdot D_{0}$ is the intersection number of $D$ and $D_{0}
$ and $su\in O\left( R\right) $ is the product of $s$ and $u$ under the
isomorphism $O\left( D\right) \otimes O\left( -K\right) \cong O\left(
R\right) $. The action of $L\mathcal{E}_{7}$ on $\mathcal{R}_{7}$ is
completely determined by these products \cite{Adams}.\qquad 

\bigskip

We remark that the homomorphism $\mathcal{L}_{n}\otimes \mathcal{L}%
_{n}\rightarrow \mathcal{R}_{n}$ we described earlier is a homomorphism
between representation bundles over $X_{n}$ which corresponds to the Lie
algebra fact that $\mathbf{R}_{n}$ is an irreducible component of $\mathbf{L}%
_{n}\otimes \mathbf{L}_{n}$.

\subsection{Decomposing $\mathcal{L}_{n}$ under $\protect\pi ^{\ast }L%
\mathcal{E}_{n-1}$}

We fix a line $L$ on $X_{n}$ and we denote the blow down morphism of $X_{n}$
along $L$ as $\pi :X_{n}\rightarrow X_{n-1}$. We also denote the image of $L$
as $p\in X_{n-1}$. In the previous subsection we decompose $L\mathcal{E}_{n}$
under the Lie subalgebra bundle $\pi ^{\ast }L\mathcal{E}_{n-1}$. Now we
discuss a similar decomposition for $\mathcal{L}_{n}$. Notice that choosing
a fix line on $X_{n}$ break the symmetry among all lines on $X_{n}$ and this
decomposition of $\mathcal{L}_{n}$ reflects the intersection properties of $%
L $ with other lines on $X_{n}$.

\bigskip

First we discuss the decomposition of the representation $\mathbf{L}_{n}$
under the Lie subalgebra $E_{n-1}$, the branching rule. As a representation
of $E_{n-1}$, we have $\mathbf{L}_{n}=\mathbf{L}_{n-1}+\mathbf{R}_{n-1}+%
\mathbf{1}$ when $n\leq 6.$ Here $\mathbf{1}$ denote the trivial
representation of dimension one. We also have $\mathbf{L}_{7}=\mathbf{L}_{6}+%
\mathbf{R}_{6}+\mathbf{1+1}$. When $n=8$ then $E_{7}+A_{1}$ is a maximal Lie
subalgebra of $E_{8}$ and as a representation of this subalgebra we have the
decomposition $\mathbf{L}_{8}=\mathbf{L}_{7}\otimes \Lambda _{1}+\mathbf{R}%
_{7}+A_{1}$ where $\Lambda _{1}$ is the standard representation of $A_{1}=%
\mathbf{sl}\left( 2\right) $.

\bigskip

First we notice that $\pi ^{\ast }\mathcal{L}_{n-1}$ is a subbundle of $%
\mathcal{L}_{n}$. This is equivalent to the statement that $p$ does not lie
on any line on $X_{n-1}$. Otherwise the strict transform of the line
containing $p$ would be a $\left( -2\right) $-curve and the existence of
such curve violates the ampleness assumption of the anti-canonical class of $%
X_{n}.$ Second, it is obvious that $O\left( L\right) $ is also a subbundle
of $\mathcal{L}_{n}$.

\bigskip

Third if $R$ determines a ruling on $X_{n-1}$, then there is a unique fiber
of the ruling which contains the point $p$. By similar reason as above, this
must be a smooth fiber of the ruling. The strict transform of this smooth
fiber containing $p$ is a line on $X_{n}$. Hence there is a line in the
divisor class $\pi ^{\ast }R-L$. Now we recall that the representation
bundle of $L\mathcal{E}_{n-1}$ corresponding to $\mathbf{R}_{n-1}$ is given
by 
\begin{eqnarray*}
\mathcal{R}_{n-1} &=&\bigoplus_{R\text{:ruling}}O_{X_{n-1}}\left( R\right) \,%
\text{when }n\leq 7 \\
\mathcal{R}_{7\,\,} &=&\bigoplus_{R\text{:ruling}}O_{X_{7}}\left( R\right)
+O\left( -K\right) ^{\oplus 7}.
\end{eqnarray*}

Therefore, when $n\leq 7$, then $\pi ^{\ast }\mathcal{L}_{n-1}+\pi ^{\ast }%
\mathcal{R}_{n-1}\otimes O\left( -L\right) +O\left( L\right) $is a subbundle
of $\mathcal{L}_{n}$. In fact if $n\leq 6$ this exhaust the whole bundle $%
\mathcal{L}_{n}$. This is because $l\cdot L=-1,0$ or $1$ for any line $l$ on 
$X_{n}$. When $l\cdot L=-1$ (resp. $0$ and $1$), then $l$ is $L$ (resp.
constructed using $\mathcal{L}_{n-1}$ and $\mathcal{R}_{n-1}$). When $n=7$,
there is a unique line $l$ on $X_{7}$ such that $l\cdot L=2$. Namely $l$ and 
$L$ forms a 2-gon. That is, $l$ is linearly equivalent to $-K-L$. When $n=8$%
, the situation is slightly more involved and we will discuss it later. In
conclusion we have 
\begin{equation*}
\begin{array}{llc}
\mathcal{L}_{n}=\pi ^{\ast }\mathcal{L}_{n-1} & +\pi ^{\ast }\mathcal{R}%
_{n-1}\otimes O\left( -L\right) +O\left( L\right) & \text{when }n\leq 6, \\ 
\mathcal{L}_{7}=\pi ^{\ast }\mathcal{L}_{6} & +\pi ^{\ast }\mathcal{R}%
_{6}\otimes O\left( -L\right) +O\left( L\right) +O\left( -K-L\right) \text{,}
&  \\ 
\mathcal{L}_{8}=\pi ^{\ast }\mathcal{L}_{7}\otimes \Lambda _{1}^{\ast } & 
+\pi ^{\ast }\mathcal{R}_{7}\otimes O\left( -L\right) +\mathcal{A}%
_{1}\otimes O\left( -K\right) . & 
\end{array}
\end{equation*}

Moreover these decompositions of vector bundles are in fact decompositions
of representation bundles of $L\mathcal{E}_{n-1}$.

\bigskip

For example the ten lines on a cubic surface $X_{6}$ intersecting $L$ forms
five pairs of intersecting lines (see e.g. \cite{Hartshorne} and \cite{Reid2}%
). It is reflecting the fact that the representation $\mathbf{R}_{5}$ of $%
E_{5}$ is the standard representation of $\mathbf{so}\left( 10\right) $
under the identification of $E_{5}$ with $D_{5}=\mathbf{so}\left( 10\right) $
and therefore it carries a natural quadratic form. See section six for
details. Moreover each line on the cubic surface determines a unique ruling
on it and vice versa. More precisely we have $\mathcal{R}_{6}=\mathcal{L}%
_{6}^{\ast }\otimes O\left( -K_{X_{6}}\right) $. This isomorphism
corresponds to the isomorphism $\mathbf{R}_{6}\cong \mathbf{L}_{6}^{\ast }$
via the outer-automorphism of the Dynkin diagram of $E_{6}$. Then we can
rewrite the above decomposition of $\mathcal{L}_{7}$ as 
\begin{equation*}
\mathcal{L}_{7}=\left( \pi ^{\ast }\mathcal{L}_{6}+O\left( -L\right) \right)
\otimes \left( O+O\left( -K\right) \right) .
\end{equation*}

To see this geometrically, we consider a projection of the cubic surface $%
X_{6}$ away from a point $q$ on it. Then the projection defines a double
cover of $\mathbb{P}^{2}$ branched along a quartic plane curve $C$. The 27
lines on $X_{6}$ together with $L$ on $X_{7}$ project to the 28 bitangents
to the quartic plane curve. The pre-image of each bitangent to $C$ is a
2-gon on $X_{7}$. These informations are neatly packaged in the above
decomposition of $\mathcal{L}_{7}$.

\bigskip

Next we want to verify that the Lie algebra bundle associated to each of the 
$E_{n}$-bundle $\mathcal{E}_{n}$ on $X_{n}$ is precisely $L\mathcal{E}_{n}$.
This is because each of those $\mathcal{E}_{n}$'s are defined as
automorphism bundles of certain fiberwise algebraic structure on
representation bundles of $L\mathcal{E}_{n}$. Moreover such algebraic
structures are also preserved by $L\mathcal{E}_{n}$ and therefore the
identification of $Lie\left( \mathcal{E}_{n}\right) $ and $L\mathcal{E}_{n}$
follows from corresponding statements in representation theory of Lie
algebra (see e.g. \cite{Adams}).

\subsection{Small $n$ examples}

When $n=1$ then $X_{1}$ is just the blowup of $\mathbb{P}^{2}$ at one point
and there is a unique line $l$ on $X_{1}$, namely the exceptional curve for
the blowup morphism. Also $X_{1}$ has a unique ruling determined by $H-l$
where $H$ is the pullback of the hyperplane class from $\mathbb{P}^{2}$. In
this case we have $L\mathcal{E}_{1}=O$, $\mathcal{L}_{1}=O\left( l\right) $
and $\mathcal{R}_{1}=O\left( H-l\right) $.

\bigskip

When $n=2$ there are 3 lines on $X_{2}$. We denote them as $l_{1}$, $l_{2}$
and $l_{3}$. They are divided into two different types. By rearranging
indexes, we can assume that $l_{1}\cdot l_{3}=l_{2}\cdot l_{3}=1$ and $%
l_{1}\cdot l_{2}=0$. Notice that this is the only time that lines on $X_{n}$
are divided into different types. This can translated into reducibility of $%
\mathcal{L}_{2}$ as a representation bundle of $L\mathcal{E}_{2}$.

\bigskip

Recall that $E_{2}=A_{1}+\mathbf{u}\left( 1\right) $ and we can write $L%
\mathcal{E}_{2}=L\mathcal{A}_{1}+O$ where $L\mathcal{A}_{1}=O+O\left(
l_{1}-l_{2}\right) +O\left( l_{2}-l_{1}\right) $ is the Lie algebra bundle
of type $A_{1}$. From definitions we have $\mathcal{L}_{2}=O\left(
l_{1}\right) +O\left( l_{2}\right) +O\left( l_{3}\right) $ and $\mathcal{R}%
_{2}=O\left( l_{1}+l_{3}\right) +O\left( l_{2}+l_{3}\right) $. As
representations of $A_{1}$ we have $\mathbf{L}_{2}=\mathbf{R}_{2}+\mathbf{1}$%
. Similarly, as representation bundles of $L\mathcal{A}_{1}$ we have $%
\mathcal{L}_{2}=\mathcal{R}_{2}\otimes O\left( -l_{3}\right) +O\left(
l_{3}\right) $.

\bigskip

If we blow down $X_{2}$ along $l_{1}$ (or $l_{2}$), then we get $X_{1}$ and
the $E_{1}$-bundle $L\mathcal{E}_{1}$ over it with $E_{1}=\mathbf{u}\left(
1\right) $ as before. However we can also blow down $X_{2}$ along $l_{3}$,
then we get the surface $\mathbb{P}^{1}\times \mathbb{P}^{1}$. Now we should
interpret $E_{1}$ as the Lie subalgebra $A_{1}$ of $E_{2}$ instead of $%
\mathbf{u}\left( 1\right) $. To prevent confusion, we call it $\overline{E}%
_{1}$ and we also write $\overline{\mathbf{L}}_{1}$ , $\overline{\mathbf{R}}%
_{1}$ for its representations. In fact $\overline{\mathbf{L}}_{1}$ is
trivial and $\overline{\mathbf{R}}_{1}$ is the standard representation of it
under the identification of $\overline{E}_{1}$ with $A_{1}=\mathbf{sl}\left(
2\right) $.

\bigskip

As for $\mathbb{P}^{1}\times \mathbb{P}^{1}$, it has no lines and two
rulings. Therefore $\overline{\mathcal{L}}_{1}$ is trivial and $\overline{%
\mathcal{R}}_{1}$ is a rank two bundle over $\mathbb{P}^{1}\times \mathbb{P}%
^{1}$. Moreover $L\overline{\mathcal{E}}_{1}=End_{0}\left( \overline{%
\mathcal{R}}_{1}\right) $. The previous decomposition of $L\mathcal{E}_{n}$
and $\mathcal{L}_{n}$ under $L\mathcal{E}_{n-1}$ also hold true in this
case. Since $\overline{\mathcal{L}}_{1}$ is trivial we have 
\begin{eqnarray*}
L\mathcal{E}_{2} &=&\pi ^{\ast }L\overline{\mathcal{E}}_{1}+O, \\
\mathcal{L}_{2} &=&\pi ^{\ast }\overline{\mathcal{R}}_{1}\otimes O\left(
-L\right) +O\left( L\right) .
\end{eqnarray*}
Here $L=l_{3}$. We can also write down the decomposition for $\mathcal{R}%
_{2} $, namely $\mathcal{R}_{2}=\pi ^{\ast }\overline{\mathcal{R}}_{1}$.

\bigskip

When $n=3$ the non-simpleness of $E_{3}$ is related to Cremona
transformations on $X_{3}$. First we represent $X_{3}$ as a blowup of $%
\mathbb{P}^{2}$ at three generic points, $\pi :X_{3}\rightarrow $ $\mathbb{P}%
^{2}$ and denote the exceptional locus as $L_{1}\cup L_{2}\cup L_{3}$.
Notice that the Lie algebra $E_{3}$ is not simple: $E_{3}=\mathbf{sl}\left(
3\right) +\mathbf{sl}\left( 2\right) $ and $\mathbf{L}_{3}$ is the tensor
product of standard representation of $\mathbf{sl}\left( 3\right) $ and $%
\mathbf{sl}\left( 2\right) $. Geometrically we have $L\mathcal{E}_{3}=%
\mathcal{A}_{3}+\mathcal{A}_{2}$ where $\mathcal{A}_{3}=2O\oplus _{i\neq
j}O\left( L_{i}-L_{j}\right) $ is a $\mathbf{sl}\left( 3\right) $ Lie
algebra bundle over $X_{3}$ and $\mathcal{A}_{2}=O+O\left(
H-L_{1}-L_{2}-L_{3}\right) +O\left( -H+L_{1}+L_{2}+L_{3}\right) $ is a $%
\mathbf{sl}\left( 2\right) $ Lie algebra bundle over $X_{3}$. For $\mathcal{L%
}_{3}$ we also have $\mathcal{L}_{3}=\mathcal{W}_{3}\otimes \mathcal{W}_{2}$
with $\mathcal{W}_{3}=O\left( L_{1}\right) +O\left( L_{2}\right) +O\left(
L_{3}\right) $ and $\mathcal{W}_{2}=O+O\left( H-L_{1}-L_{2}-L_{3}\right) $
being representation bundles of $\mathcal{A}_{3}$ and $\mathcal{A}_{2}$
respectively.

\qquad

We can also write $\mathcal{L}_{3}=\mathcal{W}_{3}^{\prime }\otimes \mathcal{%
W}_{2}^{\ast }$ with $\mathcal{W}_{3}^{\prime }=O\left( H-L_{2}-L_{3}\right)
+O\left( H-L_{3}-L_{1}\right) +O\left( H-L_{1}-L_{2}\right) $. The symmetric
relation between $\mathcal{W}_{3}$ and $\mathcal{W}_{3}^{\prime }$ is
related to the Cremona transformation.

$\qquad $

When $n=4$ the ranks of $\mathcal{L}_{4}$ and $\mathcal{R}_{4}$ equals ten
and five respectively. Notice that $E_{4}=A_{4}=\mathbf{sl}\left( 5\right) $
and the representation $\mathbf{R}_{4}$ of $E_{4}$ corresponds to the
standard representation of $\mathbf{sl}\left( 5\right) $. Therefore we
expect that $L\mathcal{E}_{4}=End_{0}\left( \mathcal{R}_{4}\right) $, the
bundle of trace-free endomorphisms of $\mathcal{R}_{4}$. To see this
directly, we represent $X_{4}$ as a blowup of $\mathbb{P}^{4}$ at four
points and we denote the exceptional locus of this blowup as $L_{1}\cup
L_{2}\cup L_{3}\cup L_{4}$, then we have 
\begin{equation*}
\mathcal{L}_{4}=\oplus _{i=1}^{4}O\left( L_{i}\right) +\oplus _{1\leq
i<j\leq 4}O\left( H-L_{i}-L_{j}\right)
\end{equation*}
and 
\begin{equation*}
\mathcal{R}_{4}=\oplus _{i=1}^{4}O\left( H-L_{i}\right) +O\left(
2H-L_{1}-L_{2}-L_{3}-L_{4}\right) .
\end{equation*}
By direct computations, we have 
\begin{align*}
End_{0}\left( \mathcal{R}_{4}\right) & =O_{X}^{\oplus 4}+\oplus _{i\neq
j}\left( L_{i}-L_{j}\right) +\oplus _{i=1}^{4}\left( O\left(
H-\tsum_{j=1}^{4}L_{j}+L_{i}\right) +O\left(
-H+\tsum_{j=1}^{4}L_{j}-L_{i}\right) \right) \\
& =O_{X}^{\oplus 4}+\bigoplus_{\substack{ D=l-l^{\prime }  \\ l,l^{\prime }%
\text{disjoint lines on }X_{4}}}O\left( D\right) \\
& =O_{X}^{\oplus 4}+\bigoplus_{\substack{ D^{2}=-2  \\ D\cdot K=0}}O\left(
D\right) .
\end{align*}

We also have the isomorphism 
\begin{equation*}
\Lambda ^{3}\mathcal{R}_{4}=\mathcal{L}_{4}\otimes O\left( -K\right) .
\end{equation*}
This corresponds to the fact that the fundamental representation $\mathbf{L}%
_{4}$ of $\mathbf{sl}\left( 5\right) $ equals the third exterior power of
its standard representation.

\newpage

\section{$D_{n-1}$-bundle over ruled surfaces}

\subsection{Ruling on $X_{n}$\ and construction of $\mathcal{D}_{n-1}$}

The $E_{n}$-bundle $\mathcal{E}_{n}$ over $X_{n}$ would have its structure
group reduces to a smaller subgroup of $E_{n}$ if $X_{n}$ admits additional
geometric structure. In this section, we study the subalgebra $\mathbf{so}%
\left( 2n-2\right) =D_{n-1}\subset E_{n}$. This corresponds to removing the
node on the right end in the Dynkin diagram of $E_{n}$. 
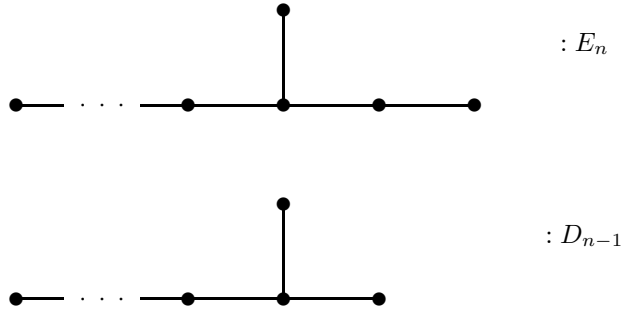
\begin{equation*}
\end{equation*}

\begin{picture}(3.4,1)(1.4,.4)
\thicklines
\put(4.3,1){\circle*{.075}}
\put(3.8,1){\circle*{.075}}
\put(3.8,1){\line(1,0){.5}}
\put(3.3,1){\circle*{.075}}
\put(3.3,1){\line(0,1){.5}}
\put(3.3,1.5){\circle*{.075}}
\put(3.3,1){\line(1,0){.5}}
\put(2.8,1){\circle*{.075}}
\put(2.8,1){\line(1,0){.5}}
\put(2.55,1){\line(1,0){.25}}
\put(2.45,1){\circle*{.02}}
\put(2.35,1){\circle*{.02}}
\put(2.25,1){\circle*{.02}}
\put(1.9,1){\circle*{.075}}
\put(1.9,1){\line(1,0){.25}}
\put(4.75,1.2){\makebox(.25,.25){$: E _ n$}}
\end{picture}

\begin{picture}(3.4,1)(1.4,.4)
\thicklines
\put(3.8,1){\circle*{.075}}
\put(3.3,1){\circle*{.075}}
\put(3.3,1){\line(0,1){.5}}
\put(3.3,1.5){\circle*{.075}}
\put(3.3,1){\line(1,0){.5}}
\put(2.8,1){\circle*{.075}}
\put(2.8,1){\line(1,0){.5}}
\put(2.55,1){\line(1,0){.25}}
\put(2.45,1){\circle*{.02}}
\put(2.35,1){\circle*{.02}}
\put(2.25,1){\circle*{.02}}
\put(1.9,1){\circle*{.075}}
\put(1.9,1){\line(1,0){.25}}
\put(4.75,1.2){\makebox(.25,.25){$: D _ {n-1}$}}
\end{picture}

Recall that this right end node corresponds to the fundamental
representation $\mathbf{R}_{n}$ of $E_{n}$. Its associated representation
bundle $\mathcal{R}_{n}$ is constructed using the set of rulings on $X_{n}.$
Now we fix one ruling on $X_{n}$. That is we pick a divisor $R$ on $X_{n}$
which satisfies $R^{2}=0$ and $R\cdot K_{X_{n}}=-2.$ We shall see that such
a geometric structure reduces the structure group of $\mathcal{E}_{n}$ from $%
E_{n}$ to $D_{n-1}$. We remark that the construction of the $D_{n-1}$-bundle
over $X_{n}$ applies to blowups of $\mathbb{P}^{2}$ at arbitrary number of
points.

Let $X=X_{n}$ be a blowup of $\mathbb{P}^{2}$ at $n$ generic points and let $%
R$ be a divisor on $X_{n}$ defining a ruling .We define a vector bundle on $%
X_{n}$ as follow: 
\begin{equation*}
\mathcal{W}_{n-1}=\bigoplus_{\substack{ C^{2}=-1  \\ CK=-1  \\ CR=0}}%
O_{X}\left( C\right) .
\end{equation*}
Then $\mathcal{W}_{n-1}$ is a subbundle of $\mathcal{L}_{n}$ for $n\leq 8$.
In general the rank of $\mathcal{W}_{n-1}$ equals $2n-2$ and it carries a
natural fiberwise non-degenerate quadratic form $q:\mathcal{W}_{n-1}\otimes 
\mathcal{W}_{n-1}\rightarrow O_{X_{n}}\left( R\right) $. We denote the
connected component of the automorphism bundle of $\mathcal{W}_{n-1}$
preserving $q$ as $\mathcal{D}_{n-1}$.

\bigskip

To see this we need to look at the ruling $\pi :X_{n}\rightarrow \mathbb{P}%
^{1}$ determined by $R$. Namely the fiber divisor class of $\pi $ equals $R$%
. Genericity of $X_{n}$ implies that each singular fiber of $\pi $ consists
of two irreducible components. Each of them is a (-1) curve and they meet at
a single point transversely. 

Since the Euler number of $X_{n}$ equals $n+3$, the number of singular
fibers is $n-1.$ We denote them as $L_{1}\cup L_{1}^{\prime },L_{2}\cup
L_{2}^{\prime },...,L_{n-1}\cup L_{n-1}^{\prime }$. For any $C=L_{i}$ or $%
L_{i}^{\prime }$, it satisfies $C^{2}=-1,C\cdot K_{X}=-1$ and $C\cdot R=0$.
It is not difficult to show that the converse is also true. Hence $\mathcal{W%
}_{n-1}=\oplus _{i=1}^{n-1}\left( O_{X_{n}}\left( L_{i}\right)
+O_{X_{n}}\left( L_{i}^{\prime }\right) \right) $. Using the isomorphisms $%
O_{X_{n}}\left( L_{i}\right) \otimes O_{X_{n}}\left( L_{i}^{\prime }\right)
\cong O_{X_{n}}\left( R\right) $, we obtain a non-degenerated fiberwise
quadratic form 
\begin{equation*}
q_{n-1}:\mathcal{W}_{n-1}\otimes \mathcal{W}_{n-1}\rightarrow
O_{X_{n}}\left( R\right) .
\end{equation*}
In terms of the above decomposition of $\mathcal{W}_{n-1}$, we have $%
q_{n-1}=\oplus ^{n-1}\mathbf{H}$ with $\mathbf{H}=\left( 
\begin{array}{cc}
0 & 1 \\ 
1 & 0
\end{array}
\right) .$

It is easy to see that the adjoint bundle of $\mathcal{D}_{n-1}$ is given by 
\begin{equation*}
L\mathcal{D}_{n-1}=\Lambda ^{2}\mathcal{W}_{n-1}\otimes O_{X_{n}}\left(
-R\right) .
\end{equation*}
Using the identities $O\left( L_{i}+L_{j}-R\right) \cong O\left(
L_{i}-L_{j}^{\prime }\right) $ and $O\left( L_{i}+L_{j}^{\prime }-R\right)
\cong O\left( L_{i}-L_{j}\right) $, we obtain 
\begin{equation*}
L\mathcal{D}_{n-1}=O_{X}^{\oplus n-1}\bigoplus_{\substack{ C^{2}=-2  \\ CK=0 
\\ CR=0}}O_{X}\left( C\right) .\footnote{%
Such classes $C$ can be characterized as those that can be writeen as the
difference of two disjoint lines, each has zero intersection with $D$.}
\end{equation*}
From this description, it is obvious that $L\mathcal{D}_{n-1}$ is a vector
subbundle of $L\mathcal{E}_{n}$ when $n\leq 8$. In fact $L\mathcal{D}_{n-1}$
is a Lie subalgebra bundle of $L\mathcal{E}_{n}$.

\subsection{Spinor bundles of $\mathcal{D}_{n-1}$}

Using $\Lambda ^{l}\mathcal{W}_{n-1}$ with $l=1,2,...,n-3,$ we can construct
associated bundles of $\mathcal{D}_{n-1}$ corresponding to every fundamental
representation of $D_{n-1}$ except the two spinor representations $S^{+}$
and $S^{-}$. In fact, these two spinor bundles $\mathcal{S}^{+}$ and $%
\mathcal{S}^{-}$ over $X_{n}$ will be needed later when we decompose $L%
\mathcal{E}_{n}$ and $\mathcal{L}_{n}$ under restriction to $L\mathcal{D}%
_{n-1}$. We leave them as exercise for readers to check directly that 
\begin{equation*}
\mathcal{S}^{+}=\bigoplus_{\substack{ S^{2}=-1  \\ SK=-1  \\ SR=1}}%
O_{X_{n}}\left( S\right) \text{ and }\mathcal{S}^{-}=\bigoplus_{\substack{ %
T^{2}=-2  \\ TK=0  \\ TR=1}}O_{X_{n}}\left( T\right) .
\end{equation*}
In particular divisor $S$ as above correspond to line on $X$ which is a
section of the ruling $R$. Similar for any divisor $T$ as above, the divisor 
$T+R$ corresponds to a ruling on $X$ which is a section of the given ruling $%
R$. In fact $n\leq 5$, then every line on $X_{n}$ is either a section of the
ruling or is a component of a singular fiber of the ruling. This gives rise
to the decomposition $\mathcal{L}_{n}=\mathcal{W}_{n-1}+\mathcal{S}^{+},$
\thinspace \thinspace \thinspace for $n\leq 5$. On the other hand, adding to 
$R$ a divisor $T$ as above gives us another ruling of $X_{n}$. When $n$ $%
\leq 4$, then determines all rulings on $X_{n}$. This gives rise to the
decomposition $\mathcal{R}_{n}=O\left( R\right) \otimes \left( O+\mathcal{S}%
^{-}\right) ,$ \thinspace \thinspace \thinspace for $n\leq 4$.

\begin{example}
$\left[ D_{3}\text{-bundle over }X_{4}\right] $ case: Let $X=X_{4}$ be the
blowup of $\mathbb{P}^{2}$ at four generic points. Recall that, under $%
Aut\left( \mathbb{P}^{2}\right) =\mathbb{P}GL\left( 3,\mathbb{C}\right) $,
any four generic points on $\mathbb{P}^{2}$ is equivalent to any others. Let
us denote the exceptional locus of this blowup as $L_{0},L_{1},L_{2}$ and $%
L_{3}$. We also denote the pullback of the hyperplane class from $\mathbb{P}%
^{2}$ as $H$. Then $R=H-L_{0}$ satisfies $R^{2}=0$ and $R\cdot K_{X}=-2$ and
hence defines a ruling on $X$.

In this ruling, there are three singular fibers $L_{1}\cup L_{1}^{\prime
},L_{2}\cup L_{2}^{\prime }$ and $L_{3}\cup L_{3}^{\prime }$. The linear
equivalent classes of $L_{i}^{\prime }$ are given by $L_{i}^{\prime
}=H-L_{0}-L_{i}$ for $i=1,2,3$. Moreover $L_{0}$ is a section of the ruling.

\begin{tabular}{|c|c|c|c|}
\hline
$\text{Bundle over }X_{4}:$ & $\mathcal{W}_{3}$ & $\mathcal{S}^{+}$ & $%
\mathcal{S}^{-}$ \\ \hline
$
\begin{array}{c}
\text{divisors for } \\ 
\text{its direct } \\ 
\text{summands}
\end{array}
$ & $
\begin{array}{c}
L_{1} \\ 
L_{2} \\ 
L_{3} \\ 
H-L_{0}-L_{1} \\ 
H-L_{0}-L_{2} \\ 
H-L_{0}-L_{3}
\end{array}
$ & $
\begin{array}{c}
L_{0} \\ 
H-L_{1}-L_{2} \\ 
H-L_{2}-L_{3} \\ 
H-L_{3}-L_{1}
\end{array}
$ & $
\begin{array}{c}
L_{0}-L_{1} \\ 
L_{0}-L_{2} \\ 
L_{0}-L_{3} \\ 
H-L_{1}-L_{2}-L_{3}.
\end{array}
$ \\ \hline
\end{tabular}
\end{example}

\begin{example}
$\left[ D_{4}\text{-bundle over }X_{5}\right] $ case: We use notations
similar to the previous example and we have the following table:

\begin{tabular}{|c|c|c|c|}
\hline
$\text{Bundle over }X_{5}:$ & $\mathcal{W}_{4}$ & $\mathcal{S}^{+}$ & $%
\mathcal{S}^{-}$ \\ \hline
$
\begin{array}{c}
\text{divisors for } \\ 
\text{its direct } \\ 
\text{summands}
\end{array}
$ & $
\begin{array}{c}
L_{1} \\ 
L_{2} \\ 
L_{3} \\ 
L_{4} \\ 
H-L_{0}-L_{1} \\ 
H-L_{0}-L_{2} \\ 
H-L_{0}-L_{3} \\ 
H-L_{0}-L_{4}
\end{array}
$ & $
\begin{array}{c}
L_{0} \\ 
H-L_{1}-L_{2} \\ 
H-L_{1}-L_{3} \\ 
H-L_{1}-L_{4} \\ 
H-L_{2}-L_{3} \\ 
H-L_{2}-L_{4} \\ 
H-L_{3}-L_{4} \\ 
2H-\tsum_{i=0}^{4}L_{i}
\end{array}
$ & $
\begin{array}{c}
L_{0}-L_{1} \\ 
L_{0}-L_{2} \\ 
L_{0}-L_{3} \\ 
L_{0}-L_{4} \\ 
H-L_{2}-L_{3}-L_{4} \\ 
H-L_{1}-L_{3}-L_{4} \\ 
H-L_{1}-L_{2}-L_{4} \\ 
H-L_{1}-L_{2}-L_{3}.
\end{array}
$ \\ \hline
\end{tabular}
\end{example}

In both these examples, we see that the union of the first two columns
exhausts all lines on $X_{4}$ and $X_{5}$. As we will explain, the
situations are more complicated for $X_{6},X_{7}$ and $X_{8}$.

\begin{center}
\textbf{Clifford multiplication}
\end{center}

In linear algebra the product of a vector and a spinor produces a spinor
with a different parity. For our bundles of spinors over $X_{n}$, we have
homomorphisms 
\begin{align*}
\mathcal{S}^{+}\otimes \mathcal{W}_{n-1}^{\ast }& \rightarrow \mathcal{S}%
^{-}, \\
\mathcal{S}^{-}\otimes \mathcal{W}_{n-1}& \rightarrow \mathcal{S}^{+},
\end{align*}
corresponding to fiberwise Clifford multiplications\footnote{$\mathcal{W}%
_{n-1}=\mathcal{W}_{n-1}^{\ast }\otimes O_{X_{n}}\left( R\right) .$}. It is
very easy to describe these homomorphisms explicitly. Let $O_{X}\left(
S\right) $ (and $O_{X}\left( -C\right) $) be a direct summand of $\mathcal{S}%
^{+}$ (and $\mathcal{W}_{n-1}^{\ast }$). If $S\cdot C=0$ then $T=S-C$
satisfies $T^{2}=-2,T\cdot K_{X}=0$ and $T\cdot R=1$. That is $O_{X}\left(
S\right) \otimes O_{X}\left( C\right) ^{\ast }=O_{X}\left( T\right) \subset 
\mathcal{S}^{-}$. We set the product of elements from other types of direct
summands to be zero. This gives the first Clifford multiplication
homomorphism. The second homomorphism can be described in the same way.

Geometrically these homomorphisms reflect the fact that an irreducible
component of a singular fiber of the given ruling together with a line
section passing through it determines another ruling which is a section of
the given ruling. Moreover all ruling sections arise in this manner.

\begin{center}
\textbf{Dualities}
\end{center}

When $n=2m$ is even, we have an isomorphism 
\begin{equation*}
\left( \mathcal{S}^{+}\right) ^{\ast }\otimes O_{X_{2m}}\left( \left(
m-4\right) R-K_{X}\right) \cong \mathcal{S}^{-}.
\end{equation*}
This follows from earlier explicit descriptions of $\mathcal{S}^{+}$ and $%
\mathcal{S}^{-}$. This isomorphism corresponds to the linear algebra fact
that the two spinor representations of $Spin\left( 4m-2\right) $ are dual to
each other.

Similarly when $n=2m+1$ is odd, we have isomorphisms 
\begin{equation*}
\left( \mathcal{S}^{+}\right) ^{\ast }\otimes O_{X_{2m+1}}\left( \left(
m-3\right) R-K_{X}\right) \cong \mathcal{S}^{+},
\end{equation*}
and 
\begin{equation*}
\left( \mathcal{S}^{-}\right) ^{\ast }\otimes O_{X_{2m+1}}\left( \left(
m-4\right) R-K_{X}\right) \cong \mathcal{S}^{-}.
\end{equation*}
These isomorphisms correspond to linear algebra facts that the two spinor
representations of $Spin\left( 4m\right) $ are self-dual. If we rewrite
these two homomorphisms as 
\begin{equation*}
\mathcal{S}^{+}\otimes \mathcal{S}^{+}\rightarrow O_{X_{2m+1}}\left( \left(
m-3\right) R-K_{X}\right) ,
\end{equation*}
and 
\begin{equation*}
\mathcal{S}^{-}\otimes \mathcal{S}^{-}\rightarrow O_{X_{2m+1}}\left( \left(
m-4\right) R-K_{X}\right) .
\end{equation*}
We obtain fiberwise non-degenerate quadratic forms on $\mathcal{S}^{+}$ and $%
\mathcal{S}^{-}$. On the one hand, it is related to the fact that spinor
representations $\mathcal{S}^{\pm }$ of $Spin\left( 4m\right) $ are of real
type. On the other hand, they give inner product structure on $\mathcal{S}%
^{+}$ and $\mathcal{S}^{-}$ which will become important later. There are
many other algebraic structures on spinor representations of $D_{n-1}$ which
can be translated to our bundle version in a similar manner. Instead of
writing out each of them, we are going to look at relationship between $%
\mathcal{E}_{n}$ and $\mathcal{D}_{n-1}$. Such relationship is very
interesting because they reflect the configuration of lines on $X_{n}$ in
relation to a fixed ruling on $X_{n}$.

\subsection{Reduction from $E_{n}$ to $D_{n-1}$}

\begin{center}
\textbf{Decompose }$\mathcal{L}_{n}\mathbf{\ }$\textbf{under }$\mathcal{D}%
_{n-1}$
\end{center}

Recall that the bundle $\mathcal{L}_{n}$ over $X_{n}$ is constructed from
the collection of lines on $X_{n}$. When $R$ determines a ruling on $X_{n}$,
every line would have non-negative intersection with $R$. However this
intersection number cannot be too big.

\begin{proposition}
If $L$ is a line on $X_{n}$ with a ruling determined by $R$, then $L\cdot
R\leq 1$ if $n\leq 5$
\end{proposition}

Proof of proposition: When $n\leq 6$ any two lines of $X_{n}$ intersect at
not more than one point. Now, from the ampleness assumption of $K^{-1}$,
every singular fiber of the ruling consists of two lines. Therefore $L\cdot
R $ is at most two. However if $L\cdot R=2$ then the union of $L$ with any
singular fiber of the ruling will form a triangle. But we know from
proposition \ref{Prop: d-gon}\ that triangle does not exist when $n\leq 5$.
Hence $L\cdot R\leq 1$. $\square $

\begin{remark}
From the proof of the proposition, we can also see that in the case of $X_{6}
$ there is only one $L$ with $L\cdot R=2$, otherwise $L\cdot R\leq 1$. It is
because the union of $L$ with any singular fiber of the ruling forms a
triangle which is necessary an anti-canonical divisor. Hence $L=-K-R$.

In fact, with more care, we can show that $L\cdot R\leq 2$ even when $n=7$.
\end{remark}

From the proposition, we have $L\cdot R=0$ or $1$ for any line $L$ on $X_{n}$
with $n\leq 5$. Line bundles associated to these divisors are precisely
direct summands of $\mathcal{W}_{n-1}$ and $\mathcal{S}^{+}$ respectively.
Therefore we have isomorphism 
\begin{equation*}
\mathcal{L}_{n}=\mathcal{W}_{n-1}+\mathcal{S}^{+}\text{ when }n\leq 5\text{.}
\end{equation*}
In fact these are isomorphisms as representation bundles of $\mathcal{D}%
_{n-1}$. Of course, these isomorphisms corresponds to the decomposition into
irreducible summands of $\mathbf{L}_{n}$ when we restrict the Lie algebra
from $E_{n}$ to $D_{n}:$%
\begin{equation*}
\mathbf{L}_{n}|_{D_{n-1}}=W_{n-1}+S^{+}
\end{equation*}
for $n\leq 5$.

\begin{example}[$D_{3}$-bundle over $X_{4}$]
In this case, after choosing a ruling $R$, we have $\mathcal{L}_{4}=\mathcal{%
W}_{3}+\mathcal{S}^{+}$. The four lines corresponding to $\mathcal{S}^{+}$
are sections of the ruling. We call them $L_{1},L_{2},L_{3}$ and $L_{4}$.
They are disjoint lines on $X_{4}$. It is because every line on $X_{4}$
intersects three other lines (from section two that $\mathcal{L}_{4}=%
\mathcal{L}_{3}+\mathcal{R}_{3}+L$ and $\dim \mathbf{R}_{3}=3$) and when it
is a section, it intersects the three singular fibers of the ruling. Each
singular fiber consists of two lines and therefore all sections are
disjoint. We can blowdown these $L_{i}$'s and we obtain $\mathbb{P}^{2}$.
Then the ruling $R$ equals $2H-L_{1}-L_{2}-L_{3}-L_{4}$ where $H=-K-R$ is
the pullback of the hyperplane class on $\mathbb{P}^{2}$. If we denote $%
L_{ij}$ the line on $X_{4}$ in the class $H-L_{i}-L_{j}$, then the three
singular fibers of the rulings are $L_{12}\cup L_{34},L_{13}\cup L_{42}$ and 
$L_{14}\cup L_{23}$. Moreover $L_{ij}$ meets $L_{k}$ if and only if $k=i$ or 
$j$.

In terms of $L\mathcal{D}_{3}$ representation bundles over $X_{4}$, we
obtain $\mathcal{W}_{3}=\Lambda ^{2}\mathcal{S}^{+}\otimes O\left(
K+2R\right) $ and $\mathcal{L}_{4}=\mathcal{S}^{+}+\Lambda ^{2}\mathcal{S}%
^{+}\otimes O\left( K+2R\right) $. This reflects the coincidence $%
D_{3}=A_{3} $ or $\mathbf{so}\left( 6\right) =\mathbf{sl}\left( 4\right) $.
It is because $S^{+}$ is the standard representation of $\mathbf{sl}\left(
4\right) $ and therefore $W_{3}=\Lambda ^{2}S^{+}$. Similarly we have $%
\mathcal{S}^{-}=\Lambda ^{3}\mathcal{S}^{+}\otimes O\left( K+R\right) $
corresponding to $S^{-}=\Lambda ^{3}S^{+}$ as representation of $D_{3}$.
\end{example}

We now come back to discuss the decomposition of $\mathbf{L}_{n}$ when $n=6$%
, we have 
\begin{equation*}
\mathbf{L}_{6}|_{D_{5}}=\mathbf{1}+W_{5}+S^{+}\text{.}
\end{equation*}
On the representation bundle level, we have 
\begin{equation*}
\mathcal{L}_{6}=O_{X_{6}}\left( -K_{X}-R\right) +\mathcal{W}_{5}+\mathcal{S}%
^{+}\text{.}
\end{equation*}
From the proof of the previous proposition, there is a unique line $L$ on $%
X_{6}$ with $L\cdot R>1$. In fact we have $L\cdot R=2$ and its divisor class
is linearly equivalent to $-K_{X_{6}}-R$. The rank decomposition of $%
\mathcal{L}_{6}$ is quite familiar: $27=1+10+16$. This reminds us the fact
that, given any line $L$ on a cubic surface, there are exactly $10$ lines
intersecting $L$ and $16$ lines disjoint from $L$. This translate into $%
\mathcal{L}_{6}=O_{X_{6}}\left( L\right) +\pi ^{\ast }\mathcal{R}_{5}\otimes
O\left( -L\right) +\pi ^{\ast }\mathcal{L}_{5}$, where $\pi
:X_{6}\rightarrow X_{5}$ is the blow down of $L$.

To explain it, we remind ourselves that there is an outer-automorphism of $%
E_{6}$ which corresponds to the $\mathbb{Z}_{2}$ symmetry of the Dynkin
diagram of $E_{6}$. In particular it interchanges the two nodes
corresponding to fundamental representations $\mathcal{R}_{6}$ and $\mathcal{%
L}_{6}$. In terms of the geometry of the cubic surface $X_{6}$, it
establishes a correspondence between rulings on $X_{6}$ and lines on $X_{6}$%
. Concretely, given any divisor $R$ on $X_{6}$ which determines a ruling,
there is a unique line $L$ with $L\cdot R=2$ as we mentioned earlier.
Moreover the two decompositions of $\mathcal{L}_{6}$ are identified under
this $\mathbb{Z}_{2}$ symmetry. Namely $\mathcal{W}_{5}\cong \pi ^{\ast }%
\mathcal{R}_{5}\otimes O\left( -L\right) $ and $\mathcal{S}^{+}\cong \pi
^{\ast }\mathcal{L}_{5}$.

When $n=7$ or $8,$ the decompositions for the representations $\mathbf{L}%
_{n} $ are different. Namely 
\begin{equation*}
\mathbf{L}_{7}|_{D_{6}\times A_{1}}=W_{6}\otimes \Lambda _{1}+S^{+}
\end{equation*}
and 
\begin{equation*}
\mathbf{L}_{8}|_{D_{8}}=W_{8}+S^{+}.
\end{equation*}
The corresponding decompositions of representation bundles are 
\begin{equation*}
\mathcal{L}_{7}=\mathcal{W}_{6}\otimes \Lambda _{1}+\mathcal{S}^{+}\text{,}
\end{equation*}
and 
\begin{equation*}
\mathcal{L}_{8}=\mathcal{W}_{8}+\mathcal{S}^{+}\text{.}
\end{equation*}
The details for these two cases will be discussed in later sections.

\begin{center}
\textbf{Decompose }$\mathcal{R}_{n}\mathbf{\ }$\textbf{under }$\mathcal{D}%
_{n-1}$
\end{center}

We have the following decomposition of $\mathcal{R}_{n}$ as a representation
bundle of $L\mathcal{D}_{n-1}:$

\begin{eqnarray*}
\mathcal{R}_{n} &=&O\left( R\right) \left( O+\mathcal{S}^{-}\right) ,\text{
\thinspace \thinspace \thinspace for }n\leq 4, \\
\mathcal{R}_{5} &=&O\left( R\right) \left( O+\mathcal{S}^{-}+O\left(
-K-2R\right) \right) , \\
\mathcal{R}_{6} &=&O\left( R\right) \left( O+\mathcal{S}^{-}+\mathcal{W}%
_{5}\otimes O\left( -K-2R\right) \right) , \\
\mathcal{R}_{7} &=&O\left( R\right) \left( S^{2}\Lambda _{1}+\mathcal{S}%
^{-}\otimes \Lambda _{1}+\Lambda ^{2}\mathcal{W}_{6}\otimes \left(
-K-2R\right) \right) \text{.}
\end{eqnarray*}
For $\mathcal{R}_{7}$ we decompose it under $L\mathcal{D}_{6}+L\mathcal{A}%
_{1}$. If we restrict to any fiber, it recovers the Lie algebra facts: 
\begin{eqnarray*}
\mathbf{R}_{n} &=&\mathbf{1}+S^{-},\text{ \thinspace \thinspace \thinspace
for }n\leq 4 \\
\mathbf{R}_{5} &=&\mathbf{1}+S^{-}+\mathbf{1}, \\
\mathbf{R}_{6} &=&\mathbf{1}+S^{-}+W_{5}, \\
\mathbf{R}_{7} &=&\mathbf{1\otimes }S^{2}\Lambda _{1}+S^{-}\otimes \Lambda
_{1}+\Lambda ^{2}W_{6}\otimes \mathbf{1}\text{.}
\end{eqnarray*}

We see from the above decompositions that, after we fix a ruling $R$ on $%
X_{n}$, then other rulings on $X_{n}$ can be obtained by adding an element
of $\mathcal{S}^{-}$. In fact, when $n$ is not bigger than four, it gives
all rulings on $X_{n}$. First if $O\left( T\right) $ is a direct summand of $%
\mathcal{S}^{-}$, namely $T$ satisfies $T^{2}=-2,T\cdot K=0$ and $T\cdot R=1$%
, then it is easy to check that $R^{\prime }=R+T$ satisfies $R^{\prime
2}=0,R^{\prime }\cdot K=-2.$ Hence $R^{\prime }$ is another ruling on $X_{n}$%
. Moreover $R\cdot R^{\prime }=1$, that is $R^{\prime }$ is a section of the
ruling determined by $R$. To prove the reverse direction, we need the
following geometric proposition.

\begin{proposition}
If we fix a ruling $R$ on $X_{n}$, then any other ruling $R^{\prime }$ is a
section of it, i.e. $R\cdot R^{\prime }=1$, provided that $n\leq 4$.

When $n=5$, then either $R^{\prime }$ is a section or $R^{\prime }=-K-R$. In
the latter case, $R^{\prime }$ is a 2-section, i.e. $R\cdot R^{\prime }=2$.
\end{proposition}

Proof of proposition: We assume that $R$ and $R^{\prime }$ are two rulings
on $X_{n}$. Without loss of generality we assume that both $R$ and $%
R^{\prime }$ are singular fibers of the corresponding ruling. In particular
each of them consists of two lines. Recall that,when $n\leq 6$, every two
lines intersect at not more than one point. Therefore if $R\cdot R^{\prime
}>1$, either $R\cup R^{\prime }$ is a rectangle or it contains a triangle
subconfiguration. There is no triangle when $n\leq 5$ and there is no
rectangle when $n\leq 4$ and the only type of rectangle on $X_{5}$ is an
anti-canonical divisor. This proves the proposition. $\square $

Now if $R^{\prime }$ is another ruling with $R\cdot R^{\prime }=1$ then, by
reversing previous arguments, $O\left( T\right) $ is a direct summand of $%
\mathcal{S}^{-}$. Therefore we have verified above decompositions of $%
\mathcal{R}_{n}$ for $n\leq 5$. For $n=6$, we have isomorphism $\mathcal{R}%
_{6}=\mathcal{L}_{6}^{\ast }\otimes O\left( -K\right) $ and hence the
decomposition for $\mathcal{L}_{6}$ under $L\mathcal{D}_{5}$ gives the one
for $\mathcal{R}_{6}$. For $n=7$, there is an isomorphism $\mathcal{R}_{7}=L%
\mathcal{E}_{7}\otimes O\left( -K\right) $. In section seven, we will
discuss a decomposition of $L\mathcal{E}_{7}$ under $L\mathcal{D}_{6}+L%
\mathcal{A}_{1}$ and hence we obtain a decomposition for $\mathcal{R}_{7}$
which is the one given above.

\begin{center}
\textbf{Decompose }$L\mathcal{E}_{n}\mathbf{\ }$\textbf{under }$L\mathcal{D}%
_{n-1}$
\end{center}

We have the following decomposition of $L\mathcal{E}_{n}$ as a
representation bundle of $L\mathcal{D}_{n-1}$ (or $L\left( \mathcal{D}_{7}+%
\mathcal{A}_{1}\right) $, or $L\mathcal{D}_{8}$):

\begin{eqnarray*}
L\mathcal{E}_{n} &=&L\mathcal{D}_{n-1}+O+\mathcal{S}^{-}+\left( \mathcal{S}%
^{-}\right) ^{\ast },\text{ \thinspace \thinspace \thinspace for }n\leq 6, \\
L\mathcal{E}_{7} &=&L\mathcal{D}_{7}+L\mathcal{A}_{1}+\mathcal{S}^{-}\otimes
\Lambda _{1}^{\ast }, \\
L\mathcal{E}_{8} &=&L\mathcal{D}_{8}+\mathcal{S}^{+}\text{.}
\end{eqnarray*}
Notice that when $n=2m$ is even, we have an isomorphism $\left( \mathcal{S}%
^{-}\right) ^{\ast }=\mathcal{S}^{+}\otimes O\left( \left( 4-m\right)
R+K\right) $ and when $n=2m+1$ is odd, we have an isomorphism $\left( 
\mathcal{S}^{-}\right) ^{\ast }=\mathcal{S}^{-}\otimes O\left( \left(
4-m\right) R+K\right) $. Moreover in the case of $X_{7}$, we have an
isomorphism $\Lambda _{1}^{\ast }=\Lambda _{1}\otimes O\left( R+K\right) $
(see section seven). So we can express the decompositions of $L\mathcal{E}%
_{n}$ in terms of representation bundles of $L\mathcal{D}_{n-1}$ using only
fundamental representations of $D_{n-1}$. More precisely, we have 
\begin{eqnarray*}
L\mathcal{E}_{n} &=&L\mathcal{D}_{n-1}+O+\mathcal{S}^{-}+\mathcal{S}%
^{-}\otimes O\left( \left( 4-m\right) R+K\right) ,\text{ \thinspace
\thinspace \thinspace for }n=2m+1\leq 5, \\
L\mathcal{E}_{n} &=&L\mathcal{D}_{n-1}+O+\mathcal{S}^{-}+\mathcal{S}%
^{+}\otimes O\left( \left( 4-m\right) R+K\right) ,\text{ \thinspace
\thinspace \thinspace for }n=2m\leq 6, \\
L\mathcal{E}_{7} &=&L\mathcal{D}_{7}+L\mathcal{A}_{1}+\mathcal{S}^{-}\otimes
\Lambda _{1}\otimes O\left( R+K\right) , \\
L\mathcal{E}_{8} &=&L\mathcal{D}_{8}+\mathcal{S}^{+}\text{.}
\end{eqnarray*}
If we restrict to any fiber, it recovers the Lie algebra facts: 
\begin{eqnarray*}
E_{n} &=&D_{n-1}+\mathbf{1}+S^{-}+S^{-},\text{ \thinspace \thinspace
\thinspace for }n=2m+1\leq 5, \\
E_{n} &=&D_{n-1}+\mathbf{1}+S^{-}+S^{+},\text{ \thinspace \thinspace
\thinspace for }n=2m\leq 6, \\
E_{7} &=&D_{7}+LA_{1}+S^{-}\otimes \Lambda _{1}, \\
E_{8} &=&D_{8}+S^{+}\text{.}
\end{eqnarray*}

Now to show that $L\mathcal{E}_{n}=L\mathcal{D}_{n-1}+O+\mathcal{S}%
^{-}+\left( \mathcal{S}^{-}\right) ^{\ast },$ \thinspace \thinspace
\thinspace for $n\leq 6$, we need the following geometric proposition.

\begin{proposition}
If we fix a ruling $R$ on $X_{n}$, then any $D$ with $D^{2}=-2$ and $D\cdot
K=0$ satisfies $\left| D\cdot R\right| \leq 1$ provided that $n\leq 6$.
\end{proposition}

Proof of proposition: We can assume that $D=l_{1}-l_{2}$ with $l_{1}\cdot
l_{2}=0$. By considering a singular fiber of the ruling, we can assume that $%
R=l_{3}+l_{4}$ with $l_{3}\cdot l_{4}=1$. Here $l_{i}$'s are all lines on $%
X_{n}$. We suppose that $D\cdot R=\left( l_{1}-l_{2}\right) \cdot \left(
l_{3}+l_{4}\right) \geq 2$. We first assume that all $l_{i}$'s are distinct.
Then $\left( l_{1}-l_{2}\right) \cdot \left( l_{3}+l_{4}\right) \leq
l_{1}\cdot \left( l_{3}+l_{4}\right) \leq 2$ because two lines on $X_{n}$
intersect at not more than one point when $n\leq 6$. Therefore all
inequality sign are equality and we get $l_{2}\cdot l_{3}=l_{2}\cdot l_{4}=0$
and $l_{1}\cdot l_{3}=l_{1}\cdot l_{4}=1$. Together with $l_{3}\cdot l_{4}=1$%
, it implies that $l_{1}\cup l_{3}\cup l_{4}$ forms a triangle. So $n$ must
be six and the triangle is anti-canonical. Hence $1=-K\cdot l_{2}=\left(
l_{1}+l_{3}+l_{4}\right) \cdot l_{2}=0$ gives the contradiction.

If some of the $l_{i}$'s are the same, the only possible way that $D\cdot R$
can be bigger than one is $l_{2}=l_{3}$ (or $l_{4}$). In this case we have 
\begin{eqnarray*}
D\cdot R &=&\left( l_{1}-l_{2}\right) \cdot \left( l_{2}+l_{4}\right) \\
&=&l_{1}\cdot l_{4}-l_{2}^{2}-l_{2}\cdot l_{4} \\
&=&l_{1}\cdot l_{4}+1-1 \\
&=&l_{1}\cdot l_{4}\leq 1.
\end{eqnarray*}
Hence we always have $D\cdot R\leq 1.$ By replacing $D$ by $-D$ if
necessary, we obtain $\left| D\cdot R\right| \leq 1$ provided that $n\leq 6$%
. Hence the proposition. $\square $

The above equality as vector bundle decomposition of $L\mathcal{E}_{n}$
follows immediately from the proposition and the definition of $L\mathcal{D}%
_{n-1}$ and $\mathcal{S}^{-}$ when $n\leq 6$. In fact, such decomposition is
a decomposition of $L\mathcal{E}_{n}$ as $L\mathcal{D}_{n-1}$-representation
bundles. The decompositions for $L\mathcal{E}_{7}$ and $L\mathcal{E}_{8}$
will be discussed in later sections.

\begin{center}
\textbf{Reconstructing }$L\mathcal{E}_{n}\mathbf{\ }$\textbf{from }$L%
\mathcal{D}_{n-1}$
\end{center}

We can also use the above vector bundle decompositions to reconstruct the
fiberwise Lie algebra structure on $L\mathcal{E}_{n}$ from $L\mathcal{D}%
_{n-1}$. This provides us with an alternative way to describe the Lie
algebra bundle $L\mathcal{E}_{n}$. When $n\leq 6$ the construction is
completely analogous to the construction of $L\mathcal{E}_{n}$ from its
decomposition under $L\mathcal{E}_{n-1}.$ We omit the details here. The $n=8$
case is different and it will be discussed in section eight.\newpage

\section{$\mathcal{A}_{n-2}$-bundle over ruled surfaces with a section}

\subsection{Reduction to $A_{n-2}$-bundle over $X_{n}$ using $\mathcal{S}%
^{+} $}

We continue our notations from the $D_{n-1}$-bundle $\mathcal{D}_{n-1}$ over 
$X_{n}$. The Dynkin diagram of $A_{n-2}=\mathbf{sl}\left( n-1\right) $ can
be obtained by removing the node in the Dynkin diagram of $D_{n-1}$ which
corresponds to the spinor representation $S^{+}$ (or $S^{-}$ which we will
discuss in the next subsection).

To reduce $\mathcal{D}_{n-1}$ to an $A_{n-2}$-bundle over $X_{n}$, we choose
a line section $S$ to the ruling $R$ over $X_{n}$. This corresponds to a
direct summand of $\mathcal{S}^{+}=\bigoplus O_{X}\left( S\right) $. Each
singular fiber of the ruling consists of two lines. One of them, say $C$,
intersect the section $S$ and it satisfies $C^{2}=-1=C\cdot K,C\cdot R=0$
and $C\cdot S=1$. Using these, we define the vector bundle 
\begin{equation*}
\mathcal{\Lambda }_{n-2}=\bigoplus_{\substack{ C^{2}=-1  \\ CK=-1  \\ %
CR=0\,\,  \\ CS=1\,\,}}O_{X_{n}}\left( C\right)
\end{equation*}
over $X_{n}$. Notice that the other line of the singular fiber containing $C$
is linear equivalent to $-C+R$. Therefore we have 
\begin{equation*}
\mathcal{W}_{n-1}=\mathcal{\Lambda }_{n-2}+\mathcal{\Lambda }_{n-2}^{\ast
}\otimes O_{X_{n}}\left( R\right) .
\end{equation*}
In particular, the rank of $\mathcal{\Lambda }_{n-2}$ equals $n-1$. It is
not difficult to check that there is an isomorphism 
\begin{equation*}
\det :\det \mathcal{\Lambda }_{n-2}\overset{\cong }{\rightarrow }%
O_{X_{n}}\left( -K_{X}-2S+\left( n-4\right) R\right) \text{.}
\end{equation*}
Now the automorphism bundle of $\mathcal{\Lambda }_{n-2}$ preserving $\det $
is an $A_{n-2}$-bundle $\mathcal{A}_{n-2}$ over $X_{n}$. If we denote the $%
l^{th}$-wedge product of $\mathcal{\Lambda }_{n-2}$ as $\mathcal{\Lambda }%
_{n-2}^{l}$. For instance we have $\mathcal{\Lambda }_{n-2}^{0}=\mathbf{1}$
and $\mathcal{\Lambda }_{n-2}^{n-1}=\det \mathcal{\Lambda }_{n-2}$. Then $%
\mathcal{\Lambda }_{n-2}^{l}$ for $1\leq l\leq n-2$ are representation
bundles of $\mathcal{A}_{n-1}$ corresponding to all the fundamental
representations of $A_{n-2}=\mathbf{sl}\left( n-1\right) $. We also have $L%
\mathcal{A}_{n-2}=End_{0}\left( \mathcal{\Lambda }_{n-2}\right) $.

We have the following numerical characterization of $L\mathcal{A}_{n-2}$,
whose proof we similar to previous arguments and therefore skipped.

\begin{proposition}
We have the following decomposition of $L\mathcal{A}_{n-2}$ over $X_{n}:$%
\begin{equation*}
L\mathcal{A}_{n-2}=O^{\oplus n-2}+\bigoplus_{\substack{ D^{2}=-2  \\ DK=0 
\\ DR=0  \\ DS=0}}O\left( D\right)
\end{equation*}
\end{proposition}

\bigskip

\begin{center}
\textbf{Decompose }$\mathcal{D}_{n-1}$-\textbf{bundles under }$\mathcal{A}%
_{n-2}$
\end{center}

From the previous isomorphism, we have $\mathcal{\Lambda }_{n-2}^{\ast }=%
\mathcal{\Lambda }_{n-2}^{n-2}\otimes \det \mathcal{\Lambda }_{n-2}^{\ast }=%
\mathcal{\Lambda }_{n-2}^{n-2}\otimes O\left( K+2S+\left( 4-n\right)
R\right) .$ Using this, we get 
\begin{equation*}
\mathcal{W}_{n-1}=\mathcal{\Lambda }_{n-2}+\mathcal{\Lambda }%
_{n-2}^{n-2}\otimes O\left( K+2S+\left( 5-n\right) R\right) .
\end{equation*}
Therefore we get similar decompositions for wedge products of $\mathcal{W}%
_{n-1}$ in terms of fundamental representation bundles of $\mathcal{A}_{n-2}$%
. For example, using the isomorphisms $L\mathcal{D}_{n-1}=\Lambda ^{2}%
\mathcal{W}_{n-1}\otimes O\left( -R\right) $ and $L\mathcal{A}%
_{n-2}=End_{0}\left( \mathcal{\Lambda }_{n-2}\right) $, we have the
following isomorphism between representation bundles of $L\mathcal{A}_{n-2}$%
: 
\begin{equation*}
L\mathcal{D}_{n-1}=L\mathcal{A}_{n-2}+O+\mathcal{\Lambda }_{n-2}^{2}\otimes
O\left( -R\right) +\left( \mathcal{\Lambda }_{n-2}^{2}\right) ^{\ast
}\otimes O\left( R\right) .
\end{equation*}

Besides wedge products of $\mathcal{W}_{n-1}$, there are two other
fundamental representation bundles of $\mathcal{D}_{n-1}$, namely $\mathcal{S%
}^{+}$ and $\mathcal{S}^{-}$. For them we have the following decompositions: 
\begin{equation*}
\mathcal{S}^{+}=\sum_{l=0}^{\left[ \frac{n-1}{2}\right] }\mathcal{\Lambda }%
_{n-2}^{2l}\otimes O\left( S-lR\right) ,
\end{equation*}
and 
\begin{equation*}
\mathcal{S}^{-}=\sum_{l=1}^{\left[ \frac{n}{2}\right] }\mathcal{\Lambda }%
_{n-2}^{2l-1}\otimes O\left( S-lR\right) .
\end{equation*}
We leave the verification of these two formulae for our readers. These
decompositions of representation bundles correspond to decomposition of
representations of $D_{n-1}$ under $A_{n-2}$: 
\begin{eqnarray*}
D_{n-1} &=&A_{n-2}+\mathbf{1}+\Lambda _{n-2}^{2}+\left( \Lambda
_{n-2}^{2}\right) ^{\ast }, \\
W_{n-1} &=&\Lambda _{n-2}+\Lambda _{n-2}^{\ast }=\Lambda _{n-2}+\Lambda
_{n-2}^{n-2}, \\
S^{+} &=&\sum_{l=0}^{\left[ \frac{n-1}{2}\right] }\Lambda _{n-2}^{2l}\text{
and }S^{-}=\sum_{l=1}^{\left[ \frac{n}{2}\right] }\Lambda _{n-2}^{2l-1}.
\end{eqnarray*}

\subsection{Reduction to $A_{n-2}$-bundle over $X_{n}$ using $\mathcal{S}%
^{-} $}

As we mentioned that the Dynkin diagram of $A_{n-2}=\mathbf{sl}\left(
n-1\right) $ can also be obtained by removing the node in the Dynkin diagram
of $D_{n-1}$ which corresponds to the spinor representation $S^{-}$, we can
also reduce $\mathcal{D}_{n-1}$ to an $A_{n-2}$-bundle over $X_{n}$ by
choosing a direct summand $O\left( T\right) $ of $\mathcal{S}^{-}$. This
means that $T$ satisfies $T^{2}=-2$, $T\cdot K=0$ and $T\cdot R=1$ or
equivalently $T+R$ is a ruling section of the ruling determined by $R$.

Similar to the previous subsection, we consider the vector bundle 
\begin{equation*}
\mathcal{\bar{\Lambda}}_{n-2}=\bigoplus_{\substack{ C^{2}=-1  \\ CK=-1  \\ %
CR=0\,\,  \\ CT=1\,\,}}O_{X_{n}}\left( C\right)
\end{equation*}
over $X_{n}$. Also we have 
\begin{equation*}
\mathcal{W}_{n-1}=\mathcal{\bar{\Lambda}}_{n-2}+\mathcal{\bar{\Lambda}}%
_{n-2}^{\ast }\otimes O_{X_{n}}\left( R\right) .
\end{equation*}
In particular, the rank of $\mathcal{\bar{\Lambda}}_{n-2}$ equals $n-1$. It
is not difficult to check that there is an isomorphism 
\begin{equation*}
\det :\det \mathcal{\bar{\Lambda}}_{n-2}\overset{\cong }{\rightarrow }%
O_{X_{n}}\left( -K_{X}-2T+\left( n-5\right) R\right) \text{.}
\end{equation*}
Now the automorphism bundle of $\mathcal{\bar{\Lambda}}_{n-2}$ preserving $%
\det $ is a principal $A_{n-2}$-bundle $\overline{\mathcal{A}}_{n-2}$ over $%
X_{n}$. If we denote the $l^{th}$-wedge product of $\mathcal{\bar{\Lambda}}%
_{n-2}$ as $\mathcal{\bar{\Lambda}}_{n-2}^{l}$. For instance we have $%
\mathcal{\bar{\Lambda}}_{n-2}^{0}=\mathbf{1}$ and $\mathcal{\bar{\Lambda}}%
_{n-2}^{n-1}=\det \mathcal{\bar{\Lambda}}_{n-2}$. Then $\mathcal{\bar{\Lambda%
}}_{n-2}^{l}$ for $1\leq l\leq n-2$ are representation bundles of $\overline{%
\mathcal{A}}_{n-2}$ corresponding to all the fundamental representations of $%
A_{n-2}=\mathbf{sl}\left( n-1\right) $. From similar reasons as before, we
have the following proposition.

\begin{proposition}
We have the following decomposition of $L\overline{\mathcal{A}}_{n-2}$ over $%
X_{n}:$%
\begin{equation*}
L\overline{\mathcal{A}}_{n-2}=O^{\oplus n-2}+\bigoplus_{_{\substack{ %
D^{2}=-2  \\ DK=0  \\ DR=0  \\ DT=0}}}O\left( D\right)
\end{equation*}
\end{proposition}

From the previous isomorphism, we have $\mathcal{\bar{\Lambda}}_{n-2}^{\ast
}=\mathcal{\bar{\Lambda}}_{n-2}^{n-2}\otimes \det \mathcal{\bar{\Lambda}}%
_{n-2}^{\ast }=\mathcal{\bar{\Lambda}}_{n-2}^{n-2}\otimes O\left(
K+2T+\left( 5-n\right) R\right) .$ Using this, we get 
\begin{equation*}
\mathcal{W}_{n-1}=\mathcal{\bar{\Lambda}}_{n-2}+\mathcal{\bar{\Lambda}}%
_{n-2}^{n-2}\otimes O\left( K+2T+\left( 6-n\right) R\right) .
\end{equation*}
Therefore we get similar decompositions for wedge products of $\mathcal{W}%
_{n-1}$ in terms of fundamental representation bundles of $\overline{%
\mathcal{A}}_{n-2}$. For example, we have the following isomorphism between
representation bundles of $L\overline{\mathcal{A}}_{n-2}$: 
\begin{equation*}
L\mathcal{D}_{n-1}=L\overline{\mathcal{A}}_{n-2}+O+\mathcal{\bar{\Lambda}}%
_{n-2}^{2}\otimes O\left( -R\right) +\left( \mathcal{\bar{\Lambda}}%
_{n-2}^{2}\right) ^{\ast }\otimes O\left( R\right) .
\end{equation*}
Besides wedge products of $\mathcal{W}_{n-1}$, there are two other
fundamental representation bundles of $\mathcal{D}_{n-1}$, namely $\mathcal{S%
}^{+}$ and $\mathcal{S}^{-}$. For them we have the following decompositions: 
\begin{equation*}
\mathcal{S}^{+}=\sum_{l=1}^{\left[ \frac{n}{2}\right] }\mathcal{\bar{\Lambda}%
}_{n-2}^{2l-1}\otimes O\left( T-\left( l-1\right) R\right) .
\end{equation*}
and 
\begin{equation*}
\mathcal{S}^{-}=\sum_{l=0}^{\left[ \frac{n-1}{2}\right] }\mathcal{\bar{%
\Lambda}}_{n-2}^{2l}\otimes O\left( T-lR\right) ,
\end{equation*}

We leave the verification of these two formulae for our readers. These
decompositions of representation bundles correspond to decomposition of
representations of $D_{n-1}$ under $A_{n-2}$: 
\begin{eqnarray*}
D_{n-1} &=&A_{n-2}+\mathbf{1}+\Lambda _{n-2}^{2}+\left( \Lambda
_{n-2}^{2}\right) ^{\ast }, \\
W_{n-1} &=&\Lambda _{n-2}+\Lambda _{n-2}^{\ast }=\Lambda _{n-2}+\Lambda
_{n-2}^{n-2}, \\
S^{+} &=&\sum_{l=0}^{\left[ \frac{n-1}{2}\right] }\Lambda _{n-2}^{2l}\text{
and }S^{-}=\sum_{l=1}^{\left[ \frac{n}{2}\right] }\Lambda _{n-2}^{2l-1}.
\end{eqnarray*}

\newpage

\section{Quartic surface $X_{5}$ in $\mathbb{P}^{4}$ and its $E_{5}$-bundle}

In this section we study the surface $X_{5}$ obtained by blowing up $\mathbb{%
P}^{2}$ at five generic points. It is a complete intersection of two quadric
hypersurfaces in $\mathbb{P}^{4}$. In earlier section, we have constructed a 
$E_{5}$-bundle $L\mathcal{E}_{5}$ over $X_{5}$ and its associated bundles $%
\mathcal{L}_{5}$ and $\mathcal{R}_{5}$ using lines and rulings on $X_{5}$.

Recall that $E_{5}$ is a classical Lie algebra, namely $E_{5}=D_{5}=\mathbf{%
so}\left( 10\right) $. Moreover the $E_{5}$ fundamental representation $%
\mathbf{R}_{5}$ is just the standard representation of $\mathbf{so}\left(
10\right) $. Correspondingly the bundle $\mathcal{R}_{5}$ on $X_{5}$ has a
fiberwise non-degenerate quadratic form 
\begin{equation*}
q_{5}:\mathcal{R}_{5}\otimes \mathcal{R}_{5}\rightarrow O\left( -K\right) 
\text{.}
\end{equation*}
We now describe this quadratic form on $\mathcal{R}_{5}$. If we choose two
rulings $R_{1}$ and $R_{2}$ on $X_{5}$ and write each of them as a sum of
two intersecting lines, i.e. $R_{1}=l_{1}+l_{2}$ and $R_{2}=l_{3}+l_{4}$
with $l_{1}\cdot l_{2}=l_{3}\cdot l_{4}=1$. In general we have $R_{1}\cdot
R_{2}\leq 2$. When the equality sign holds, then these four lines must form
a rectangle on $X_{5}$ because there is no triangle on $X_{n}$ when $n<6$.
Moreover every rectangle on $X_{5}$ is an anti-canonical divisor. This
determines $q_{5}$ on $O\left( R_{1}\right) \otimes O\left( R_{2}\right) $.
Then the quadratic form $q_{5}$ on the whole $\mathcal{R}_{5}$ is obtained
by extending these linearly.

The automorphism bundle of $\mathcal{R}_{5}$ preserving $q_{5}$ is a $E_{5}$
bundle over $X_{5}$ which was denoted $\mathcal{E}_{5}$ in the introduction.
Moreover the associated Lie algebra of $\mathcal{E}_{5}$ is just $L\mathcal{E%
}_{5}$. If we regard $L\mathcal{E}_{5}$ as a Lie algebra bundle of type $%
\mathbf{so}\left( 10\right) $, then $\mathcal{L}_{5}$ is its representation
bundle corresponding to a spinor representation of $\mathbf{so}\left(
10\right) $, which has rank equals 16$.$

If we fix a ruling $R$ on $X_{5}$, then we break the symmetry from $E_{5}$
to $D_{4}=\mathbf{so}\left( 8\right) $. Their corresponding representation
bundles $\mathcal{W}_{4}$, $\mathcal{S}^{+}$ and $\mathcal{S}^{-}$ have been
discussed in section three. Nevertheless, there is an additional
homomorphism among them. 
\begin{equation*}
f:\mathcal{S}^{+}\otimes \mathcal{S}^{-}\otimes \mathcal{W}_{4}\rightarrow
O_{X_{5}}\left( -K-R\right) .
\end{equation*}
This homomorphism $f$ is constructed via isomorphisms $O_{X}\left(
-K_{X}-R\right) $ and tensor products of various direct summands of $%
\mathcal{S}^{+},\mathcal{S}^{-}$ and $\mathcal{W}_{4}$. For instance, using
notations in the example of $D_{4}$ in last section, $O\left(
H-L_{1}-L_{2}\right) ,O\left( L_{0}-L_{3}\right) $ and $O\left(
H-L_{0}-L_{4}\right) $ are direct summands of $\mathcal{S}^{+},\mathcal{S}%
^{-}$ and $\mathcal{W}_{4}$ respectively. We have isomorphism 
\begin{equation*}
O\left( H-L_{1}-L_{2}\right) \otimes O\left( L_{0}-L_{3}\right) \otimes
O\left( H-L_{0}-L_{4}\right) \cong O\left( -K-R\right)
\end{equation*}
where $R=H-L_{0}.$ We expect that this homomorphism $f$ would be used to
describe triality among these bundles in a manner similar to the triality
for the Lie algebra $\mathbf{so}\left( 8\right) $.

\subsection{Reduction to $\mathbf{sl}\left( 4\right) \times \mathbf{sl}%
\left( 2\right) \times \mathbf{sl}\left( 2\right) $-bundle}

In this subsection, we discuss a degeneration of $X_{5}$ into nonnormal del
Pezzo surface which is a normal crossing variety. In terms of $L\mathcal{E}%
_{5}$, it corresponds to reducing the structure group from $E_{5}=\mathbf{so}%
\left( 10\right) $ to $A_{3}\times A_{1}\times A_{1}=\mathbf{sl}\left(
4\right) \times \mathbf{sl}\left( 2\right) \times \mathbf{sl}\left( 2\right)
=\mathbf{so}\left( 6\right) \times \mathbf{so}\left( 4\right) $. Recall that 
$X_{5}$ is a complete intersection of two quadric hypersurfaces $Q$ and $%
Q^{\prime }$ in $\mathbb{P}^{4}$. Now we degenerate one of these quadric
hypersurfaces, say $Q^{\prime }$, into a union of two hyperplanes $H_{1}\cup
H_{2}$ and denote the degenerating family by $Q^{\prime }\left( t\right) $
with $Q^{\prime }\left( 0\right) =H_{1}\cup H_{2}$. We define the
corresponding family of del Pezzo surfaces as $X_{5}\left( t\right) =Q\cap
Q^{\prime }\left( t\right) $\ \textbf{. }So $X_{5}\left( 0\right) =Q_{\left(
1\right) }\cup Q_{\left( 2\right) }$ where $Q_{\left( i\right) }=Q\cap H_{i}$%
. Each $Q_{\left( i\right) }$ is a quadric surface in $\mathbb{P}^{3}\simeq
H_{i}$, so it has two rulings by lines in $H_{i}\subset \mathbb{P}^{4}.$

We want to understand what happens to these lines on $X_{5}\left( t\right) $
as $t$ approaches zero. The first step is to see which members of the two
rulings on $Q_{\left( i\right) }$ can be deformed to lines in $X_{5}\left(
t\right) $ for small nonzero $t$. Let $C$ be the singular locus of $%
X_{5}\left( 0\right) $. Then $C$ is a conic curve $C=Q\cap H_{1}\cap H_{2}$
inside $H_{1}\cap H_{2}\simeq \mathbb{P}^{2}$. As we vary $t$ away from
zero, $Q^{\prime }\left( t\right) $ intersects $C$ at four points\footnote{%
It should really be a zero dimensional scheme of length four. For simplicity
we assume that $Q$ is generic so that the four points are distinct.}
infinitesimally. We call them $q_{1},q_{2},q_{3}$ and $q_{4}$. More
precisely, if we let $Z\left( t\right) $ be the intersection of $C$ with $%
Q^{\prime }\left( t\right) $ for nonzero $t.$ Then $Z\triangleq
\lim_{t\rightarrow 0}Z\left( t\right) =\left\{
q_{1},q_{2},q_{3},q_{4}\right\} $. Notice that $Z\subset Q_{\left( 1\right)
}\cap Q_{\left( 2\right) }=C$. Now it is not difficult to see that a line in 
$Q_{\left( 1\right) }\cup Q_{\left( 2\right) }$ has a first order
deformation inside $X_{5}\left( t\right) $ if and only if it passes through
one of the $q_{i}$'s. In fact, all these infinitesimal deformations are
unobstructed and exhaust all possible lines in $X_{5}\left( t\right) $ for
small nonzero $t$.

\bigskip

For example, there are two lines in $Q_{\left( 1\right) }$ passes through
each $q_{j}$ which come from the two rulings on $Q_{\left( 1\right) }$. The
same count holds for $Q_{\left( 2\right) }$. Therefore the number of lines
on $X_{5}$ equals $4\times 2+4\times 2=16$. We want to obtain a
corresponding decomposition for $\mathcal{L}_{5}\left( t\right) $ as $t$
approach zero. On each of the $Q_{\left( i\right) }$'s, there is a rank two
vector bundle $\mathcal{R}_{\left( i\right) }$ over it given by the two
ruling on it. Namely

\begin{equation*}
\mathcal{R}_{\left( i\right) }=\bigoplus_{\substack{ R^{2}=0  \\ RK=-2}}%
O_{Q_{\left( i\right) }}\left( R\right) \text{.}
\end{equation*}
This determines a $SL\left( 2\right) $-bundle $\mathcal{A}_{\left( i\right)
} $ over $Q_{\left( i\right) }$ because $\Lambda ^{2}\mathcal{R}_{\left(
i\right) }\cong K_{Q_{\left( i\right) }}^{-1/2}$. Now giving a ruling on $%
Q_{\left( 1\right) }$ determined by $R$, a line in this ruling correspond to
the zero set of a section of $O_{Q_{\left( 1\right) }}\left( R\right) $.
Such line passes through $q_{j}$ if and only if the corresponding section
can be lift to $\frak{I}_{\left\{ q_{j}\right\} }\otimes O_{Q_{\left(
1\right) }}\left( R\right) $. Hence it is natural to define 
\begin{equation*}
\mathcal{L}_{5}\left( 0\right) =\frak{I}\otimes \mathcal{R}_{\left( 1\right)
}+\frak{I}\otimes \mathcal{R}_{\left( 2\right) },
\end{equation*}
with $\frak{I}=\bigoplus_{q\in Z}\frak{I}_{\left\{ q\right\} }$ a rank four
coherent sheaf on $X_{5}\left( 0\right) $. One should compare this
decomposition of $\mathcal{L}_{5}\left( 0\right) $ with the decomposition of
the $E_{5}$ representation $\mathbf{L}_{5}$ when we view it as a
representation of $A_{3}\times A_{1}\times A_{1}$: 
\begin{equation*}
\mathbf{L}_{5}|_{A_{3}\times A_{1}\times A_{1}}=\Lambda _{3}\otimes \Lambda
_{1}\otimes \mathbf{1}+\Lambda _{3}^{\ast }\otimes \mathbf{1}\otimes \Lambda
_{1}.
\end{equation*}
Here $\Lambda _{n}$ denotes the standard representation of $A_{n}=\mathbf{sl}%
\left( n\right) $. It would be useful to know exactly how $\mathcal{L}%
_{5}\left( t\right) $ degenerates as $t$ goes to zero and compare the limit
with $\mathcal{L}_{5}\left( 0\right) $.

Next we want to understand the behavior of $\mathcal{R}_{5}\left( t\right) $
as $t$ approach zero. Equivalently we want to study the behaviors of rulings
on $X_{5}\left( t\right) $ as $t$ approach zero. In Lie algebra term, this
will correspond to the decomposition of $E_{5}$ representation $\mathbf{R}%
_{5}$ when we view it as a representation of $A_{3}\times A_{1}\times A_{1}$%
: 
\begin{equation*}
\mathbf{R}_{5}|_{A_{3}\times A_{1}\times A_{1}}=\Lambda _{3}^{2}\otimes 
\mathbf{1}\otimes \mathbf{1}+\mathbf{1}\otimes \Lambda _{1}\otimes \Lambda
_{1}.
\end{equation*}
Here $\Lambda _{n}^{2}$ denotes the fundamental representation of $A_{n}$
which is the wedge product of $\Lambda _{n}$ with itself.

Now if we fix a ruling on $Q_{\left( 1\right) }$ and a ruling on $Q_{\left(
2\right) }$, then the sum of the corresponding lines passes through any
point on $C=Q_{\left( 1\right) }\cap Q_{\left( 2\right) }$ can be deformed
to a smooth rational curve $R$ on $X_{5}\left( t\right) $ for small nonzero $%
t$. Moreover $R$ satisfies $R^{2}=0$ and $R\cdot K_{X_{5}\left( t\right)
}=-2 $ and therefore it defines a ruling on $X_{5}\left( t\right) $.

On the other hand, if we fix two of the $q_{i}$'s, say $q_{1}$ and $q_{2}$,
then they define another type of ruling on nearby $X_{5}\left( t\right) $.
To see this, we pick one of the $Q_{\left( i\right) }$'s, say $Q_{\left(
1\right) }$. We choose the line in $Q_{\left( 1\right) }$ that passes
through $q_{1}$ in one of its ruling and choose another line in $Q_{\left(
1\right) }$ that passes through $q_{2}$ in the other ruling. Then their sum
can be deformed to a smooth rational curve on $X_{5}\left( t\right) $ for
small nonzero $t$ which determines a ruling on $X_{5}\left( t\right) $ as
before. It is not difficult to show that this ruling on $X_{5}\left(
t\right) $ is independent of the choice of $Q_{\left( i\right) }$ and the
choice of the particular ruling on $Q_{\left( i\right) }$ that gives the
line through $q_{1}$. Namely this ruling on $X_{5}\left( t\right) $ depends
only on the choice of the two points among the $q_{i}$'s.

One can also verify that this exhausts all possible ruling in a nearby $%
X_{5}\left( t\right) $. Notice that there are also four pairs of lines that
determines a particular ruling on $X_{5}$.

Therefore we define a coherent sheaf on $X_{5}\left( 0\right) $ as 
\begin{equation*}
\mathcal{R}_{5}\left( 0\right) =\frak{I}\wedge \frak{I}+\mathcal{R}_{\left(
1\right) }\otimes \mathcal{R}_{\left( 2\right) }\text{.}
\end{equation*}
It would be interesting to know if $\mathcal{R}_{5}\left( t\right) $ does
degenerate to $\mathcal{R}_{5}\left( 0\right) $ as $t$ approaches zero.

The above degeneration of $X_{5}$ into union of two copies of $\mathbb{P}%
^{1}\times \mathbb{P}^{1}$ can be described from a different viewpoint: If $%
X $ is a double cover of $\mathbb{P}^{1}\times \mathbb{P}^{1}$%
\begin{equation*}
\delta :X\rightarrow \mathbb{P}^{1}\times \mathbb{P}^{1},
\end{equation*}
branched along a smooth curve $B$ of bidegree (2,2). Then we have (See e.g 
\cite{BPV}) 
\begin{align*}
K_{X}& =\delta ^{\ast }\left( K_{\mathbb{P}^{1}\times \mathbb{P}^{1}}\otimes
O\left( 1,1\right) \right) \\
& =\delta ^{\ast }O\left( 1,1\right) .
\end{align*}
So $K^{-1}$ is ample and $K^{2}=4$. That is $X$ is a blowup of $\mathbb{P}%
^{2}$ at five points, $X_{5}$. On the other hand, requiring a degree four
del Pezzo surface in $\mathbb{P}^{4}$ to be a double cover of $\mathbb{P}%
^{1}\times \mathbb{P}^{1}$ is a codimension one condition, by dimension
counting. Each ruling on $\mathbb{P}^{1}\times \mathbb{P}^{1}$ gives $X$ a
ruling. Moreover $\delta ^{-1}\left( \mathbb{P}^{1}\times p\right) $ is a
union of two lines on $X$ if and only in $\mathbb{P}^{1}\times p$ passes
through a branched point of the double cover $B\overset{\iota }{\rightarrow }%
\mathbb{P}^{1}\times \mathbb{P}^{1}\overset{pr_{2}}{\rightarrow }\mathbb{P}%
^{1}$, where the first morphisms is the inclusion and the second morphism is
the projection to the second factor. Since $B$ is a genus one curve by the
adjunction formula, $B\overset{pr_{2}\circ \iota }{\longrightarrow }\mathbb{P%
}^{1}$ has four branched points. This gives rise to eight lines on $X$. The
other eight lines on $X$ comes from considering the other ruling on $\mathbb{%
P}^{1}\times \mathbb{P}^{1}$.

Now we degenerate $B$ to a double curve. For example $B=2\Delta $ where $%
\left\{ \Delta =\left( x,x\right) \in \mathbb{P}^{1}\times \mathbb{P}%
^{1}:x\in \mathbb{P}^{1}\right\} $ is the diagonal of $\mathbb{P}^{1}\times 
\mathbb{P}^{1}$. Then $X$ becomes a union of two copies of $\mathbb{P}%
^{1}\times \mathbb{P}^{1}$, identified along their diagonals. In the
degeneration of $B$ to $2\Delta $, an infinitesimal near $B$ will intersect $%
\Delta $ at four points $q_{1},q_{2},q_{3}$ and $q_{4}$. There are distinct
points on $\Delta $ if the degeneration is generic. In fact, these points
are the limit of the branched points of the double cover of $B$ to $\mathbb{P%
}^{1}$ for either of the two rulings. Now the appearance of the reduction to 
$\mathbf{sl}\left( 4\right) \times \mathbf{sl}\left( 2\right) \mathbf{\times
sl}\left( 2\right) $ is clear.

\begin{remark}
Recall that $X_{5}$ is the intersection of two quadric hypersurfaces $Q$ and 
$Q^{\prime }$ in $\mathbb{P}^{4}$. Instead of degenerating one of them to
hyperplanes, we can also degenerate both of them to union of hyperplanes and
study the configuration of line and ruling under such degeneration. This
degeneration of $X_{5}$ corresponds to the reduction of $E_{5}$ to $%
A_{1}\times A_{1}\times A_{1}\times A_{1}$.
\end{remark}

\newpage

\section{Cubic surface $X_{6}$ and its $E_{6}$-bundle}

When $n=6$ then $X_{6}$ is a cubic surface in $\mathbb{P}^{3}$. Studying of
lines on cubic surface has a long history in algebraic geometry. As we
mentioned earlier, every line on $X_{6}$ determines a ruling and vice versa.
Since once a line $L$ on $X_{6}$ is chosen then we have the decomposition $%
\mathcal{L}_{6}=\pi ^{\ast }\mathcal{L}_{5}+\pi ^{\ast }\mathcal{R}%
_{5}\otimes O\left( -L\right) +O\left( L\right) $ from section two. The
component $\pi ^{\ast }\mathcal{R}_{5}\otimes O\left( -L\right) $
corresponds to those lines intersecting $L$. However $E_{5}=D_{5}=\mathbf{so}%
\left( 10\right) $ and $\mathbf{R}_{5}$ corresponds to the standard
representation of $\mathbf{so}\left( 10\right) $ and therefore carries a
natural quadratic form as described in section five. There are five pairs of
intersecting lines that intersect $L$ on $X_{6}$. That is $L$ constitutes
one side of five different triangles. Remaining sides of any of these five
triangles determines a ruling on $X_{6}$ and this ruling is independent of
the choice of such triangle. As vector bundles we have isomorphism $\mathcal{%
R}_{6}\cong \mathcal{L}_{6}^{\ast }\otimes O\left( -K\right) $. In terms of
representations of $E_{6}$, we have $\mathbf{R}_{6}=\mathbf{L}_{6}^{\ast }$
which is given by the outer-automorphism of $E_{6}$ responding to the $%
\mathbb{Z}/\mathbb{Z}_{2}$ symmetry of its Dynkin diagram.

Recall that we have representation bundles homomorphism 
\begin{equation*}
c_{6}:\mathcal{L}_{6}\otimes \mathcal{L}_{6}\rightarrow \mathcal{R}_{6}
\end{equation*}
which exists for any $n$. On $X_{6}$ there is another representation bundle
homomorphism 
\begin{equation*}
c_{6}^{\ast }:\mathcal{R}_{6}\otimes \mathcal{R}_{6}\rightarrow \mathcal{L}%
_{6}\otimes O\left( -K\right) .
\end{equation*}
This is because $\mathcal{R}_{6}$ is in fact a direct summand of $\mathcal{L}%
_{6}\otimes \mathcal{L}_{6}$. The homomorphism $c_{6}^{\ast }$ is simply the
inclusion of $\mathcal{R}_{6}$ inside $\mathcal{L}_{6}\otimes \mathcal{L}%
_{6} $ using the isomorphism $\mathcal{R}_{6}\cong \mathcal{L}_{6}^{\ast
}\otimes O\left( -K\right) $.

As we mentioned before, the automorphism bundle of the pair $\left( \mathcal{%
L}_{6},\mathcal{R}_{6}\right) $ preserving $c_{6}$ and $c_{6}^{\ast }$ is a
bundle of type $E_{6}$ over $X_{6}$ which was denoted as $\mathcal{E}_{6}$.
Moreover it associated Lie algebra bundle is precisely $L\mathcal{E}_{6}$.
This follows from a direct construction of $E_{6}$ using corresponding
products on $\mathbf{L}_{6}$ and $\mathbf{R}_{6}$ (see e.g. \cite{Adams}).

If we use the isomorphism $\mathcal{R}_{6}\cong \mathcal{L}_{6}^{\ast
}\otimes O\left( -K\right) $ on the image of the homomorphism $c_{6}$. Then
we obtain a triple product 
\begin{equation*}
\mathcal{L}_{6}\otimes \mathcal{L}_{6}\otimes \mathcal{L}_{6}\rightarrow
O\left( -K\right) .
\end{equation*}
This triple product is symmetry in all three variables. In terms of the
geometry of the cubic surface, this product can be described by determining
whether three given lines on $X_{6}$ forms a triangle or not. Similarly we
also have 
\begin{equation*}
\mathcal{R}_{6}\otimes \mathcal{R}_{6}\otimes \mathcal{R}_{6}\rightarrow
O\left( -2K\right) \text{.}
\end{equation*}

\subsection{Reduction to $\mathbf{sl}\left( 3\right) \times \mathbf{sl}%
\left( 3\right) \times \mathbf{sl}\left( 3\right) $ bundle}

In this subsection we degenerate the cubic surface into union of three
generic planes in $\mathbb{P}^{3}$. We are going to see that the symmetry
group of $\mathcal{L}_{6},\mathcal{R}_{6}$ and $L\mathcal{E}_{6}$ would be
reduced from $E_{6}$ to $\mathbf{sl}\left( 3\right) \times \mathbf{sl}\left(
3\right) \times \mathbf{sl}\left( 3\right) =A_{2}\times A_{2}\times A_{2}$.

Let $X\left( t\right) $ be the family of cubic surfaces in $\mathbb{P}^{3}$
parametrized by $t$. We assume that $X\left( 0\right) $ is a union of three
planes $X\left( 0\right) =H_{1}\cup H_{2}\cup H_{3}$ and $X\left( t\right) $
is smooth for nonzero $t$.Genericity of $H_{i}$ means that their common
intersection consists of a single point. For example if $f\left(
z_{0},z_{1},z_{2},z_{3}\right) $ is a generic homogenerous cubic polynomial,
then 
\begin{equation*}
X\left( t\right) =\left\{ tf\left( z_{0},z_{1},z_{2},z_{3}\right)
+z_{1}z_{2}z_{3}=0\right\} \subset \mathbb{P}^{3},
\end{equation*}
would satisfy our assumption. We denote $C_{i}=H_{j}\cap H_{k}$ with $%
\left\{ i,j,k\right\} =\left\{ 1,2,3\right\} $. As we various $t$, $X\left(
t\right) $ intersects $C_{i}$ at three points infinitesimally. We call them $%
Z_{i}=\left\{ p_{i1},p_{i2},p_{i3}\right\} $. Then a line in $H_{3}$ which
can be deformed to one in $X\left( t\right) $ for small $t$ if and only if
it intersects $C_{1}$ and $C_{2}$ at $p_{1\alpha }$ and $p_{2\beta }$
respectively, for some $\alpha ,\beta \in \left\{ 1,2,3\right\} $. We denote
it as $\overline{p_{1\alpha }p_{2\beta }}$. The number of such lines in $%
H_{3}$ equals $9=3\times 3$. If we replace $H_{3}$ by either $H_{1}$ or $%
H_{2}$ then analogous statements would hold true. Conversely any one
parameter family of lines on $X\left( t\right) $ parametrized by $t$ would
converge to one of these lines (see \cite{Segre}). The total number of these
lines equals $27=3\times 3+3\times 3+3\times 3$.

In fact the decomposition of the representation $\mathbf{L}_{6}$ of $E_{6}$
under $A_{2}\times A_{2}\times A_{2}$ has the similar structure. If we
denote the standard representation of $A_{2}$ by $\Lambda _{2}$ and its dual
representation by $\Lambda _{2}^{\ast }$. Then 
\begin{equation*}
\mathbf{L}_{6}|_{A_{2}\times A_{2}\times A_{2}}=\Lambda _{2}\otimes \Lambda
_{2}^{\ast }\otimes \mathbf{1}+\mathbf{1}\otimes \Lambda _{2}\otimes \Lambda
_{2}^{\ast }+\Lambda _{2}^{\ast }\otimes \mathbf{1}\otimes \Lambda _{2}.
\end{equation*}
Before we discuss $\mathcal{L}_{6}\left( 0\right) $ we need to define the
coherent sheaf $\mathcal{I}_{\left( i\right) }=\oplus _{p\in Z_{i}}\mathcal{I%
}_{\left\{ p\right\} }$ where $\mathcal{I}_{\left\{ p\right\} }$ is the
ideal sheaf of $\left\{ p\right\} $ in $X_{6}\left( 0\right) $. We will use
the homomorphism $\mathcal{I}_{\left( i\right) }\otimes \mathcal{I}_{\left(
i\right) }\rightarrow O_{X_{6}\left( 0\right) }$ defined by $\left( \oplus
_{p\in Z_{i}}s_{p}\right) \otimes \left( \oplus _{p\in Z_{i}}s_{p}^{\prime
}\right) \mapsto \sum_{p\in Z_{i}}s_{p}s_{p}^{\prime }.$

Also $O_{H_{i}}\left( 1\right) $ denotes the hyperplane bundle on $H_{i}$
and we treat it as a rank one coherent sheaf on $X_{6}\left( 0\right) $
whose restriction on $X_{6}\left( 0\right) \backslash H_{i}$ is trivial. Now
sections of $\mathcal{I}_{\left( 1\right) }\otimes \mathcal{I}_{\left(
2\right) }\otimes O_{H_{3}}\left( 1\right) $ would correspond to line on $%
H_{3}$ meeting $C_{1}$ and $C_{2}$ at $p_{1\alpha }$ and $p_{2\beta }$
respectively, for some $\alpha ,\beta \in \left\{ 1,2,3\right\} $. Therefore
it is reasonable to propose that 
\begin{equation*}
\mathcal{L}_{6}\left( 0\right) =\mathcal{I}_{\left( 1\right) }\otimes 
\mathcal{I}_{\left( 2\right) }\otimes O_{H_{3}}\left( 1\right) +\mathcal{I}%
_{\left( 2\right) }\otimes \mathcal{I}_{\left( 3\right) }\otimes
O_{H_{1}}\left( 1\right) +\mathcal{I}_{\left( 3\right) }\otimes \mathcal{I}%
_{\left( 1\right) }\otimes O_{H_{2}}\left( 1\right) .
\end{equation*}

Next we discuss the triple product on $\mathcal{L}_{6}\left( 0\right) $, 
\begin{equation*}
\mathcal{L}_{6}\left( 0\right) \otimes \mathcal{L}_{6}\left( 0\right)
\otimes \mathcal{L}_{6}\left( 0\right) \rightarrow O_{X_{6}\left( 0\right)
}\left( 1\right) .
\end{equation*}
Here $O_{X_{6}\left( 0\right) }\left( 1\right) $ can be interpreted as the
anti-canonical sheaf on $X_{6}\left( 0\right) $ even though $X_{6}\left(
0\right) $ is nonnormal. Using the homomorphism $\mathcal{I}_{\left(
2\right) }\otimes \mathcal{I}_{\left( 2\right) }\rightarrow O_{X_{6}\left(
0\right) }$ defined above, we obtain a homomorphism 
\begin{equation*}
\left( \mathcal{I}_{\left( 1\right) }\otimes \mathcal{I}_{\left( 2\right)
}\otimes O_{H_{3}}\left( 1\right) \right) \otimes \left( \mathcal{I}_{\left(
2\right) }\otimes \mathcal{I}_{\left( 3\right) }\otimes O_{H_{1}}\left(
1\right) \right) \rightarrow \mathcal{I}_{\left( 1\right) }\otimes \mathcal{I%
}_{\left( 3\right) }\otimes O_{H_{1}}\left( 1\right) \otimes O_{H_{1}}\left(
1\right) .
\end{equation*}
We can use the same homomorphism with $3$ (and $1$) replacing $2$ to further
compose with $\mathcal{I}_{\left( 3\right) }\otimes \mathcal{I}_{\left(
1\right) }\otimes O_{H_{2}}\left( 1\right) $ and obtain a homomorphism from
the tensor product of the three components of $\mathcal{L}_{6}\left(
0\right) $ into $O_{H_{1}}\left( 1\right) \otimes O_{H_{1}}\left( 1\right)
\otimes O_{H_{2}}\left( 1\right) =O_{X_{6}\left( 0\right) }\left( 1\right) $%
. This induces the triple product on $\mathcal{L}_{6}\left( 0\right) $ in
terms of its individual components. Of course, this is complete analogous to
the description of the triple product on $\mathbf{L}_{6}$, 
\begin{equation*}
\mathbf{L}_{6}\otimes \mathbf{L}_{6}\otimes \mathbf{L}_{6}\rightarrow 
\mathbb{C}\text{.}
\end{equation*}
in terms of its decomposition as representation of $A_{2}\times A_{2}\times
A_{2}$: 
\begin{equation*}
\mathbf{L}_{6}|_{A_{2}\times A_{2}\times A_{2}}=\Lambda _{2}\otimes \Lambda
_{2}^{\ast }\otimes \mathbf{1}+\mathbf{1}\otimes \Lambda _{2}\otimes \Lambda
_{2}^{\ast }+\Lambda _{2}^{\ast }\otimes \mathbf{1}\otimes \Lambda _{2}.
\end{equation*}
In that case the triple product is given as follows: Let $s_{1}\otimes
t_{2}\in \Lambda _{2}\otimes \Lambda _{2}^{\ast }\otimes \mathbf{1\subset L}%
_{6}$,$s_{2}\otimes t_{3}\in \mathbf{1}\otimes \Lambda _{2}\otimes \Lambda
_{2}^{\ast }\subset \mathbf{L}_{6}$ and $s_{3}\otimes t_{1}\in \Lambda
_{2}^{\ast }\otimes \mathbf{1}\otimes \Lambda _{2}\subset \mathbf{L}_{6}$.
Then their product is 
\begin{equation*}
\left( s_{1}\otimes t_{2}\right) \otimes \left( s_{2}\otimes t_{3}\right)
\otimes \left( s_{3}\otimes t_{1}\right) \mapsto \left( t_{1},s_{1}\right)
\left( t_{2},s_{2}\right) \left( t_{3},s_{3}\right) .
\end{equation*}
Here $\left( t,s\right) $ is the natural pairing between $\Lambda _{2}^{\ast
}$ and $\Lambda _{2}$. And the product for all other combinations are zero
except those obtained by permutation. (see \cite{Segre})

Once the situation for $\mathcal{L}_{6}\left( 0\right) $ is settled, the
decomposition for $\mathcal{R}_{6}\left( 0\right) $ is immediately. It is
because if $l$ is a line on $X_{6}$, then $R=-l-K$ determines a ruling on
it. Hence 
\begin{eqnarray*}
\mathcal{R}_{6}\left( 0\right) &=&\mathcal{I}_{\left( 1\right) }\otimes 
\mathcal{I}_{\left( 2\right) }\otimes O_{H_{1}}\left( 1\right) \otimes
O_{H_{2}}\left( 1\right) \\
&&+\mathcal{I}_{\left( 2\right) }\otimes \mathcal{I}_{\left( 3\right)
}\otimes O_{H_{2}}\left( 1\right) \otimes O_{H_{3}}\left( 1\right) \\
&&+\mathcal{I}_{\left( 3\right) }\otimes \mathcal{I}_{\left( 1\right)
}\otimes O_{H_{3}}\left( 1\right) \otimes O_{H_{1}}\left( 1\right) .
\end{eqnarray*}
The triple product on $\mathcal{R}_{6}\left( 0\right) $ can be described in
similar manner.

It would be useful to have a rigorous description of the degeneration of $%
\mathcal{L}_{6}\left( t\right) $ and $\mathcal{R}_{6}\left( t\right) $ as $t$
goes to zero. Such description should also give the decomposition of $L%
\mathcal{E}_{6}\left( 0\right) $. We recall that the decomposition of $E_{6}$
under $A_{2}\times A_{2}\times A_{2}$ is as follow: 
\begin{equation*}
E_{6}=A_{2}+A_{2}+A_{2}+\Lambda _{2}\otimes \ \Lambda _{2}\otimes \ \Lambda
_{2}+\Lambda _{2}^{\ast }\otimes \Lambda _{2}^{\ast }\otimes \Lambda
_{2}^{\ast }\text{.}
\end{equation*}
The description of $E_{6}$ on $\mathbf{L}_{6}=\Lambda _{2}\otimes \Lambda
_{2}^{\ast }\otimes \mathbf{1}+\mathbf{1}\otimes \Lambda _{2}\otimes \Lambda
_{2}^{\ast }+\Lambda _{2}^{\ast }\otimes \mathbf{1}\otimes \Lambda _{2}$ in
terms of this decomposition is very simple. Each $A_{2}$ component in $E_{6}$
acts on the corresponding factor in each component of $\mathbf{L}_{6}$. To
describe the action of $\Lambda _{2}\otimes \ \Lambda _{2}\otimes \ \Lambda
_{2}$ on $\Lambda _{2}\otimes \Lambda _{2}^{\ast }\otimes \mathbf{1,}$ we
first notice that there is a homomorphism from the tensor product of $%
\Lambda _{2}$ with itself to its dual space. It is because, being the
standard representation of $A_{2}$, the third wedge power of $\Lambda _{2}$
is $\mathbf{1}$. Using this and the natural pairing between $\Lambda _{2}$
and $\Lambda _{2}^{\ast }$, we obtain the homomorphism 
\begin{equation*}
\left( \Lambda _{2}\otimes \ \Lambda _{2}\otimes \ \Lambda _{2}\right)
\bigotimes \left( \Lambda _{2}\otimes \Lambda _{2}^{\ast }\otimes \mathbf{1}%
\right) \rightarrow \left( \Lambda _{2}^{\ast }\otimes \mathbf{1}\otimes
\Lambda _{2}\right) ,
\end{equation*}
which determines the action of $\Lambda _{2}\otimes \ \Lambda _{2}\otimes \
\Lambda _{2}$ on $\Lambda _{2}\otimes \Lambda _{2}^{\ast }\otimes \mathbf{1}$%
. The action of $\Lambda _{2}\otimes \ \Lambda _{2}\otimes \ \Lambda _{2}$
on other components is similar. Also the action of $\Lambda _{2}^{\ast
}\otimes \Lambda _{2}^{\ast }\otimes \Lambda _{2}^{\ast }$ on $\mathbf{L}%
_{6} $ is analogous. This gives the complete description of $E_{6}$ on $%
\mathbf{L}_{6}$ in this decomposition.

\subsection{Reduction to $\mathbf{sl}\left( 2\right) \times \mathbf{sl}%
\left( 6\right) $ bundle}

In this section we degenerate the cubic surface into union of a plane $H$
and a smooth quadric surface $Q$ in $\mathbb{P}^{3}$. We are going to see
that the symmetry group of $\mathcal{L}_{6},\mathcal{R}_{6}$ and $L\mathcal{E%
}_{6}$ would be reduced from $E_{6}$ to $\mathbf{sl}\left( 6\right) \times 
\mathbf{sl}\left( 2\right) =A_{5}\times A_{1}$.

Let $X\left( t\right) $ be the family of cubic surfaces in $\mathbb{P}^{3}$
parametrized by $t$. We assume $X\left( 0\right) =H\cup Q$ and $X\left(
t\right) $ is smooth for nonzero $t$. We denote $C=H\cap Q$ and $Z=\left\{
p_{1},p_{2},p_{3},p_{4},p_{5},p_{6}\right\} $ consists of those points where 
$C$ intersects $X\left( t\right) $ infinitesimally as we various $t$ away
from zero.

Now a line in $H$ which can be deformed to one in $X\left( t\right) $ for
small $t$ if and only if it intersects $C$ at two of the six points in $Z$.
On $Q$, there are two rulings. A member in one of the two rulings on $Q$ can
be deformed to a line in $X\left( t\right) $ for small $t$ if and only if it
intersects $C$ at one of the six points in $Z$. Moreover any one parameter
family of lines on $X\left( t\right) $ parametrized by $t$ would converge to
one of these lines. The total number of these lines equals $27=\binom{6}{2}%
+6\times 2$.

In fact the decomposition of the representation $\mathbf{L}_{6}$ of $E_{6}$
under $A_{5}\times A_{1}$ has the similar structure. If we denote the
standard representation of $A_{n}$ by $\Lambda _{n}$ and the fundamental
representation of $A_{n}$ given by the second wedge product of its standard
representation by $\Lambda _{n}^{2}$. Then 
\begin{equation*}
\mathbf{L}_{6}|_{A_{5}\times A_{1}}=\Lambda _{5}^{2}\otimes \mathbf{1}%
+\Lambda _{5}\otimes \Lambda _{2}.
\end{equation*}

We define the coherent sheaf $\mathcal{I}=\oplus _{p\in Z}\mathcal{I}%
_{\left\{ p\right\} }$ and we denote the hyperplane bundle on $H$ by $%
O_{H}\left( 1\right) $. We treat $O_{H}\left( 1\right) $ as a rank one
coherent sheaf on $X_{6}\left( 0\right) $ whose restriction on $X_{6}\left(
0\right) \backslash H$ being trivial. Moreover we define the rank two vector
bundle 
\begin{equation*}
\mathcal{R}_{Q}=\tbigoplus_{_{\substack{ R^{2}=0  \\ RK_{Q}=-2}}}O_{Q}\left(
R\right) ,
\end{equation*}
on $Q$ in terms of it rulings. We also treat $\mathcal{R}_{Q}$ as a rank two
coherent sheaf on $X_{6}\left( 0\right) $ whose restriction on $X_{6}\left(
0\right) \backslash Q$ being trivial. Therefore it is reasonable to propose
that 
\begin{equation*}
\mathcal{L}_{6}\left( 0\right) =\Lambda ^{2}\mathcal{I}\otimes O_{H}\left(
1\right) +\mathcal{I}\otimes \mathcal{R}_{Q}.
\end{equation*}
In order to describe the triple product on $\mathcal{L}_{6}\left( 0\right) $%
, we consider several homomorphisms. The first one is a homomorphism from $%
\bigotimes^{3}\left( \Lambda ^{2}\mathcal{I}\right) $ to $O_{X_{6}\left(
0\right) }$ defined by the composition of naturally defined homomorphisms: 
\begin{equation*}
\bigotimes^{3}\left( \Lambda ^{2}\mathcal{I}\right) \rightarrow \Lambda
^{3}\left( \Lambda ^{2}\mathcal{I}\right) \rightarrow \Lambda ^{6}\mathcal{I}%
=\mathcal{I}_{Z}\rightarrow O_{X_{6}\left( 0\right) }.
\end{equation*}
The second one is a homomorphism from $\left( \Lambda ^{2}\mathcal{I}\right)
\otimes \mathcal{I}\otimes \mathcal{I}$ to $O_{X_{6}\left( 0\right) }$
defined by

\begin{equation*}
\left( \oplus _{p\neq q\in Z}s_{pq}\right) \otimes \left( \oplus _{p\in
Z}s_{p}^{\prime }\right) \otimes \left( \oplus _{p\in Z}s_{p}^{\prime \prime
}\right) \mapsto \sum_{p\neq q\in Z}s_{pq}\left( s_{p}^{\prime
}s_{q}^{\prime \prime }+s_{q}^{\prime }s_{p}^{\prime \prime }\right) .
\end{equation*}
The last one is the isomorphism between $\Lambda ^{2}\mathcal{R}_{Q}$ and $%
O_{Q}\left( 1\right) $. We leave the details for the construction of the
triple product to our readers. The construction of $\mathcal{R}_{6}\left(
0\right) $ and its triple product is also similar (see the previous section).

\newpage

\section{Branched cover $X_{7}$ of $\mathbb{P}^{2}$ along a quartic curve
and its $E_{7}$-bundle}

The quartic plane curve comes into our picture because every degree two Del
Pezzo surface $X$ is a double cover of $\mathbb{P}^{2}$ branched along a
quartic plane curve $C$. It is a famous old fact that there are 28
bitangents to $C$. This can be seen either by analyzing the dual curve or by
identifying each bitangent with an odd theta characteristic where there are
28 of them. (\cite{Clemens})

\subsection{Degree two del Pezzo surface and its $E_{7}$ bundle}

\begin{center}
\textbf{The representation bundle }$\mathcal{R}_{7}$
\end{center}

On $X_{7}$ we have $\mathcal{L}_{7}=\bigoplus_{\substack{ l^{2}=-1  \\ %
l\cdot K=-1}}O\left( l\right) $ and $\mathcal{R}_{7}=O\left( -K\right)
^{\oplus 7}+\bigoplus_{\substack{ R^{2}=0  \\ R\cdot K=-2}}O\left( R\right)
. $ Notice that the extra summand $O\left( -K\right) ^{\oplus 7}$ in $%
\mathcal{R}_{7}$ does not correspond to rulings on $X_{7}$. We now explain
the appearance of this summand from two different points of view. First we
consider the representation bundles homomorphism 
\begin{equation*}
\mathcal{L}_{n}\otimes \mathcal{L}_{n}\rightarrow \mathcal{R}_{n}\text{.}
\end{equation*}
For $n\leq 6,$ the construction of this homomorphism is related to the fact
that if $l_{1}$ and $l_{2}$ are distinct lines on $X_{6}$, then $l_{1}+l_{2}$
determines a ruling on $X_{6}$ provided that $l_{1}\cdot l_{2}=1$. The
anti-canonical class of $X_{n}$ is very ample and defines an embedding of $%
X_{n}$ inside $\mathbb{P}^{9-n}$. Therefore $l_{1}\cdot l_{2}$ is either
zero or one. When $n=7$, the anti-canonical line bundle is no longer very
ample and its sections only defines $X_{7}$ as a double cover of $\mathbb{P}%
^{2}$. Therefore $l_{1}\cdot l_{2}$ can be zero, one or two. Explicitly, if
we represent $X_{7}$ as blowup of $\mathbb{P}^{2}$ at seven points with
exceptional locus $\bigcup_{i=1}^{7}L_{i}$. Then $l_{1}\equiv
3H-\sum_{i=1}^{7}L_{i}-L_{j}$ and $l_{2}=L_{j}$ are two distinct lines on $%
X_{7}$ with $l_{1}\cdot l_{2}=2$. When this happens, we have $O\left(
l_{1}+l_{2}\right) =O\left( -K\right) $ for $j=1,...,7$. This explains why
we need the extra factor $O\left( -K\right) ^{\oplus 7}$ in $\mathcal{R}_{7}$%
.

Second point of view: Such an extra factor can also be explained in terms of
Lie algebra theory. As a representation of $E_{7}$, the adjoint
representation of $E_{7}$ coincides with $\mathbf{R}_{7}$. Therefore we
expect $L\mathcal{E}_{7}$ and $\mathcal{R}_{7}$ to be equal up to tensoring
with a line bundle over $X_{7}$. Recall that $L\mathcal{E}%
_{7}=O_{X_{7}}^{\oplus 7}+O\left( D\right) $ where we sum over divisors $D$
with $D^{2}=-2$ and $D\cdot K=0$. For any such $D$, we have $\left(
D-K\right) ^{2}=0$ and $\left( D-K\right) \cdot K=-2$ because $K^{2}=2$.
This implies that $D-K$ determines a ruling on $X_{7}$. That is 
\begin{equation*}
L\mathcal{E}_{7}\otimes O\left( -K\right) =\mathcal{R}_{7}\text{.}
\end{equation*}
So the $O_{X_{7}}^{\oplus 7}$ summand in $L\mathcal{E}_{7}$ turns into $%
O\left( -K\right) ^{\oplus 7}$summand in $\mathcal{R}_{7}$.

\begin{center}
\textbf{Quadratic structure on }$\mathcal{L}_{7}$
\end{center}

For any line $l$ on $X_{7}$, using $K^{2}=2$, we have $\left( -l-K\right)
^{2}=-1$ and $\left( -l-K\right) \cdot K=-1$. That is, there is another line 
$l^{\prime }$ on $X_{7}$ which is linearly equivalent to $-l-K$. Using
isomorphisms $O\left( l\right) \otimes O\left( l^{\prime }\right) \overset{%
\cong }{\rightarrow }O\left( -K\right) $ with various line $l$, we obtain a
fiberwise non-degenerate quadratic form on $\mathcal{L}_{7}$, 
\begin{equation*}
q_{7}:\mathcal{L}_{7}\otimes \mathcal{L}_{7}\rightarrow O\left( -K\right) 
\text{.}
\end{equation*}
This quadratic form is preserved by the action of $L\mathcal{E}_{7}$.

Notice that $l\cdot l^{\prime }=l\cdot \left( -l-K\right) =2$ and $l^{\prime
}$ can be characterized as the only line with this property.

\begin{proposition}
If $l$ is a line on $X_{7}$, then there is exactly one line $l^{\prime }$ on 
$X_{7}$ such that $l\cdot l^{\prime }=2.$ Moreover $l^{\prime }\cong -l-K$.
\end{proposition}

Proof of proposition: If $l$ is a line on $X_{7}$, using $K^{2}=2$, we get $%
\left( -l-K\right) ^{2}=\left( -l-K\right) \cdot K=-1$. Hence there is a
line $l^{\prime }$ on $X_{7}$ linearly equivalent to $-l-K$. Moreover $%
l\cdot l^{\prime }=2$.

On the other hand if $l_{1}$ is any line with $l\cdot l_{1}=2$, then 
\begin{equation*}
l_{1}\cdot l^{\prime }=l_{1}\cdot \left( -l-K\right) =-1\text{.}
\end{equation*}
Since both $l_{1}$ and $l^{\prime }$ are irreducible effective divisors on $%
X_{7}$, their intersection must be nonnegative unless they are the same. So
we get uniqueness result. $\square $

Because of this proposition, the set of 56 lines on $X_{7}$ is naturally
divided into 28 pairs of lines. If we consider $X_{7}$ as a double cover of $%
\mathbb{P}^{2}$, 
\begin{equation*}
\delta :X_{7}\rightarrow \mathbb{P}^{2}.
\end{equation*}
The branched locus $B$ of $\delta $ is a quartic curve on $\mathbb{P}^{2}$.
It is classically known that there are 28 bitangents to $B$. The inverse
image of each such bitangent consists precisely a pair of lines $l,l^{\prime
}$ on $X_{7}$ with $l\cdot l^{\prime }=2$. The above quadratic form $q_{7}$
is just reflecting this property because $l+l^{\prime }=-K_{X}$ which is the
pullback of the hyperplane bundle of $\mathbb{P}^{2}$ via $\delta $: $%
K_{X}^{-1}=\delta ^{\ast }\left( K_{\mathbb{P}^{2}}\otimes O\left( 2\right)
\right) ^{-1}=\delta ^{\ast }O_{\mathbb{P}^{2}}\left( 1\right) $.

\begin{remark}
Consider the action of $L\mathcal{E}_{7}$ on $\mathcal{L}_{7}$ given by the
homomorphism $L\mathcal{E}_{7}\otimes \mathcal{L}_{7}\rightarrow \mathcal{L}%
_{7}.$ If we use the quadratic form $q_{7}$ to identify $\mathcal{L}_{7}$
with $\mathcal{L}_{7}^{\ast }\otimes O\left( -K\right) $ and the fact that $L%
\mathcal{E}_{7}\cong \left( L\mathcal{E}_{7}\right) ^{\ast }$, then the
action of $L\mathcal{E}_{7}$ on $\mathcal{L}_{7}$ gives rise to $\mathcal{L}%
_{7}\otimes \mathcal{L}_{7}\rightarrow L\mathcal{E}_{7}\otimes O\left(
-K\right) $. Using the isomorphism between $L\mathcal{E}_{7}\otimes O\left(
-K\right) $ and $\mathcal{R}_{7},$ we obtain $\mathcal{L}_{7}\otimes 
\mathcal{L}_{7}\rightarrow \mathcal{R}_{7}$ which is just the homomorphism $%
c_{7}$.
\end{remark}

\begin{center}
\textbf{Quartic structure on }$\mathcal{L}_{7}$
\end{center}

Next we want to describe a symmetric quartic form on $\mathcal{L}_{7}$: 
\begin{equation*}
f_{7}:\mathcal{L}_{7}\otimes \mathcal{L}_{7}\otimes \mathcal{L}_{7}\otimes 
\mathcal{L}_{7}\rightarrow O\left( -2K\right)
\end{equation*}
which is invariant under the action of $L\mathcal{E}_{7}$. As we have
mentioned in the introduction, (the identity component of) the bundle of
automorphisms of $\mathcal{L}_{7}$ preserving $q_{7}$ and $f_{7}$ determined
a $E_{7}$-bundle $\mathcal{E}_{7}$ over $X_{7}$ whose associated Lie algebra
bundle is just $L\mathcal{E}_{7}$.

To understand this, we notice that this quartic form $f_{7}$ globalizing a
corresponding quartic form $f$ on $\mathbf{L}_{7}$ which is invariant under $%
E_{7}$ (see chapter 12 of \cite{Adams}). Moreover the identity component of
the group of automorphism of $\mathbf{L}_{7}$ preserving it and the Killing
form is a simply connected Lie group of type $E_{7}$.

We are going to give two different descriptions of $f_{7}$.

First description: Given three distinct lines $l_{1},l_{2}$ and $l_{3}$ on $%
X_{7}$, we have 
\begin{equation*}
\left( -l_{1}-l_{2}-l_{3}-2K\right) \cdot K=-1,
\end{equation*}
and 
\begin{equation*}
\left( -l_{1}-l_{2}-l_{3}-2K\right) ^{2}=-7+2\left( l_{1}\cdot
l_{2}+l_{2}\cdot l_{3}+l_{1}\cdot l_{3}\right) .
\end{equation*}
Hence there is a line $l_{4}$ in the linearly equivalent class of $%
-l_{1}-l_{2}-l_{3}-2K$ if and only if $l_{1}\cdot l_{2}+l_{2}\cdot
l_{3}+l_{1}\cdot l_{3}=3$. We have $0\leq l_{i}\cdot l_{j}\leq 2$ for $i\neq
j$. So, up to permutation of their indexes, there are only two possibilities
for us to have $l_{1}\cdot l_{2}+l_{2}\cdot l_{3}+l_{1}\cdot l_{3}=3$.
Namely, case (1): $l_{1}\cdot l_{2}=l_{2}\cdot l_{3}=l_{1}\cdot l_{3}=1$ and
case (2) $l_{1}\cdot l_{2}=2,l_{1}\cdot l_{3}=1,l_{2}\cdot l_{3}=0$. In case
(1), we also have $l_{i}\cdot l_{4}=1$ for $i=1,2,3.$ That is $l_{i}\cdot
l_{j}=1$ for all $i\neq j$ between $1$ and $4$. In case (2), we also have $%
l_{3}\cdot l_{4}=2,l_{2}\cdot l_{4}=1$ and $l_{1}\cdot l_{4}=0.$

In either cases, we have isomorphism 
\begin{equation*}
O\left( l_{1}\right) \otimes O\left( l_{2}\right) \otimes O\left(
l_{3}\right) \overset{\cong }{\rightarrow }O\left( -l_{4}\right) \otimes
O\left( -2K\right) ,
\end{equation*}
or equivalently 
\begin{equation*}
O\left( l_{1}\right) \otimes O\left( l_{2}\right) \otimes O\left(
l_{3}\right) \otimes O\left( l_{4}\right) \overset{\cong }{\rightarrow }%
O\left( -2K\right) .
\end{equation*}
Combining all such homomorphisms with $l_{1},l_{2},l_{3}$ satisfying either
(1) or (2) and with $l_{4}\equiv -l_{1}-l_{2}-l_{3}-2K$, we obtain a
homomorphism 
\begin{equation*}
f_{7}:\mathcal{L}_{7}\otimes \mathcal{L}_{7}\otimes \mathcal{L}_{7}\otimes 
\mathcal{L}_{7}\rightarrow O\left( -2K\right) .
\end{equation*}
By construction, it is clear that $f_{7}$ is symmetric. One can also verify
directly that $f_{7}$ is invariant under the action of $L\mathcal{E}_{7}$.
Moreover, on each fiber of $\mathcal{L}_{7}$, the restriction of $f_{7}$ to
it is isomorphic to $f$ on $\mathbf{L}_{7}$. Therefore (the identity
component of) the bundle of automorphisms of $\mathcal{L}_{7}$ preserving
the quadratic form $q_{7}$ and the quartic form $f_{7}$ is a principal
bundle $\mathcal{E}_{7}$ of type $E_{7}$. This also implies that its
associated Lie algebra bundle is just $L\mathcal{E}_{7}$ because of the fact
that both $q_{7}$ and $f_{7}$ are invariant under the action of $L\mathcal{E}%
_{7}$.

Second description: Recall that, for any $n\leq 8$, we have a representation
bundles homomorphism: $c_{n}:\mathcal{L}_{n}\otimes \mathcal{L}%
_{n}\rightarrow \mathcal{R}_{n}$. However, by the previous section, we also
have an isomorphism $\mathcal{R}_{7}=L\mathcal{E}_{7}\otimes O\left(
-K\right) $. Combining these two, we have the following homomorphism: 
\begin{equation*}
\mathcal{L}_{7}\otimes \mathcal{L}_{7}\rightarrow L\mathcal{E}_{7}\otimes
O\left( -K\right) .
\end{equation*}
Also $\mathcal{L}_{7}$ is a representation bundle of $L\mathcal{E}_{7}$, so
it gives another homomorphism 
\begin{equation*}
L\mathcal{E}_{7}\otimes \mathcal{L}_{7}\rightarrow \mathcal{L}_{7}.
\end{equation*}
On the other hand, there is an isomorphism 
\begin{equation*}
\mathcal{L}_{7}\cong \mathcal{L}_{7}^{\ast }\otimes O\left( -K\right) ,
\end{equation*}
given by the quadratic form $q_{7}$. Combining these, we obtain a
homomorphism 
\begin{equation*}
\mathcal{L}_{7}\otimes \mathcal{L}_{7}\otimes \mathcal{L}_{7}\rightarrow 
\mathcal{L}_{7}^{\ast }\otimes O\left( -2K\right) ,
\end{equation*}
which is equivalent to the above quartic form 
\begin{equation*}
f_{7}:\mathcal{L}_{7}\otimes \mathcal{L}_{7}\otimes \mathcal{L}_{7}\otimes 
\mathcal{L}_{7}\rightarrow O\left( -2K\right) .
\end{equation*}
We leave the verifications to our readers. In this description, it is not
obvious that $f_{7}$ is symmetric.

\subsection{Reduction to $\mathbf{so}\left( 12\right) \times \mathbf{sl}%
\left( 2\right) $-bundle}

We recall that $\mathbf{so}\left( 12\right) \times \mathbf{sl}\left(
2\right) =D_{6}\times A_{1}$.

Now we fix a ruling $R$ on $X_{7}$. A line $l$ on $X_{7}$ would have $l\cdot
R$ equals zero, one or two. However, if $l\cdot R=0$ then there is another
line $l^{\prime }$ on $X_{7}$ with $l^{\prime }\cdot R=2$. In fact $%
l^{\prime }\equiv l-R-K$ which can be verified directly using $K^{2}=2$.
This give 
\begin{equation*}
\mathcal{L}_{7}=\mathcal{S}^{+}+\mathcal{W}_{6}\otimes \Lambda _{1}\text{.}
\end{equation*}
Here $\Lambda _{1}$ is the rank two vector bundle $O_{X_{7}}+O_{X_{7}}\left(
-R-K\right) $%
\begin{equation*}
\Lambda _{1}=O_{X_{7}}+O_{X_{7}}\left( -R-K\right) \text{.}
\end{equation*}
The automorphism bundle of $\Lambda _{1}$ preserving $\det :\det \Lambda _{1}%
\overset{\cong }{\rightarrow }O\left( -R-K\right) $ is a $SL\left( 2\right) $%
-bundle over $X_{7}$ which we call $\mathcal{A}_{1}^{X_{7}}$, or simply $%
\mathcal{A}_{1}$. Its corresponding $\mathbf{sl}\left( 2\right) $-Lie
algebra bundle will be called $L\mathcal{A}_{1}^{X_{7}}$, or simply $L%
\mathcal{A}_{1}$. Then the above decomposition is a decomposition of $L%
\mathcal{D}_{6}+L\mathcal{A}_{1}$-representation bundles corresponding to
the following decomposition of $E_{7}$ representation: 
\begin{equation*}
\mathbf{L}_{7}|_{D_{6}\times A_{1}}=S^{+}\otimes \mathbf{1}+W_{6}\otimes
\Lambda _{1}\text{,}
\end{equation*}
where $\Lambda _{1}$ is the standard representation of $A_{1}=\mathbf{sl}%
\left( 2\right) $.

Next we want to decompose $L\mathcal{E}_{7}$ under $L\mathcal{D}_{6}+L%
\mathcal{A}_{1}$. First we have 
\begin{eqnarray*}
L\mathcal{A}_{1} &=&End_{0}\left( O+O\left( -R-K\right) \right) \\
&=&O\left( R+K\right) +O+O\left( -R-K\right) .
\end{eqnarray*}
If we write $D=R+K$ or $-R-K$, then we have $D^{2}=-2$ and $D\cdot K=0$.
That is $L\mathcal{A}_{1}$ is a vector subbundle of $L\mathcal{E}_{7}$.
Recall that $\mathcal{S}^{-}=\bigoplus O\left( T\right) $ with $T$
satisfying $T^{2}=-2,T\cdot K=0$ and $T\cdot R=1$, is always a subbundle of $%
L\mathcal{E}_{7}$. By direct computations, we have 
\begin{equation*}
\mathcal{S}^{-}\otimes O\left( R+K\right) =\bigoplus_{\substack{ R^{2}=-2 
\\ RK=0  \\ RR=-1}}O\left( R\right) ,
\end{equation*}
and therefore also a subbundle of $L\mathcal{E}_{7}$. Moreover 
\begin{equation*}
L\mathcal{E}_{7}=L\mathcal{D}_{6}+L\mathcal{A}_{1}+\mathcal{S}^{-}+\mathcal{S%
}^{-}\otimes O\left( R+K\right) .
\end{equation*}

As a decomposition of $L\mathcal{D}_{6}+L\mathcal{A}_{1}$ representation
bundles, we have 
\begin{equation*}
L\mathcal{E}_{7}=L\mathcal{D}_{6}+L\mathcal{A}_{1}+\mathcal{S}^{-}\otimes
\Lambda _{1}^{\ast }.
\end{equation*}
This corresponds to the following decomposition of $E_{7}$-representation
under $D_{6}\times A_{1}:$%
\begin{equation*}
E_{7}|_{D_{6}\times A_{1}}=D_{6}+A_{1}+S^{-}\otimes \Lambda _{1}^{-1}\text{.}
\end{equation*}
We obtain a similar decomposition for $\mathcal{R}_{7}$ because $\mathcal{R}%
_{7}=L\mathcal{E}_{7}\otimes O\left( -K\right) $. Namely 
\begin{equation*}
\mathcal{R}_{7}=O\left( R\right) \left( S^{2}\Lambda _{1}+\mathcal{S}%
^{-}\otimes \Lambda _{1}+\Lambda ^{2}\mathcal{W}_{6}\otimes \left(
-K-2R\right) \right) \text{.}
\end{equation*}

\subsection{Reduction to $\mathbf{sl}\left( 8\right) $ bundle}

In this section we degenerate $X_{7}$ into a normal crossing variety which
consists of two copies of $\mathbb{P}^{2}$ joining along a conic curve. We
shall see that the $E_{7}$ structure on $\mathcal{L}_{7}$ and $\mathcal{R}%
_{7}$ reduces to $SL\left( 8\right) =A_{7}$ structure under such
degeneration.

To begin we consider $X_{7}$ as a double cover of $\mathbb{P}^{2}$ branched
along a quartic plane curve. Now we want deform the quartic curve into a
double conic. Let $t$ be the deformation parameter and $B\left( t\right) $
be a family of quartic plane curve with $B\left( 0\right) =2C$ with $C$ is
smooth conic in $\mathbb{P}^{2}$. We also assume that $B\left( t\right) $ is
smooth when $t$ is not zero. Let us denote the double cover of $\mathbb{P}%
^{2}$ branched along $B\left( t\right) $ as 
\begin{equation*}
\delta \left( t\right) :X\left( t\right) \rightarrow \mathbb{P}^{2}.
\end{equation*}
In particular $X\left( 0\right) $ consists of two copies of $\mathbb{P}^{2}$
joining along $C$. We call the two $\mathbb{P}^{2}$ as $\mathbb{P}_{\left(
1\right) }^{2}$ and $\mathbb{P}_{\left( 2\right) }^{2}$. We let $Z=\left\{
p_{1},p_{2},...,p_{8}\right\} \subset C$ be the set consisting of points
where $C$ meets $B\left( t\right) $ as we various $t$ away from zero
infinitesimally.

Now a line in $\mathbb{P}_{\left( i\right) }^{2}$ which can be deformed to
one in $X\left( t\right) $ for small $t$ if and only if it intersects $C$ at
two of the eight points in $Z$. Moreover any one parameter family of lines
on $X\left( t\right) $ parametrized by $t$ would converge to one of these
lines. The total number of these lines equals 
\begin{equation*}
56=\left( 
\begin{array}{c}
8 \\ 
2
\end{array}
\right) +\left( 
\begin{array}{c}
8 \\ 
2
\end{array}
\right) .
\end{equation*}
In fact these 28 line on $\mathbb{P}^{2}$ corresponds precisely to the limit
of the 28 bitangents to the quartic plane curve on $B\left( t\right) $ (see 
\cite{Clemens}) We saw in section two, that if $l$ is a line on $X_{7}$ then
there is exactly one other line $l^{\prime }$ which intersect $l$ at two
points and moreover $l^{\prime }$ is in the linearly equivalent class of $%
-l-K$. This divides the 56 lines into 28 pairs. In our situation here, each
pair consists of one line in $\mathbb{P}_{\left( 1\right) }^{2}$ and one
line in $\mathbb{P}_{\left( 2\right) }^{2}$.

The decomposition of the representation $\mathbf{L}_{7}$ of $E_{7}$ under $%
A_{7}$ has the similar structure. If we denote the standard representation
of $A_{n}$ by $\Lambda _{n}$ and the fundamental representation of $A_{n}$
given by the $k^{th}$ wedge product of its standard representation by $%
\Lambda _{n}^{k}$. Then 
\begin{equation*}
\mathbf{L}_{7}|_{A_{7}}=\Lambda _{7}^{2}+\Lambda _{2}^{6}.
\end{equation*}

We define the coherent sheaf $\mathcal{I}=\oplus _{p\in Z}\mathcal{I}%
_{\left\{ p\right\} }$ and we denote the hyperplane bundle on $\mathbb{P}%
_{\left( i\right) }^{2}$ by $O_{\mathbb{P}_{\left( i\right) }^{2}}\left(
1\right) $. We treat $O_{\mathbb{P}_{\left( i\right) }^{2}}\left( 1\right) $
as a rank one coherent sheaf on $X\left( 0\right) $ whose restriction on $%
X\left( 0\right) \backslash \mathbb{P}_{\left( i\right) }^{2}$ being
trivial. Therefore it is natural to define 
\begin{equation*}
\mathcal{L}_{7}\left( 0\right) =\Lambda ^{2}\mathcal{I}\otimes O_{\mathbb{P}%
_{\left( 1\right) }^{2}}\left( 1\right) +\Lambda ^{2}\mathcal{I}\otimes O_{%
\mathbb{P}_{\left( 2\right) }^{2}}\left( 1\right) .
\end{equation*}

We now also describe the quartic form on $\mathcal{L}_{7}\left( 0\right) $, 
\begin{equation*}
f_{7}:\mathcal{L}_{7}\left( 0\right) \otimes \mathcal{L}_{7}\left( 0\right)
\otimes \mathcal{L}_{7}\left( 0\right) \otimes \mathcal{L}_{7}\left(
0\right) \rightarrow O\left( -2K\right) .
\end{equation*}
In terms of lines on $X_{7}\left( 0\right) $, nontrivial product in $f_{7}$
corresponds to either (1) four lines in one of the $\mathbb{P}_{\left(
i\right) }^{2}$'s and they pass though all eight points of $Z$ or (2) two
pairs of lines.

As in other cases of degeneration into nonnormal surfaces, it is important
to have a better understanding of the degeneration of $\mathcal{L}_{7}\left(
t\right) $ (and $\mathcal{R}_{7}\left( t\right) $) and compare the limit
with $\mathcal{L}_{7}\left( 0\right) $.\newpage

\section{$E_{8}$-bundles over $X_{8}$}

Among exceptional Lie algebra, $E_{8}$ has many exceptional properties not
share by other $E_{n}$'s. These properties can also be seen on the geometry
of $X_{8}$. For example the number of lines on $X_{8}$ equals 240 which is
different from the dimension of $\mathbf{L}_{8}$ which equals 248. In fact 
\begin{equation*}
\mathcal{L}_{8}=O\left( -K\right) ^{\oplus 8}\bigoplus_{\substack{ l^{2}=-1 
\\ lK=-1}}O\left( l\right) .
\end{equation*}

To explain this extra summand $O\left( -K\right) ^{\oplus 8}$, we first
notice that if $l$ and $l^{\prime }$ are lines on $X_{8}$, then $\left(
l+l^{\prime }+K\right) ^{2}=-5+2l\cdot l^{\prime }$ and $\left( l+l^{\prime
}+K\right) \cdot K=-1$. When $l\cdot l^{\prime }=2$, there is a line in the
linearly equivalent class of $l+l^{\prime }+K$. This suggests a product
structure on $\bigoplus_{\substack{ l^{2}=-1  \\ lk=-1}}O\left( l\right) $
and $\mathcal{L}_{8}$: 
\begin{equation*}
\mathcal{L}_{8}\otimes \mathcal{L}_{8}\rightarrow \mathcal{L}_{8}\otimes
O\left( -K\right) \footnote{%
This product is symmetric. However if we suitably introduce signs to varies
components of this homomorphism, we can make it coincides with the Lie
algebra structure on \ $E_{8}$.}.
\end{equation*}
In fact $\mathcal{L}_{8}$ is isomorphic to the bundle $L\mathcal{E}_{8}$ up
to tensoring with a line bundle: 
\begin{equation*}
\mathcal{L}_{8}=L\mathcal{E}_{8}\otimes O\left( -K\right) .
\end{equation*}
Using $K^{2}=1$, we obtain $\left( l+K\right) ^{2}=-1$ and $\left(
l+K\right) \cdot K=0$. It implies the above isomorphism as vector bundles.
In fact, one can check that they are isomorphic as representation bundles.
In terms of representations of $E_{8}$, it corresponds simply to the fact
that $\mathbf{L}_{8}$ coincides with the adjoint representation of $E_{8}$.
This also explains the extra summand $O\left( -K\right) ^{\oplus 8}$ in $%
\mathcal{L}_{8}$. The representation bundle $\mathcal{R}_{8}$ over $X_{8}$
is more complicated and we omit the discussion of it except to mention that
its rank equal 3875.

Recall that, in the introduction of this paper, we discussed the seven
types, (i),...,(vii), of lines on $X_{8}$ when we write $X_{8}$ as the
blowup of $\mathbb{P}^{2}$ at eight points. These various types can be
described in a more uniform way if we use the isomorphism $\mathcal{L}_{8}=L%
\mathcal{E}_{8}\otimes O\left( -K\right) $ and the earlier description of $L%
\mathcal{E}_{8}$ in terms of $H$ and various $L_{i}$'s. Here we denote the
exceptional locus of the blowup as $L_{1}\cup ...\cup L_{8}$ and the
pullback of the hyperplane class of $\mathbb{P}^{2}$ as $H$. In the
following table we display the various component types of $\mathcal{L}_{8}$
(except $O\left( -K\right) ^{\oplus 8}$), the corresponding type of lines
and the number of such lines. 
\begin{equation*}
\begin{tabular}{||lll||}
\hline\hline
$\left( L_{i}-L_{j}\right) -K$ & (iv) & 56 \\ \hline
$\left( H-L_{i}-L_{j}-L_{k}\right) -K$ & (v) & 56 \\ \hline
$\left( 2H-\sum_{m=1}^{6}L_{i_{m}}\right) -K$ & (vi) & 28 \\ \hline
$\left( 3H-\sum_{j=1}^{8}L_{j}-L_{i}\right) -K$ & (vii) & 8 \\ \hline\hline
\end{tabular}
\begin{tabular}{||lll||}
\hline\hline
&  &  \\ \hline
$\left( -H+L_{i}+L_{j}+L_{k}\right) -K$ & (iii) & 56 \\ \hline
$\left( -2H+\sum_{m=1}^{6}L_{i_{m}}\right) -K$ & (ii) & 28 \\ \hline
$\left( -3H+\sum_{j=1}^{8}L_{j}+L_{i}\right) -K$ & (i) & 8 \\ \hline\hline
\end{tabular}
\end{equation*}
with

\begin{tabular}{ll}
(i) & $\pi \left( D\right) =p_{i}$ \\ 
(ii) & $\pi \left( D\right) $ is a line passes through $p_{i}$ and $p_{j}$
\\ 
(iii) & $\pi \left( D\right) $ is a conic passes through five of the $p_{i}$%
's \\ 
(iv) & $\pi \left( D\right) $ is a cubic passes through seven of the $p_{i}$%
's and with one being a double point \\ 
(v) & $\pi \left( D\right) $ is a quartic passes through 8 of the $p_{i}$'s
and three being double points \\ 
(vi) & $\pi \left( D\right) $ is a quintic passes through 8 of the $p_{i}$'s
and six being double points \\ 
(vii) & $\pi \left( D\right) $ is a sextic passes through 8 of the $p_{i}$'s
\\ 
& and seven being double points and one triple point.
\end{tabular}

\subsection{Reduction to $E_{7}\times \mathbf{sl}\left( 2\right) $-bundle}

If we fix a line $L$ on $X_{8}$, then blowing down $L$ gives us a morphism $%
\pi :X_{8}\rightarrow X_{7}$ and the relationship between the $E_{8}$-bundle
on $X_{8}$ and $E_{7}$-bundle on $X_{7}$ has been discussed in section two
and we get

\begin{align*}
L\mathcal{E}_{8}& =\pi ^{\ast }L\mathcal{E}_{7}+O+\pi ^{\ast }\mathcal{L}%
_{7}\otimes O\left( -L\right) +\pi ^{\ast }\mathcal{L}_{7}^{\ast }\otimes
O\left( L\right) \\
& +O\left( -K-L\right) +O\left( K+L\right) .
\end{align*}
Now we consider the rank two vector bundle $\Lambda _{1}=O+O\left(
L+K\right) $ over $X_{8}$. Notice that $\Lambda _{1}\otimes O\left(
-K\right) $ is a subbundle of $\mathcal{L}_{8}$. Let $\mathcal{A}%
_{1}^{X_{8}} $, or simply $\mathcal{A}_{1}$, be the automorphism bundle of $%
\Lambda _{1}$ preserving $\det :\det \Lambda _{1}\overset{\cong }{%
\rightarrow }O\left( L+K\right) $. Then $\mathcal{A}_{1}$ is an $A_{1}$%
-bundle (or $\mathbf{sl}\left( 2\right) $-bundle). The corresponding Lie
algebra bundle is simply 
\begin{equation*}
L\mathcal{A}_{1}=O\left( -L-K\right) +O+O\left( L+K\right) ,
\end{equation*}
which is a Lie subalgebra bundle of $L\mathcal{E}_{8}$.

Next we use the quadratic form $q_{7}$ on $\mathcal{L}_{7}$: 
\begin{equation*}
q_{7}:\mathcal{L}_{7}\otimes \mathcal{L}_{7}\rightarrow O_{X_{7}}\left(
-K_{X_{7}}\right) \text{,}
\end{equation*}
to identify $\mathcal{L}_{7}^{\ast }$ with $\mathcal{L}_{7}\otimes O\left(
K_{X_{7}}\right) $. Using the adjunction formula $K_{X_{8}}=\pi ^{\ast
}K_{X_{7}}+L$, we obtain an isomorphism 
\begin{equation*}
\pi ^{\ast }\mathcal{L}_{7}^{\ast }\otimes O\left( L\right) \cong \pi ^{\ast
}\mathcal{L}_{7}\otimes O\left( K_{X_{8}}\right) .
\end{equation*}
Hence 
\begin{eqnarray*}
&&\pi ^{\ast }\mathcal{L}_{7}\otimes O\left( -L\right) +\pi ^{\ast }\mathcal{%
L}_{7}^{\ast }\otimes O\left( L\right) \\
&\cong &\pi ^{\ast }\mathcal{L}_{7}\otimes O\left( -L\right) \otimes \left(
O+O\left( L+K\right) \right) \\
&\cong &\pi ^{\ast }\mathcal{L}_{7}\otimes O\left( -L\right) \otimes \Lambda
_{1}\text{.}
\end{eqnarray*}

Combining these isomorphisms, we have 
\begin{equation*}
L\mathcal{E}_{8}=\pi ^{\ast }L\mathcal{E}_{7}+\mathcal{A}_{1}+\pi ^{\ast }%
\mathcal{L}_{7}\otimes O\left( -L\right) \otimes \Lambda _{1}\text{.}
\end{equation*}

Since $\mathcal{L}_{8}=L\mathcal{E}_{8}\otimes O\left( -K\right) $, a
decomposition of $\mathcal{L}_{8}$ is equivalent to a decomposition of $L%
\mathcal{E}_{8}$. Namely 
\begin{equation*}
\mathcal{L}_{8}=\pi ^{\ast }\mathcal{L}_{7}\otimes \Lambda _{1}^{\ast }+\pi
^{\ast }\mathcal{R}_{7}\otimes O\left( -L\right) +\mathcal{A}_{1}\otimes
O\left( -K\right) .
\end{equation*}

\subsection{Reduction to $D_{8}$-bundle}

Again, a decomposition of $\mathcal{L}_{8}$ is equivalent to a decomposition
of $L\mathcal{E}_{8}$ because $\mathcal{L}_{8}=L\mathcal{E}_{8}\otimes
O\left( -K\right) $,

Recall that there is a $D_{n-1}$-bundle over $X_{n}$ for any given choice of
ruling on $X_{n}$. On $X_{8}$, we have a $D_{8}$-bundle instead of a $D_{7}$%
-bundle over it. To construct this $D_{8}$-bundle $\mathcal{D}_{8}$ over $%
X_{8}$, we need to choose eight disjoint lines on $X_{8}$. We denote them as 
$L_{1},...,L_{8}$. Notice that such a choice is equivalent to represents $%
X_{8}$ as the blowup of $\mathbb{P}^{2}$ at eight distinct points. Let $H$
be the pullback of the hyperplane bundle on $\mathbb{P}^{2}$ to $X_{8}$. We
define the following rank sixteen vector bundle on $X_{8}$: 
\begin{equation*}
\mathcal{W}_{8}=\bigoplus_{i=1}^{8}\left( O\left( L_{i}\right) +O\left(
-L_{i}-K-H\right) \right) \text{.}
\end{equation*}
The isomorphism $O\left( L_{i}\right) \otimes O\left( -L_{i}-K-H\right)
\cong O\left( -K-H\right) $ defines a fiberwise quadratic form on $\mathcal{W%
}_{8}:$%
\begin{equation*}
\mathcal{W}_{8}\otimes \mathcal{W}_{8}\rightarrow O\left( -K-H\right) \text{.%
}
\end{equation*}
\qquad

Now we define the $D_{8}$-bundle $\mathcal{D}_{8}$ on $X_{8}$ as the
automorphism bundle of $\mathcal{W}_{8}$ which preserves the above fiberwise
quadratic form. We denote its associated Lie algebra bundle as $L\mathcal{D}%
_{8}$ which is itself a vector bundle of rank 120. Explicitly we have 
\begin{equation*}
L\mathcal{D}_{8}=\Lambda ^{2}\mathcal{W}_{8}\otimes O\left( K+H\right) \text{%
.}
\end{equation*}
It is not difficult to check the following equality: 
\begin{equation*}
L\mathcal{D}_{8}=O^{\oplus 8}+\bigoplus_{\substack{ D^{2}=-2  \\ DK=0  \\ DH%
\text{ even }}}O\left( D\right) \text{.}
\end{equation*}
In particular $L\mathcal{D}_{8}$ is a vector subbundle of 
\begin{equation*}
\mathcal{L}_{8}\otimes O\left( -K\right) =L\mathcal{E}_{8}=O^{\oplus
8}+\bigoplus_{\substack{ D^{2}=-2  \\ DK=0  \\ \text{ }}}O\left( D\right) .
\end{equation*}
In fact $L\mathcal{D}_{8}$ is a Lie sub-algebra bundle of $L\mathcal{E}_{8}$.

Next we introduce the positive spinor bundle for $L\mathcal{D}_{8}$. We
consider the vector bundle 
\begin{equation*}
\mathcal{S}^{+}=\bigoplus_{\substack{ l^{2}=-1  \\ lK=-1  \\ lH\text{ even } 
}}O\left( l\right) ,
\end{equation*}
which is a vector subbundle of $\mathcal{L}_{8}$. To see that $\mathcal{S}%
^{+}$ is a representation bundle of $L\mathcal{D}_{8}$ we need a
homomorphism 
\begin{equation*}
L\mathcal{D}_{8}\otimes \mathcal{S}^{+}\rightarrow \mathcal{S}^{+}\text{.}
\end{equation*}
If a divisor $D$ satisfies $D^{2}=-2,DK=0$ and a line $l$ with $l\cdot H$
even, then the divisor $l^{\prime }=D+l$ satisfies (i) $l^{\prime
2}=-3+2D\cdot l$, (ii) $l^{\prime }\cdot K=-1$ and (iii) $l^{\prime }\cdot H$
is even. Hence $l^{\prime }$ determines a line if and only if $D\cdot l=1$.
Using these various homomorphisms $O\left( D\right) \otimes O\left( l\right)
\rightarrow O\left( l^{\prime }\right) $, we obtain the above homomorphism $L%
\mathcal{D}_{8}\otimes \mathcal{S}^{+}\rightarrow \mathcal{S}^{+}$ which
makes $\mathcal{S}^{+}$ a representation bundle of $L\mathcal{D}_{8}$.

It is easy to see that the rank of $\mathcal{S}^{+}$ equals 128 and the
above homomorphism corresponds to a spinor representation of $D_{8}$. Now $L%
\mathcal{D}_{8}\otimes O\left( -K\right) $ and $\mathcal{S}^{+}$ are two
subbundle of $\mathcal{L}_{8}$ with trivial intersection. Moreover 
\begin{equation*}
rank\mathcal{L}_{8}=248=120+128=rank\left( L\mathcal{D}_{8}\otimes O\left(
-K\right) \right) +rank\mathcal{S}^{+}.
\end{equation*}
Therefore $\mathcal{L}_{8}$ is the direct sum bundle of $L\mathcal{D}%
_{8}\otimes O\left( -K\right) $ and $\mathcal{S}^{+}$, 
\begin{equation*}
\mathcal{L}_{8}=L\mathcal{D}_{8}\otimes O\left( -K\right) +\mathcal{S}^{+}%
\text{.}
\end{equation*}
Similarly we have 
\begin{equation*}
L\mathcal{E}_{8}=L\mathcal{D}_{8}+\mathcal{S}^{+}\otimes O\left( K\right) 
\text{.}
\end{equation*}
The above two decompositions are decompositions as representation bundle of $%
L\mathcal{D}_{8}$ corresponding to 
\begin{equation*}
\mathbf{L}_{8}=E_{8}=D_{8}+S^{+}\text{,}
\end{equation*}
for the Lie algebra $D_{8}=\mathbf{so}\left( 16\right) $.

\subsection{Construction of Lie algebra bundle $L\mathcal{E}_{n}$ revisited}

Now we present a different approach to describe the fiberwise Lie algebra
structure on $L\mathcal{E}_{n}$. Since $E_{8}$ is the biggest exceptional
Lie algebra, we will first use the decomposition $L\mathcal{E}_{8}=L\mathcal{%
D}_{8}+\mathcal{S}^{+}\otimes O\left( K\right) $ to describe the Lie algebra
structure on $L\mathcal{E}_{8}$. Then we can use this to describe $L\mathcal{%
E}_{n}$ with $n<8$.

Before we begin, we need two homomorphisms between representation bundles of 
$L\mathcal{D}_{8}$. The first one is a homomorphism 
\begin{equation*}
\mathcal{S}^{+}\otimes \mathcal{S}^{+}\rightarrow O\left( -2K\right) .
\end{equation*}
To describe this homomorphism, we notice that if $l$ is a line on $X_{8}$
with $l\cdot H$ even, then $\left( -l-2K\right) ^{2}=-1=\left( -l-2K\right)
\cdot K$. Namely there is a line $l^{\prime }$ in the class $-l-2K$.
Moreover $l^{\prime }\cdot H=-l\cdot H-2K\cdot H$ is also even. That is $%
O\left( l^{\prime }\right) $ is a subbundle of $\mathcal{S}^{+}$. We then
define the homomorphism $\mathcal{S}^{+}\otimes \mathcal{S}^{+}\rightarrow
O\left( -2K\right) $ using these various isomorphism $O\left( l\right)
\otimes O\left( l^{\prime }\right) \rightarrow O\left( -2K\right) $.

Notice that $l^{\prime }$ can also be characterized as the unique line on $%
X_{8}$ which intersects $l$ at three points. Therefore the above
homomorphism also gives a fiberwise non-degenerate quadratic form on $%
\mathcal{S}^{+}$ and gives an isomorphism between $\mathcal{S}^{+}$ and $%
\left( \mathcal{S}^{+}\right) ^{\ast }\otimes O\left( -2K\right) $.

The second homomorphism is an \textit{anti-symmetric} bilinear homomorphism: 
\begin{equation*}
\mathcal{S}^{+}\otimes \mathcal{S}^{+}\rightarrow L\mathcal{D}_{8}\otimes
O\left( -2K\right) \text{.}
\end{equation*}
If $l$ and $l^{\prime }$ are two lines on $X_{8}$ with $l\cdot H$ and $%
l^{\prime }\cdot H$ even. Then $D=l+l^{\prime }+2K$ satisfies (i) $%
D^{2}=-6+2l\cdot l^{\prime }$,(ii) $D\cdot K=0$ and (iii) $D\cdot H$ is
even. Therefore $O\left( D\right) $ is a subbundle of $L\mathcal{D}_{8}$ if
and only if $l\cdot l^{\prime }=2$. Using these various isomorphisms $%
O\left( l\right) \otimes O\left( l^{\prime }\right) \otimes O\left(
2K\right) \rightarrow O\left( D\right) $, we obtain the above homomorphism
up to sign. To avoid the sign problem, we use a different description. We
rewrite the spinor representation $L\mathcal{D}_{8}\otimes \mathcal{S}%
^{+}\rightarrow \mathcal{S}^{+}$ as the homomorphism $\left( \mathcal{S}%
^{+}\right) ^{\ast }\otimes \mathcal{S}^{+}\rightarrow L\mathcal{D}_{8}$
since $L\mathcal{D}_{8}$ is self-dual. Now we use the non-degenerate pairing 
$\mathcal{S}^{+}\otimes \mathcal{S}^{+}\rightarrow O\left( -2K\right) $ to
identify $\left( \mathcal{S}^{+}\right) ^{\ast }$ with $\mathcal{S}%
^{+}\otimes O\left( 2K\right) $, then we obtain $\mathcal{S}^{+}\otimes 
\mathcal{S}^{+}\rightarrow L\mathcal{D}_{8}\otimes O\left( -2K\right) $ the
anti-symmetric homomorphism which is invariant under $L\mathcal{D}_{8}$.

Now we can use the Lie algebra structure on $L\mathcal{D}_{8}$%
\begin{equation*}
\alpha :L\mathcal{D}_{8}\otimes L\mathcal{D}_{8}\rightarrow L\mathcal{D}_{8},
\end{equation*}
the spinor representation 
\begin{equation*}
\beta :L\mathcal{D}_{8}\otimes \left( \mathcal{S}^{+}\otimes O\left(
K\right) \right) \rightarrow \left( \mathcal{S}^{+}\otimes O\left( K\right)
\right) ,
\end{equation*}
and the above anti-symmetric homomorphism 
\begin{equation*}
\gamma :\left( \mathcal{S}^{+}\otimes O\left( K\right) \right) \otimes
\left( \mathcal{S}^{+}\otimes O\left( K\right) \right) \rightarrow L\mathcal{%
D}_{8},
\end{equation*}
to give an alternative way to define the Lie algebra structure on 
\begin{equation*}
L\mathcal{E}_{8}=L\mathcal{D}_{8}+\mathcal{S}^{+}\otimes O\left( K\right) .
\end{equation*}

Namely if $a+u$ and $b+v$ are two local section of $L\mathcal{E}_{8}$ with $%
a,b\in L\mathcal{D}_{8}$ and $u,v\in \mathcal{S}^{+}\otimes O\left( K\right) 
$, then we define their bracket as 
\begin{eqnarray*}
&&\left[ a+u,b+v\right] \\
&=&\alpha \left( a\otimes b\right) +\gamma \left( u\otimes v\right) +\beta
\left( a\otimes v-b\otimes u\right) .
\end{eqnarray*}
It is shown in chapter six of \cite{Adams} that fiberwise this algebraic
structure is isomorphic the Lie algebra structure of $E_{8}$.

To obtain $L\mathcal{E}_{n}$ for $n$ less than eight, all we need to do is
to choose $8-n$ disjoint lines $L_{j}$'s on $X_{8}$ and take away all those
line bundle summands of $L\mathcal{E}_{8}$ which are written as $O\left(
l-L_{j}\right) $ or $O\left( L_{j}-l\right) $ for $l$ and $L_{j}$ disjoint
lines, we also take away the summand spanned by homology classes of $L_{j}$%
's in $\Lambda \otimes _{\mathbb{Z}}O_{X}=O_{X}^{\oplus 8}$. Then such
bundle $L\mathcal{E}_{n}$ is pullback of a $E_{n}$ bundle from the blowdown
of $X_{8}$ along $L_{j}$'s.

There is another way to obtain $L\mathcal{E}_{7}$: We first choose two
disjoint lines $L_{1}$ and $L_{2}$ on $X_{8}$. We consider the $A_{1}$%
-bundle 
\begin{equation*}
\mathcal{A}_{1}=O\left( L_{1}-L_{2}\right) \oplus O\oplus O\left(
L_{2}-L_{1}\right)
\end{equation*}
on $X_{8}$ which is a Lie algebra subbundle of $L\mathcal{D}_{8}\subset L%
\mathcal{E}_{8}$. Then the centralizer bundle of $\mathcal{A}_{1}$ inside $L%
\mathcal{E}_{8}$ is a $E_{7}$-bundle on $X_{8}$ which in fact comes from
pullback of a $E_{7}$ bundle on the blow down of $X_{8}$ along the line in
the class $H-L_{1}-L_{2}$. However we do not know if a similar method would
produce other $E_{n}$ bundle with $n<7$.

\begin{center}
\textit{Acknowledgments:}
\end{center}

\textrm{The author would like to express his gratitude to Scot Adams, Jim
Bryan, Ron Donagi, Igor Dolgachev, Brendan Hassett, David Morrison, Levin
Norman, Xiaowei Wang and Siye Wu.}


\begin{thebibliography}{Hartshorne}
\bibitem[Adams]{Adams}  J.F. Adams, \textit{Lectures on Exceptional Lie
groups,} The University of Chicago Press, 1996.

\bibitem[BPV]{BPV}  W. Barth, C. Peters and A. Van de Ven, \textit{Compact
complex surfaces, }Springer-Verlag, 1984.

\bibitem[Beauville]{Beauville}  A. Beauville, \textit{Complex algebraic
surfaces, Asterisque }\textbf{54 }(1978) and Cambridge University Press,
1983.

\bibitem[Clemens]{Clemens}  H. Clemens, \textit{A scrapbook of complex curve
theory,} Plenum Press.

\bibitem[Dolgachev]{Dolgachev}  I. Dolgachev, \textit{Weyl groups and
Cremona transformations}, Proc. Symp. Pure Math. Vol. \textbf{40} (1983)
283-294.

\bibitem[Donagi]{Donagi}  R. Donagi, \textit{Principal bundles on elliptic
fibrations}. Asian J. Math.\textbf{\ 1} (1997), no. 2, 214--223.

\bibitem[FM]{FM}  R. Friedman, J. Morgan, \textit{Exceptional groups and del
Pezzo surfaces, }math.AG/0009155.

\bibitem[FMW]{FMW}  R. Friedman, J. Morgan, E. Witten, \textit{Vector
bundles and F- theory}. Comm. Math. Phys. \textbf{187} (1997), no. 3,
679--743.

\bibitem[FH]{FH}  W. Fulton, J. Harris, \textit{Representation theory},
Springer-Verlag, 1991.

\bibitem[Jacobson]{Jacobson}  N. Jacobson, \textit{Exceptional Lie algebras, 
}Marcel Dekker, 1971.

\bibitem[Hartshorne]{Hartshorne}  R. Hartshorne,\textit{\ Algebraic
geometry, }Springer-Verlag, 1977.

\bibitem[Looijenga]{Looijenga}  E. Looijenga, \textit{Rational surfaces with
an anti-canonical cycle}, Ann. of Math., \textbf{114} (1981) 267-322.

\bibitem[MPR]{MPR}  W.G. McKay, J. Patera, D.W. Rand, \textit{Tables of
representations of simple Lie algebras, Volume I. Exceptional simple Lie
algebra, }CRM, 1990.

\bibitem[Manin]{Manin}  Y. Manin, \textit{Cubic forms: algebra, geometry,
arithmetic,} Nauka, Moscow 1972; English translation, North-Holland,
Amsterdam London, 1974, second edition, 1986.

\bibitem[Reid1]{Reid1}  M. Reid, \textit{Nonnormal del Pezzo surfaces},
Publ. RIMS, Kyoto Univ. \textbf{30}, (1994), 695-727.

\bibitem[Reid2]{Reid2}  M. Reid, \textit{Chapters on algebraic surfaces}.
Complex algebraic geometry (Park City, UT, 1993), Amer. Math. Soc.,
Providence, RI, 1997.

\bibitem[Segre]{Segre}  B. Segre, \textit{The non-singular cubic surfaces,}
Oxford Press, 1942.

\bibitem[Slansky]{Slansky}  R. Slansky, \textit{Group theory for unified
model building, }Physics reports, \textbf{79}, No.1 (1981) p.1-128.
\end{thebibliography}
\end{document}